\documentclass[3p,final]{elsarticle}

\usepackage{amssymb}
\usepackage{amsmath}
\usepackage{tikz}
\usepackage{pgfplots} 
\usetikzlibrary{calc}
\newcommand{\Nx}{3}        
\newcommand{\Ny}{3}        

\pgfmathtruncatemacro{\m}{\Nx*\Ny}

\newif\ifShowGrid     \ShowGridtrue
\newif\ifShowDofs     \ShowDofsfalse   

\tikzset{
  element/.style={draw, line width=1pt},
  nodept/.style={circle, fill=black, inner sep=1pt},
  elabel/.style={font=\scriptsize},
  nlabel/.style={font=\scriptsize, anchor=south west},
  dofarrow/.style={->, blue!70!black, line width=1pt}
}

\usepackage{dsfont}
\usepackage{hyperref}
\usepackage{cleveref}
\usepackage{amsthm}
\usepackage{numprint}
\usepackage{babel}[english]
\usepackage{caption}
\usepackage{subcaption}
\usepackage{ gensymb }
\usepackage{lineno}
\modulolinenumbers[5]
\usepackage{ragged2e}
\usepackage{natbib}
\usepackage{graphicx}
\usepackage{multicol}
\usepackage{multirow}
\usepackage{relsize}
\usepackage{derivative}
\usepackage{booktabs}
\usepackage{bm}
\usepackage{soul}  
\DeclareMathOperator*{\argmin}{argmin} 

\usepackage{cellspace} 
\newtheorem{definition}{Definition}

\newtheorem{theorem}{Theorem}
\newcommand{\mc}{\textcolor[rgb]{0.00,0.00,0.00}} 

\begin{document}

\begin{frontmatter}



\title{Novel insights into Pareto fronts in multiobjective topology optimization and a comparative study of scalarization strategies}


\author[Diepenbeek,FlandersMake,TUIL]{Tom De Weer}
\ead{tom.deweer at kuleuven.be}
\author[DTU]{Ole Sigmund}
\author[TUIL]{Gabriele Eichfelder}

\affiliation[Diepenbeek]{organization={KU Leuven Campus Diepenbeek, Department of Mechanical Engineering},
            addressline={Wetenschapspark 27}, 
            city={Diepenbeek},
            postcode={3590}, 
            country={Belgium}}


\affiliation[FlandersMake]{organization={Flanders Make@KU Leuven},
            country={Belgium}}

\affiliation[DTU]{organization={Department of Civil and Mechanical Engineering, Technical University of Denmark},
            city={Kongens Lyngby},
            postcode={2800}, 
            country={Denmark}}

\affiliation[TUIL]{organization={Institute of Mathematics, Technische Universität Ilmenau},
            city={Ilmenau},
            postcode={D-98684}, 
            country={Germany}}

\begin{abstract}
\mc{Topology optimization seeks structures that minimize an objective function subject to constraints.\
While extensive research has focused on single-objective formulations, design is usually a trade-off between conflicting criteria.\
This naturally leads topology optimization to the multiobjective optimization field.\ 
Unfortunately, despite their long coexistence, interaction between these fields has remained limited.\ 
This work aims to bridge this gap.\ 
First, we provide a theoretical comparison of two common scalarization techniques, the weighted-sum and $\varepsilon$-constraint method, to the established Pascoletti-Serafini scalarization.\ 
Then, we perform extensive bi-objective numerical experiments, including minimization of volume, compliance, maximum stress and dynamic compliance.\
For each experiment, the four scalarization methods are aggregated, often alongside design space sampling and continuation, to better resolve the Pareto front.\
The major contribution of this work is the observation that the Pareto front consists of segments belonging to local fronts.\ 
Contrary to the smooth, convex fronts often reported in the literature, we show that this property can lead to discontinuous and nonconvex fronts.\
Finally, we numerically compare the approximation quality of the scalarization methods.\
The numerical examples demonstrate the robustness of the Pascoletti-Serafini scalarization against fronts with pronounced discontinuities, nonconvexities and high-curvature regions.\
In contrast, the weighted-sum and $\varepsilon$-constraint methods can produce clustering and gaps, depending on the shape of the front.\
}
\end{abstract}




\end{frontmatter}

\section{Introduction}

Topology optimization is concerned with the challenge of automatic design.\
Instead of a manual process by a designer or architect, the problem is posed as an optimization problem and solved with numerical tools, most notably simulation and optimization software.\
Since the seminal paper by Bendsoe and Kikuchi~\cite{bendsoe_kikuchi_1988}, topology optimization has distinguished itself from shape and size optimization by its ability to add and remove holes from prospective designs.\
This increased design freedom is beneficial for a wide range of applications.\
In addition to compliance minimization, the field has seen structural optimization considering plastic yielding~\cite{DuysinxBendsoe1992,daSilva2019}, buckling~\cite{Neves1995GeneralizedTD,Ferrari2019} and dynamic eigenvalues~\cite{DiazKikuchi1992}.\
Extensions to other physics, such as fluid dynamics~\cite{BorPet2003}, vibro-acoustics~\cite{cool_TO_VA} and electromagnetics~\cite{JenSig2011}, are widespread.\\

It is common practice to formulate topology optimization problems as the minimization of a single objective subject to one or more geometric or physical constraints.\
However, design typically involves balancing multiple, conflicting objectives.\
For example, bridges are designed as a trade-off between (i) weight, (ii) rigidity and (iii) stability against failure.\
A typical solution strategy is then to select weight as the primary objective and put constraints on rigidity and stability.\
For example, Ferrari and Sigmund~\cite{Ferrari2019} consider compliance minimization with a constraint on the volume and the buckling load factor.\
If a trade-off or Pareto curve is sought, then a common approach is to repeatedly solve this problem for varying constraint values.\
The process of turning a multiobjective problem into a series of single-objective problems is known as \emph{scalarization}, and the approach described so far is known as $\varepsilon$-constraint scalarization~\cite{Ehrgott2005}.\\

A second scalarization, called weighted-sum scalarization,  is to consider a weighted sum of the objectives and solve the resulting single-objective problem for various weight values.\
The literature review in this work illustrates that most Pareto curves in topology optimization are produced with either $\varepsilon$-constraint or weighted-sum scalarization.\
Although these methods have well-documented theoretical drawbacks~\cite{Eichfelder2021}, the topology optimization literature lacks a rigorous numerical comparison of weighted-sum, $\varepsilon$-constraint and more advanced scalarizations.\\

Topology optimization problems are known for their nonconvexity, which prohibits finding a global optimum without very expensive global optimization techniques, if at all possible.\
In contrast, many multiobjective optimization algorithms rely on global optimality.\
The effect of the strong nonconvexity \mc{on the Pareto front structure remains poorly understood}.\\

\mc{Together, the lack of a systematic comparison of scalarization methods and the limited understanding of the effect of nonconvexity severely complicate the development of tailored multiobjective topology optimization algorithms.\
Moreover, these challenges are self-reinforcing: poorly tailored algorithms provide limited insight, and limited insight hinders well-tailored algorithm development.\ }
To break this cycle, this paper provides three contributions.\

First, a theoretical comparison is made between (i) weighted-sum, (ii) $\varepsilon$-constraint and (iii) Pascoletti-Serafini~\cite{PascolettiSerafini1984} scalarization.\
The first two are chosen for their wide usage in topology optimization, despite their theoretical drawbacks.\
The latter is chosen for its simplicity and clear theoretical benefits.\

Second, a novel insight is provided into the structure of the Pareto front in topology optimization.\
We show that local optima in single-objective topology optimization translate to local Pareto fronts in a multiobjective context.\
Local fronts can both dominate and be dominated by each other, yielding a (global) Pareto front that can be nonconvex and even disconnected.\
This is first illustrated via a simple truss optimization example and then illustrated extensively in bi-objective topology optimization examples involving minimization of volume, compliance, maximum stress and dynamic compliance.\

\mc{Third, the work provides a qualitative comparison of the approximation quality obtained by each of the mentioned scalarization strategies, revealing their advantages and drawbacks.\
Together, the contributions allow for the development of tailored algorithms for multiobjective topology optimization. } \\

This paper is built up as follows.\
\Cref{sec:MOO,sec:TO} introduce multiobjective and topology optimization, respectively, to define the needed concepts.\
The intended audience is topology optimization experts, thus \Cref{sec:TO} is relatively brief, with many details moved to the appendix.\ 
\Cref{sec:MOTO} discusses multiobjective topology optimization, providing a more in-depth overview of the literature and the challenges induced by combining the two fields.\
Then, \Cref{sec:numerical_examples} provides numerical examples and \Cref{sec:Conclusions} summarizes the conclusions.\

\begin{table}[h!]
\caption{General notation.}
\label{tab:general_notation}
\begin{tabular}{l l}
\toprule
\textbf{Symbol} & \textbf{Name}  \\
\midrule
$i,j \in \mathbb{N}$ & Natural number, usually indices \\[0.5ex]
$\mathbb{R}_+$ & Non-negative real numbers \\[0.5ex]
$\mathbb{R}_{++}$ & Positive real numbers \\[0.5ex]
$\mathbb{R}_+^i$ & Non-negative real orthant of $\mathbb{R}^i$ \\[0.5ex]
$\mathbb{R}_{++}^i$ & Positive real orthant of $\mathbb{R}^i$ \\[0.5ex]
$A \cap B$ & Intersection of sets $A$ and $B$ \\[0.5ex]
\bottomrule
\end{tabular}
\end{table}

\begin{table}[h!]
\caption{Multiobjective notation, defined in \Cref{sec:MOO}.}
\label{tab:moo_notation}
\begin{tabular}{l l}
\toprule
\textbf{Symbol} & \textbf{Name}  \\
\midrule
$f_i : \mathbb{R}^n \rightarrow \mathbb{R}$ & Objective function component \\[0.5ex]
$\mathbf{f}  : \mathbb{R}^n \rightarrow \mathbb{R}^m$ & Objective function  \\[0.5ex]
$\mathbf{f}(X) \subseteq \mathbb{R}^m$ & Image set of $X$ under $\mathbf{f}$ \\[0.5ex]
$m \in \mathbb{N}$ & Number of objectives \\[0.5ex]
$n  \in \mathbb{N}$ & Number of optimization variables   \\[0.5ex]
$\mathbb{R}^n$ & Decision space, variable space, pre-image space \\[0.5ex]
$\mathbb{R}^m$ & Image space \\[0.5ex]
$\mathbf{x} \in \mathbb{R}^n$ & Optimization variable vector  \\[0.5ex]
$X \subseteq \mathbb{R}^n$ & Feasible set \\[0.5ex]
\bottomrule
\end{tabular}
\end{table}

\begin{table}[h!]
\caption{Scalarization methods notation, defined in \Cref{sec:MOO}.}
\label{tab:scalarization_notation}
\begin{tabular}{l l}
\toprule
\textbf{Symbol} & \textbf{Name}  \\
\midrule
$\mathbf{a}  \in \mathbb{R}^m$ & Reference point   \\[0.5ex] 
$\varepsilon_i \in \mathbb{R}$ & Constraint bound w.r.t objective $f_i$ \\[0.5ex]
$\mathbf{r}  \in  \mathbb{R}^m_{++}$ & Direction vector   \\[0.5ex]
$w_i \in \mathbb{R}_+$ & Weight w.r.t. objective $f_i$  \\[0.5ex]
$\mathbf{w} \in \mathbb{R}^m_+$ & Weights vector    \\[0.5ex]
$\mathbf{x}^{\min,i} \in \mathbb{R}^n$ & Minimizer for $f_i$    \\[0.5ex]
\bottomrule
\end{tabular}
\end{table}

\begin{table}[h!]
\caption{Topology optimization notation, defined in \Cref{sec:TO}.}
\label{tab:TO_notation}
\begin{tabular}{l l}
\toprule
\textbf{Symbol} & \textbf{Name}  \\
\midrule
$b,e,d$ & Used as subscript to refer to the blueprint ($b$), eroded ($e$) and dilated ($d$) designs  \\[0.5ex]
$f : \mathbb{R}^n \rightarrow \mathbb{R}$ & Single objective function  \\[0.5ex]
$\mathbf{g} : \mathbb{R}^n \rightarrow \mathbb{R}^k$ & Constraint function  \\[0.5ex]
$\mathbf{h} : \mathbb{R}^n \rightarrow \mathbb{R}^n$ & Filter function  \\[0.5ex]
$k \in \mathbb{N}$ & Number of constraints \\[0.5ex]
$l_x, l_y \in \mathbb{R}_{++}$ & Length of design domain in $x$ and $y$ direction \\[0.5ex]
$n_x, n_y \in \mathbb{N}$ & Number of elements in $x$ and $y$ direction \\[0.5ex]
$N_i$ & Node $i$ \\[0.5ex]
$\Omega \subseteq \mathbb{R}^2$ & Design domain  \\[0.5ex]
$p : \mathbb{R} \rightarrow \mathbb{R}$ & Material interpolation law or penalization function  \\[0.5ex]
$P \in \mathbb{R}$ & SIMP penalization constant  \\[0.5ex]
$q \in \mathbb{N}$ & Number of degrees of freedom \\[0.5ex]
$\Delta q \in \mathbb{N}$ & Number of fixed degrees of freedom \\[0.5ex]
$u_{x,i}, u_{y,i} \in \mathbb{R}$ & $x$ and $y$ degree of freedom of node $N_i$    \\[0.5ex]
$\mathbf{u} \in \mathbb{R}^m$ & Displacement vector    \\[0.5ex]
$\mathbf{x}_P \in \mathbb{R}^n$ & (Physical) design vector  \\[0.5ex]
\bottomrule
\end{tabular}
\end{table}

\paragraph{Notation} \Cref{sec:MOO,sec:TO} introduce multiobjective and topology optimization.\
For historical reasons, the two fields have different notation and terminology.\
To prevent confusion, this work uses the consistent notation listed in \Cref{tab:general_notation,tab:moo_notation,tab:scalarization_notation,tab:TO_notation}.\
For most symbols, the following conventions are followed.\
\begin{itemize}
\item Scalars and functions that map to scalars are denoted in lowercase (e.g., $a$).\
\item Vectors and functions that map to vectors are denoted in lowercase and bold (e.g., $\mathbf{a}$).\
\item Matrices and functions that map to matrices are denoted in uppercase and bold (e.g., $\mathbf{A}$).\
\end{itemize}
Symbols that appear only in the appendix are not listed in the tables.

\section{Multiobjective optimization} \label{sec:MOO}

This section introduces basic concepts and theorems from multiobjective optimization.\
For more information, the interested reader is referred to standard reference works such as \cite{Miettinen1998,Ehrgott2005,Eichfelder2008}.\
Additionally, one can search for \emph{vector optimization}, \emph{multicriteria optimization} and \emph{decision making}, as these keywords are related to multiobjective optimization.\\

Consider the multiobjective optimization problem
\begin{equation}
  \min_{\mathbf{x} \in X} \, \mathbf{f}(\mathbf{x}) = (f_1(\mathbf{x}),\dots,f_m(\mathbf{x}))^\top,
  \label{eq:mop}
  \tag{MOP}
\end{equation}
where the feasible set $X \subseteq \mathbb{R}^n$ encodes the constraints and $f_i : \mathbb{R}^n \to \mathbb{R}$ are objective
functions.\
This study considers a topology optimization context, where objectives and constraints are usually continuous and differentiable with respect to the optimization variables $\mathbf{x}$.\
For $m=1$, \ref{eq:mop} reduces to a single-objective optimization problem.\
Problems with $m=2$ and $m=3$ objectives are denoted as bi- and tri-objective problems, respectively.\
The feasible set $X$ is nonempty but for now unspecified.\
We define 

\begin{equation}
\mathbf{f}(X) := \left\{ \mathbf{f}(\mathbf{x})  \, | \, \mathbf{x} \in X  \right\}
\end{equation}

to be the \emph{image set}.\

\subsection{Optimality}

Generally, there is no single solution $\mathbf{x}$ that minimizes all objectives simultaneously.\
A solution is considered \emph{(weakly) efficient} via the concept of Edgeworth-Pareto efficiency~\cite{Edgeworth1881,Pareto1898}, defined as follows.

\begin{figure}
    \centering
    \includegraphics[width=\linewidth]{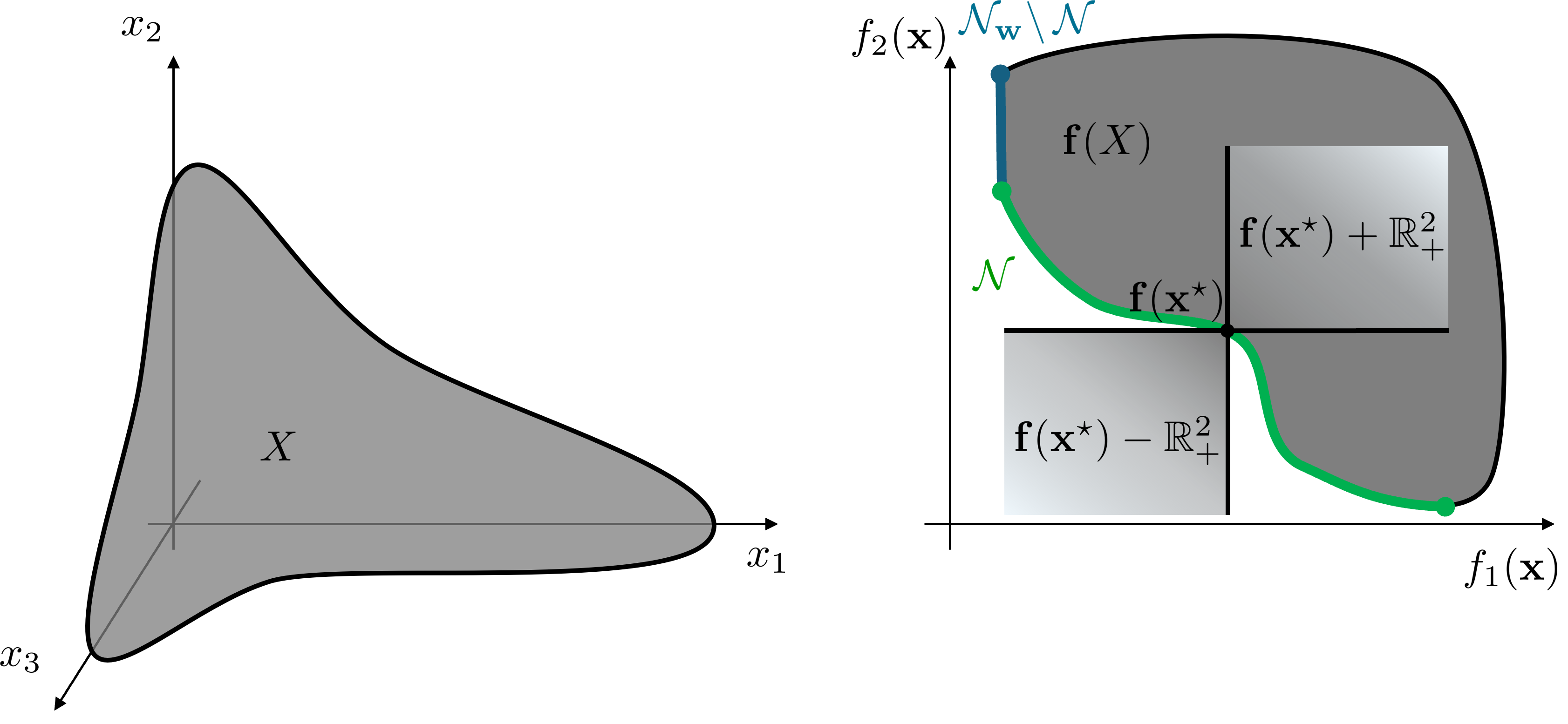}
    \caption{Concepts in multiobjective optimization, illustrated for three-dimensional ($n=3$) pre-image space (left) and a bi-objective ($m=2$) image set (right).}
    \label{fig:MOO_Illustration}
\end{figure}

\begin{definition} \label{def:efficient}
A feasible point $\mathbf{x}^\star \in X$ is \emph{efficient} if there is no $\mathbf{x} \in X$ with $f_i(\mathbf{x}) \leq f_i(\mathbf{x}^\star)$ for all $i \in \{ 1, \ldots, m\}$ and $f_j(\mathbf{x}) < f_j(\mathbf{x}^\star)$ for at least one $j \in \{ 1,\ldots, m\}$.
\end{definition}
\begin{definition} \label{def:weakly_efficient}
A feasible point $\mathbf{x}^\star \in X$ is \emph{weakly efficient} if there is no $\mathbf{x} \in X$ with $f_i(\mathbf{x}) < f_i(\mathbf{x}^\star)$ for all $i \in \{ 1, \ldots, m\}$.
\end{definition}
Definition \ref{def:efficient} reflects that an efficient solution cannot be improved for one objective function value without worsening another.\
Definition \ref{def:weakly_efficient} reflects that a weakly efficient solution cannot be improved for all objectives simultaneously.\
The terms (weak) (Edgeworth)-Pareto minimality/optimality are used as synonyms in the literature, but for clarity not further used in this paper.\

When a solution $\mathbf{x}^\star$ is (weakly) efficient, its image $\mathbf{f}(\mathbf{x}^\star)$ is called (weakly) \emph{nondominated}.\
A more geometric definition of a nondominated point is as follows.

\begin{definition} \label{def:nondominated}
A point $\mathbf{f}(\mathbf{x}^\star) \in \mathbf{f}(X)$ is nondominated if $\left( \mathbf{f}(\mathbf{x}^\star) - \mathbb{R}_{+}^m \right) \cap \mathbf{f}(X) = \mathbf{f}(\mathbf{x}^\star)$.
\end{definition}
\begin{definition}
A point $\mathbf{f}(\mathbf{x}^\star) \in \mathbf{f}(X)$ is weakly nondominated if $\left( \mathbf{f}(\mathbf{x}^\star) - \mathbb{R}_{++}^m \right) \cap \mathbf{f}(X) = \emptyset$.
\end{definition}

Here, the notation $\mathbf{u} - \mathbb{R}_{+}^m$ and $\mathbf{u} - \mathbb{R}_{++}^m$ denotes the set of points $\mathbf{v} \in \mathbb{R}^m$ such that $\mathbf{u} - \mathbf{v}$ is in $\mathbb{R}_{+}^m$ and $\mathbb{R}_{++}^m$, respectively.\
The set of all efficient solutions is denoted the efficient set $\mathcal{E}(\mathbf{f}, X) \subseteq \mathbb{R}^n$ and its image
\begin{equation}
\mathcal{N} := \mathbf{f}(\mathcal{E}(\mathbf{f}, X) ) = \{ \mathbf{f}(\mathbf{x}) \, | \, \mathbf{x} \in  \mathcal{E}(\mathbf{f}, X) \} 
\end{equation}
is called the \emph{nondominated set}.\
In engineering contexts, this is more often called the \emph{Pareto frontier/front/curve}, although technically, the nondominated set is not necessarily a curve.\
Similarly, the set of all weakly efficient solution is denoted as the weakly efficient set $\mathcal{E}_w(\mathbf{f}, X) \subseteq \mathbb{R}^n$ and its image $\mathcal{N}_w:= \mathbf{f}(\mathcal{E}_w(\mathbf{f}, X) )$ the weakly nondominated set.\
Weakly nondominated solutions that are not nondominated, i.e., the pre-image of the set $\mathcal{N}_w \backslash \mathcal{N}$, are (for numerical reasons) not found in any of the examples in this work and hence of lesser importance.\
The notion is only introduced to support the theorems on scalarization methods (\Cref{sec:scalarization}).\\

\Cref{fig:MOO_Illustration} illustrates the concepts introduced in this section for a bi-objective example ($m=2$) with three optimization variables ($n=3$).\
Here, the feasible set $X \subseteq \mathbb{R}^3$ is mapped by the objective function to the image set $\mathbf{f}(X) \subseteq \mathbb{R}^2$.\
Definition \ref{def:nondominated} is illustrated as well: the point $\mathbf{f}(\mathbf{x}^\star)$ is nondominated since there are no feasible points in the set $\mathbf{f}(\mathbf{x}^\star)-\mathbb{R}_{+}^2$ except $\mathbf{f}(\mathbf{x}^\star)$ itself.\
It is therefore part of the nondominated set $\mathcal{N}$, shown in green.\
Additionally, the blue line illustrates $\mathcal{N}_w \backslash \mathcal{N}$. \\

To support the theorems in the next section, we recall the definition of a convex set.\
\begin{definition}\label{def:convex_set}
A set $Y \subseteq \mathbb{R}^m$ is convex if for all $\mathbf{u}, \mathbf{v}$ in $Y$, it holds that $t\mathbf{u} + (1-t)\mathbf{v} \in Y$, for any $t \in [0,1]$.
\end{definition}
We then define convexity for the nondominated set as follows.\
\begin{definition}\label{def:convex_pf}
The nondominated set $\mathcal{N}$ to \ref{eq:mop} is convex if $\mathcal{N}+\mathbb{R}_+^m = \left\{  \mathbf{u} + \mathbf{v} \, | \, \mathbf{u} \in \mathcal{N}, \mathbf{v} \in \mathbb{R}_+^m \right\}$ is a convex set.
\end{definition}
Under this definition, the nondominated set of \Cref{fig:MOO_Illustration} is nonconvex since it bends inwards around the illustrated point $\mathbf{f}(\mathbf{x}^\star)$.

\subsection{Scalarization methods} \label{sec:scalarization}

\begin{figure}
    \centering
    \includegraphics[width=\linewidth]{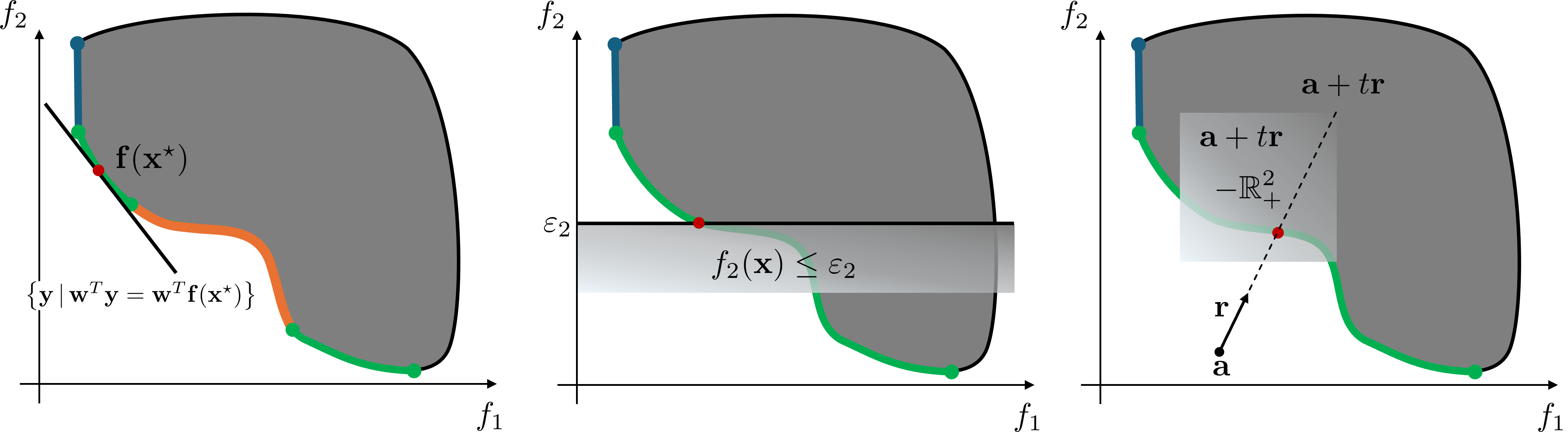}
    \caption{Illustration of weighted-sum method (left), $\varepsilon$-constraint method (middle) and Pascoletti-Serafini scalarization (right). Blue and green denote (weakly) nondominated solutions. Orange denotes nondominated solutions that cannot be found by weighted-sum scalarization as they lie in the concave region of the nondominated set.}
    \label{fig:ScalarizationMethods}
\end{figure}

Scalarization methods transform \ref{eq:mop} into a single-objective optimization problem that depends on one or more scalarization parameters.\
By varying these parameters, different efficient points are obtained.\
Scalarization thus replaces a multiobjective problem with a family of single-objective problems.\
A crucial property is if $\mathcal{E}$ can be completely recovered by varying the parameters.\
To this end, we introduce the recovery properties of the three scalarization methods studied in this work, namely weighted-sum scalarization~\cite{Zadeh1963}, $\varepsilon$-constraint scalarization~\cite{Haimes1971,Marglin1969} and the scalarization proposed by Pascoletti and Serafini~\cite{PascolettiSerafini1984}.\
The former two are selected due to their wide usage in topology optimization, as further discussed in \Cref{sec:MOTO}.
The latter is chosen because of its attractive theoretical properties.\
All three are illustrated in \Cref{fig:ScalarizationMethods}.

\subsubsection{Weighted-sum scalarization}

The weighted-sum scalarization is defined as follows.

\begin{definition}
Given a weight vector $\mathbf{w} \in \mathbb{R}^m_+$ such that $\sum_{i=1}^m w_i = 1$, the weighted-sum
scalarization to \ref{eq:mop} is
\begin{equation} \label{eq:WS_formulation}
  \min_{\mathbf{x} \in X} \, \mathbf{w}^\top \mathbf{f}(\mathbf{x})=\sum_{i=1}^m w_i f_i(\mathbf{x}).
  \tag{$\mathrm{WS}(\mathbf{w})$}
\end{equation}
\end{definition}

The weighted-sum scalarization has the following properties (see \cite[Prop. 3.9\&3.10]{Ehrgott2005}).

\begin{theorem}
If $\mathbf{x}^\star$ is a minimizer of \ref{eq:WS_formulation}, then $\mathbf{x}^\star$ is a weakly efficient solution of \ref{eq:mop}.\
Furthermore, if $w_i > 0$ for all $i \in \left\{ 1,\ldots, m \right\}$ then $\mathbf{x}^\star$ is an efficient solution of \ref{eq:mop}.
\end{theorem}

However, the inverse is not always true.

\begin{theorem}
If $\mathbf{f}(X)+\mathbb{R}_+^m$ is convex according to \Cref{def:convex_set} and $\mathbf{x}^\star$ is a weakly efficient solution of \ref{eq:mop}, then there exists a weight vector $\mathbf{w} \in \mathbb{R}^m_+$ with $\sum_{i=1}^m w_i=1$ such that $\mathbf{x}^\star$ is a minimizer of \ref{eq:WS_formulation}.
\end{theorem}

A sufficient condition for convexity of $\mathbf{f}(X)+\mathbb{R}_+^m$ is when both $X$ and all $f_i$ are convex.\
Unfortunately, this is not the case in many applications, including topology optimization.\
The result is that portions of the nondominated set cannot be recovered by any combination of the weights.\
The example in \Cref{fig:ScalarizationMethods} illustrates this issue: the orange subset of $\mathcal{N}$ represents nondominated points that are not images of minimizers of \ref{eq:WS_formulation}.\

\subsubsection{Epsilon-constraint scalarization}

The idea behind the $\varepsilon$-constraint scalarization is to obtain efficient solutions by solving single-objective problems with varying bounds on the constraints.\ 

\begin{definition}
For an objective $f_i$, $i \in \{ 1,\ldots, m\}$ and constraint bounds $\varepsilon_j$, $j \in J := \{ 1, \ldots, m \} \backslash \{ i \}$, the $\varepsilon$-constraint scalarization to \ref{eq:mop} is 
\begin{equation} \label{eq:EC_formulation}
  \min_{x \in X} f_i(\mathbf{x}) 
  \quad \text{s.t.} \quad f_j(\mathbf{x}) \le \varepsilon_j \quad \forall j \in J.
    \tag{$\mathrm{EC}(\bm{\varepsilon})$}	
\end{equation}
\end{definition}

For conciseness, we collect the constraint bounds $\varepsilon_j$ into $\bm{\varepsilon} \in \mathbb{R}^{m-1}$.\
The following results hold (see \cite[Prop. 4.4\& Theorem 4.5]{Ehrgott2005}).

\begin{theorem}
If $\mathbf{x}^\star$ is a minimizer of \ref{eq:EC_formulation}, then $\mathbf{x}^\star$ is a weakly efficient solution of \ref{eq:mop}.\
\end{theorem}
\begin{theorem}
Let $i \in \left\{ 1, \ldots, m \right\}$.\
 If $\mathbf{x}^\star$ is an efficient solution of \ref{eq:mop}, then there exists a $\bm{\varepsilon}$ such that $\mathbf{x}^\star$ is a minimizer of \ref{eq:EC_formulation}.\
Specifically, this holds for the choice $\varepsilon_j=f_j(\mathbf{x}^\star)$.
\end{theorem}

A drawback of the $\varepsilon$-constraint scalarization is that it cannot find all weakly nondominated points (see \cite[Cor.2.28]{Eichfelder2008}).\
However, as discussed, such points are not of practical importance in topology optimization.\
A more serious drawback is that \ref{eq:EC_formulation} does not always have a minimizer.\
That is, for some $\varepsilon_j$, the problem can become infeasible: no $\mathbf{x}^\star$ exists that satisfies both $\mathbf{x}^\star \in X$ and $f_j(\mathbf{x}^{\star}) \le \varepsilon_j$ for all $j \in J$.\
\mc{For example, in \Cref{fig:ScalarizationMethods}, if $\varepsilon_2 < \min \{ f_2(\mathbf{x}) \, | \, \mathbf{x} \in X \}$ then there is no feasible solution.\
Additionally, if $\varepsilon_2 > \max \{ f_2(\mathbf{x}) \, | \, \mathbf{x} \in X \}$ then the $\varepsilon_2$ constraint is not active.\
As the numerical examples in \Cref{sec:numerical_examples} will illustrate, these issues also persist in topology optimization contexts, complicating the choice of appropriate values for $\varepsilon_j$.\ }

\subsubsection{Pascoletti-Serafini scalarization}

The scalarization of Pascoletti and Serafini~\cite{PascolettiSerafini1984} employs a slack variable $t \in \mathbb{R}$ and directly incorporates Definition \ref{def:nondominated}.

\begin{definition}
Given parameter vectors $\mathbf{a} \in \mathbb{R}^m$ and $\mathbf{r} \in \mathbb{R}^m_{++}$, (a variant of) the scalarization to \ref{eq:mop} proposed by Pascoletti and Serafini is
\begin{equation} \label{eq:PS_formulation}
  \min_{x \in X, \, t \in \mathbb{R}}
    \ t
  \quad \text{s.t.} \quad
    \mathbf{f}(\mathbf{x}) \in \mathbf{a} + t\mathbf{r} - \mathbb{R}_{+}^m.
\tag{$\mathrm{PS}(\mathbf{a},\mathbf{r})$}
\end{equation}
\end{definition}

Note that the vector constraint $\mathbf{f}(\mathbf{x}) \in \mathbf{a} + t\mathbf{r} - \mathbb{R}_{+}^m$ yields $m$ scalar constraints $a_i + tr_i -f_i(\mathbf{x})  \geq 0$ for all $i \in \{ 1, \ldots, m \}$.\
In this work, we vary $\mathbf{a}$ and fix $\mathbf{r}$ in a manner introduced in \Cref{sec:numerical_examples}.\
Note that if instead of $\mathbf{r} \in \mathbb{R}^m_{++}$, the direction $\mathbf{r}$ is a unit vector in $ \mathbb{R}^m_{+}$ then the $\varepsilon$-constraint scalarization is a special case of \ref{eq:PS_formulation}~\cite{Eichfelder2009}.\
That is, \ref{eq:PS_formulation} reduces to \ref{eq:EC_formulation} if $\mathbf{r}$ lies along a unit vector.\\

The Pascoletti-Serafini scalarization has all the important properties a scalarization method should have (see e.g. \cite[Theorem 2.1\&2.11]{Eichfelder2008}).\

\begin{theorem}
If $(t^\star,\mathbf{x}^\star)$ is a minimizer of \ref{eq:PS_formulation}, then $\mathbf{x}^\star$ is a weakly efficient solution of \ref{eq:mop}.
\end{theorem}

\begin{theorem}
If $\mathbf{x}^\star$ is a weakly efficient solution of \ref{eq:mop}, then there is some $t^\star \in \mathbb{R}$ and $\mathbf{a} \in H$ such that $(t^\star, \mathbf{x}^\star)$ is a minimizer of \ref{eq:PS_formulation} for any $\mathbf{r} \in \mathbb{R}^m_{++}$.\
\end{theorem}

From these properties it follows that solving \ref{eq:PS_formulation} yields a weakly efficient solution of \ref{eq:mop} and that all such solutions can be found by fixing $\mathbf{r} \in \mathbb{R}^m_{++}$ and varying $\mathbf{a}$.\
Furthermore, from a numerical point of view, the Pascoletti-Serafini scalarization is particularly powerful since even a ``bad choice'' for $\mathbf{a}$ yields a weakly efficient solution, avoiding the infeasibility problem that occurs for $\varepsilon$-constraint scalarization.\

\Cref{tab:scalarizations_overview} summarizes this section.

\begin{table}[h!]
\caption{Summary of key properties of the weighted-sum, $\varepsilon$-constraint and Pascoletti-Serafini scalarizations.}
\label{tab:scalarizations_overview}
\centering
\begin{tabular}{m{0.2\linewidth} m{0.2\linewidth} m{0.2\linewidth} m{0.4\linewidth}}
\toprule
\textbf{Scalarization} & \textbf{Parameters} & \textbf{Minimizing yields...} & \textbf{Varying parameters yields...} \\[0.5ex]
 \midrule \\
Weighted-sum & $\mathbf{w} \in \mathbb{R}^m_+$ & weakly efficient $\mathbf{x}$ & $\mathcal{N}_w$ if $\mathbf{f}(X)+\mathbb{R}^m_+$ is convex  \\[0.5ex]
$\varepsilon$-constraints & $i \in \mathbb{N}$ and $\varepsilon_j \in \mathbb{R}$ for all $j \in \left\{ 1,\ldots,m\right\} \backslash i$ & weakly efficient $\mathbf{x}$ & $\mathcal{N}$, but can yield infeasible problems  \\[0.5ex]
Pascoletti-Serafini & $(\mathbf{a},\mathbf{r}) \in \mathbb{R}^m \times \mathbb{R}^m_{++}$ & weakly efficient $\mathbf{x}$ & $\mathcal{N}_w$ \\[0.5ex]
\bottomrule
\end{tabular}
\end{table}

\section{Topology optimization}\label{sec:TO}

This section discusses the main features of the density-based topology optimization procedure used in this work, with numerical details and extensions moved to the appendices.\
Note that the results of this paper are independent of the chosen design representation: the same conclusions should hold if, e.g., a level set method were used.\

\subsection{Optimization variables}

Consider a rectangular domain $\Omega \subseteq \mathbb{R}^2$ of size  $l_x$ by $l_y$, with $l_x,l_y \in \mathbb{R}_{++}$, meshed with $n=n_{x} \times n_{y}$ square elements.\ 
The goal of topology optimization is to choose the physical design vector $\mathbf{x}_P \in \{0,1\}^n$, where an entry $x_{P,i} \in \{0,1\}$ corresponds to the density of element $i$.\
Elements with $x_{P,i}=1$ denote the presence of solid material and, conversely, elements with $x_{P,i}=0$ denote void, i.e., the absence of material.\\

\begin{figure}[h]
   \centering
   \includegraphics[width=\linewidth]{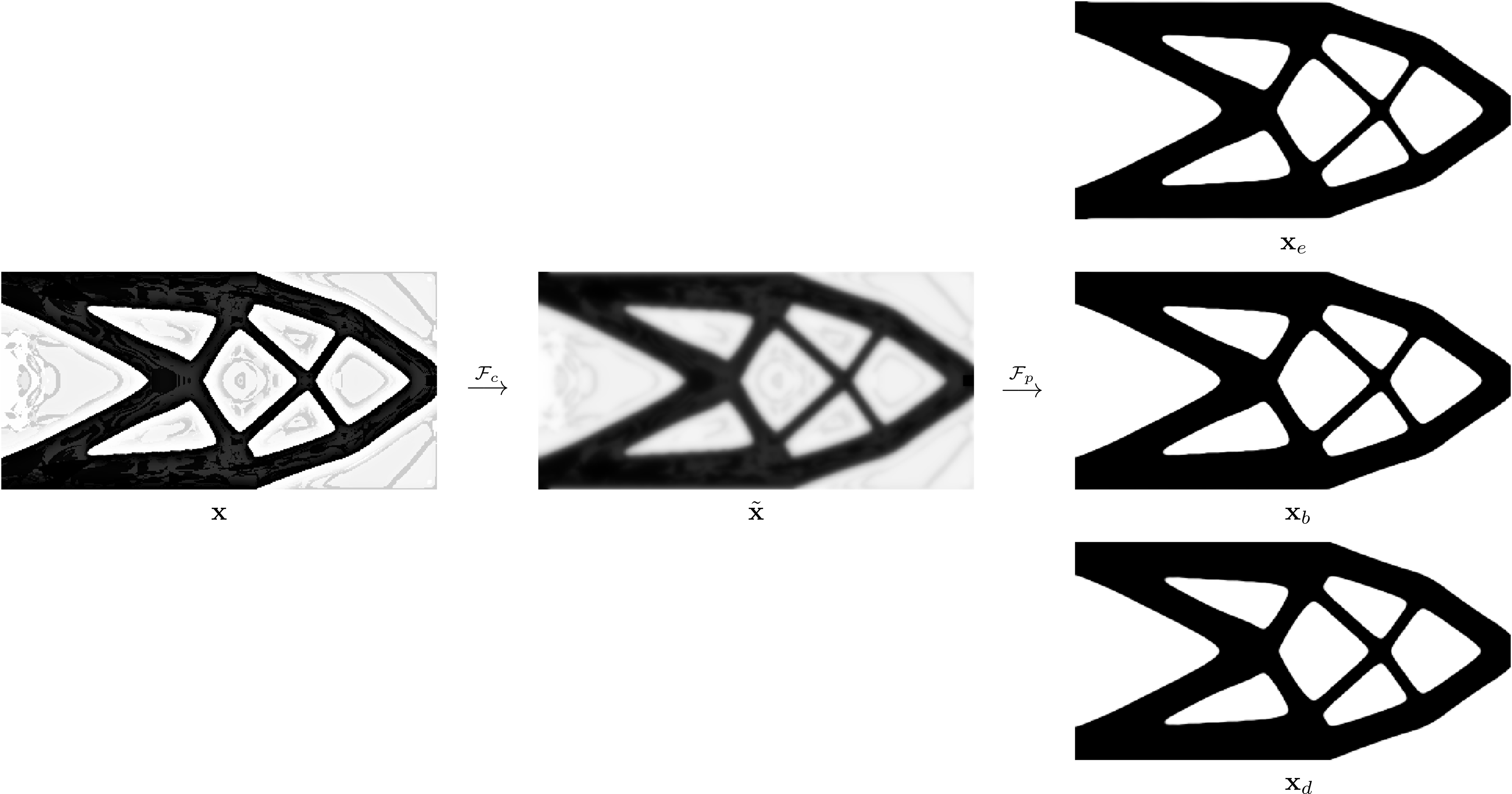}
   \caption{Filtering approach used in this work, producing the blueprint $\mathbf{x}_b$, eroded $\mathbf{x}_e$ and dilated $\mathbf{x}_d$ designs.}
   \label{fig:Filtering}
\end{figure}

To avoid solving an integer problem and allow the use of efficient gradient-based algorithms, the problem is relaxed to allow for continuous design variables: $\mathbf{x}_P \in [0,1]^n$.\
To prevent checkerboard patterns and enforce a length scale, we use the filtering scheme illustrated in \Cref{fig:Filtering}.\
The scheme employs the PDE filter from Lazarov and Sigmund~\cite{Lazarov2011} with the consistent boundary conditions from Wallin et al.~\cite{Wallin2020} and the smoothed Heaviside projection~\cite{Guest2004}, combined through the robust approach from Wang et al.~\cite{wang2011projection}.\
Because of the robust approach, the Heaviside operator is performed at three different levels and hence produces three physical design vectors: the blueprint design $\mathbf{x}_b$, the eroded design $\mathbf{x}_e$ and the dilated design $\mathbf{x}_d$.\
The value of the objective and constraint functions are then defined by taking the maximum of the three designs.\ 
That is,
\begin{equation} \label{eq:robustness}
\begin{aligned}
f(\mathbf{x}_P(\mathbf{x})) := \max (f(\mathbf{x}_b(\mathbf{x})), f(\mathbf{x}_e(\mathbf{x})), f(\mathbf{x}_d(\mathbf{x}))), \\
\mathbf{g}(\mathbf{x}_P(\mathbf{x})) := \max (\mathbf{g}(\mathbf{x}_b(\mathbf{x})), \mathbf{g}(\mathbf{x}_e(\mathbf{x})), \mathbf{g}(\mathbf{x}_d(\mathbf{x}))), \\
\end{aligned}
\end{equation}
where the $\max$-operator is applied elementwise and $f$ and $\mathbf{g}$ are generic objective and constraint functions.\
A generic single-objective density-based topology optimization problem is thus formulated as

\begin{equation} \label{eq:density_based_topology_optimization}
    \min_{\mathbf{x} \in [0,1]^n} f(\mathbf{x}_P) \quad \text{s.t.} \quad \mathbf{g}(\mathbf{x}_P) \leq 0,
\end{equation}

where $f: \mathbb{R}^n \rightarrow \mathbb{R}$ is the objective function and $\mathbf{g}: \mathbb{R}^n \rightarrow \mathbb{R}^k$ is a $k$-dimensional constraint.\
\mc{Both objective and constraint function values depend on $\mathbf{x}_P$, which is either the blueprint, dilated or eroded design representation as determined by the $\max$ operator in \Cref{eq:robustness}.} \\

For the volume and compliance objectives, discussed below, the $\max$-operator reduces to the dilated volume and the eroded compliance.\
However, for the more challenging objectives studied in this work, namely dynamic compliance and maximum stress, the nondifferentiable $\max$ function is dealt with via a bound formulation and additional optimization variables.\
This is discussed in \ref{app:SlackVariableTreatment}.\
The compliance and volume objective functions are now briefly introduced.

\subsection{Volume}

The volume $V: \mathbb{R}^n \rightarrow \mathbb{R}$ is found by summing the density values,

\begin{equation}
V(\mathbf{x}) = \sum_{i=1}^{n} x_i.
\end{equation}

The volume of the design domain, $V_{\max} \in \mathbb{R}$, is used to define the volume fraction $v: \mathbb{R}^n \rightarrow \mathbb{R}$

\begin{equation}
v(\mathbf{x}) =\frac{V(\mathbf{x})}{V_{\max}}.
\end{equation}

Thus, the constraint $v(\mathbf{x}) \leq 0.5$ means that at most $50\%$ of the design domain can contain solid material.\
The robust filtering approach means the volume fraction is the maximum of the volume of the eroded, dilated and blueprint design.\
However, due to filter's structure, the dilated design always has the largest volume~\cite{Lazarov2016}.\
Thus, $\max (v(\mathbf{x}_b(\mathbf{x})), v(\mathbf{x}_e(\mathbf{x})), v(\mathbf{x}_d(\mathbf{x}))) = v(\mathbf{x}_d(\mathbf{x})) :=v_d \in \mathbb{R}_+$.

\subsection{Compliance}

Compliance is the oldest objective function in structural optimization, going back to the work of Michell~\cite{Michell1904}.\
It equals the work of the external forces and hence represents the inverse of the stiffness of the design.\
Compliance is here computed numerically via the finite element method.\
Consider again the rectangular finite element grid with $n=n_x \times n_y$ elements and thus $(n_x+1) \times (n_y+1)$ nodes.\
Each node $N_i$ has two associated degrees of freedom, its $x$ translation $u_{i,x}$ and $y$ translation $u_{i,y}$, for a total of $q=2(n_x+1) \times (n_y+1)$ degrees of freedom.
These are collected in the displacement vector $\mathbf{u} \in \mathbb{R}^{q}$.\
We also introduce the external force vector $\mathbf{f}_{\mathrm{ext}} \in \mathbb{R}^q$.\
Finally, the stiffness matrix $\mathbf{K}(\mathbf{x}) \in \mathbb{R}^{q \times q}$ is found by evaluating a map $\mathbb{R}^n \rightarrow \mathbb{R}^{q \times q}$ that sums the $n$ elemental stiffness matrices $\mathbf{K}_i \in \mathbb{R}^{q \times q}$:
\begin{equation}
\mathbf{K}(\mathbf{x}) = \sum_{i=1}^{n}  p(x_i) \mathbf{K}_i(x_i).
\end{equation}
Here, $p: \mathbb{R} \rightarrow \mathbb{R}$ is the Solid Isotropic Material with Penalization (SIMP) interpolation $p := x \rightarrow \delta + (1-\delta)x^P$~\cite{Bendsoe1989}.\
Here,  $\delta$ is a small number to ensure $p(0)=\delta \neq 0$ and prevent ill-conditioning of $\mathbf{K}(\mathbf{x})$.
The penalization $P \in \mathbb{R}_+$ is increased during optimization to drive the design towards black-white solutions.\
Numerical details are in \ref{app:ComplianceMinimization}.\\

Using the above, we now write the state equation:

\begin{equation} \label{eq:Ku=f}
\mathbf{K}(\mathbf{x}) \mathbf{u}=\mathbf{f}_{\mathrm{ext}}.
\end{equation}
The state equation enforces linear elasticity and is enforced implicitly.\
That is, after the filtering operation provides a design vector $\mathbf{x}$, we construct $\mathbf{K}(\mathbf{x})$ and solve \Cref{eq:Ku=f}.\
As a result, the state vector $\mathbf{u}$ depends implicitly on the design vector $\mathbf{x}$.\
To incorporate the Dirichlet boundary conditions, we select the nodes which lie on a fixed edge and set their corresponding degrees of freedom to zero.\
That is, for every fixed node $N_i$, we set $u_{i,x} = u_{i,y} =0$.\
Static condensation of $\mathbf{K}$ then yields a reduced system
\begin{equation}
\hat{\mathbf{K}} \hat{\mathbf{u}} = \hat{\mathbf{f}}_{\mathrm{ext}}
\end{equation}

where $\hat{\mathbf{u}},\hat{\mathbf{f}}_{\mathrm{ext}} \in \mathbb{R}^{q-\Delta q}$ and $\hat{\mathbf{K}} \in \mathbb{R}^{(q-\Delta q) \times (q-\Delta q)}$ with $\Delta q$  the number of fixed degrees of freedom.\
The compliance objective function is then defined as $c:=\mathbf{u}^\top \mathbf{f}_{\mathrm{ext}} \in \mathbb{R}_+$.\
Note that the compliance of the eroded design is always the largest~\cite{Lazarov2016} and therefore

\begin{equation}
\max (c(\mathbf{x}_b(\mathbf{x})), c(\mathbf{x}_e(\mathbf{x})), c(\mathbf{x}_d(\mathbf{x}))) = c(\mathbf{x}_e(\mathbf{x})).
\end{equation}

This paper uses the shorthand notation $ c(\mathbf{x}_e(\mathbf{x}))=c(\mathbf{x}_e)=c_e \in \mathbb{R}$.

\subsection{Formulation}

The compliance minimization problem is stated as follows.

\begin{equation}
\begin{aligned}
\min_{\mathbf{x}} \quad &c(\mathbf{x}_e(\mathbf{x})) \\
\text{s.t.} \quad &v(\mathbf{x}_d(\mathbf{x})) \leq v_{\max} \\
&\mathbf{0} \leq \mathbf{x} \leq \mathbf{1} \\
\end{aligned}
\end{equation}

Here, the state equation is satisfied implicitly and the filter produces $\mathbf{x}_d$ for the volume and $\mathbf{x}_e$ for the compliance.\
All topology optimization examples employ the Method of Moving Asymptotes (MMA)~\cite{Svanberg1987}.\
Its hyperparameters and interface with the employed scalarizations is explained in \ref{app:SlackVariableTreatment}.\

\section{Multiobjective topology optimization} \label{sec:MOTO}

In the topology optimization literature, many works showcase their single-objective optimization algorithms by producing a Pareto curve.\
This is often achieved with either weighted-sum or $\varepsilon$-constraint scalarization.\
Below is a non-exhaustive list of examples.\\

Chen and Wu~\cite{ChenWu1998} use weighted-sum to minimize compliance and maximize the fundamental eigenvalue.\
De Leon et al.~\cite{deLeon2015} use an $\varepsilon$-constraint scalarization for stress-constrained compliant mechanisms.\
Sato et al.~\cite{Sato2017} use an adaptive weighted-sum scalarization and points selection scheme to minimize compliance and mass.\
Ferrari and Sigmund~\cite{Ferrari2019} use $\varepsilon$-constraint scalarization to find trade-offs between compliance and buckling stability.\
Wang et al.~\cite{Wang2022} use weighted-sum to optimize a microfluidic reactor.\
Christensen et al.~\cite{Christensen2023} use weighted-sum for finding Pareto curves between compliance and (local and global buckling) stability.\
Cool et al.\cite{cool_TO_VA} use $\varepsilon$-constraint scalarization for metamaterial optimization of sound transmission loss and connectivity.\
Chen et al.\cite{Chen2024} use weighted-sum to optimize the drag and heat transfer.\
Almeida et al.~\cite{Almeida2025} use weighted-sum to study the trade-off between negative Poisson's ratio and negative thermal expansion.\
Recently, Wotten and Il Yong~\cite{Wotten2026} used weighted-sum to optimize the compliance and build orientation of additively manufactured parts.\
As mentioned, the list is not exhaustive: many more works based on weighted-sum and $\varepsilon$-constraint scalarization exist in the literature.\\

Another class of multiobjective topology optimization methods are based on a so-called ``tracing'' of local optima.\
The idea, also known as homotopy or continuation, is to start from a local optimum and use that solution as the initial guess for successive optimizations for slightly different problems, which in the case of multiobjective optimization involves different scalarization parameters.\
In this way, tracing methods produces locally nondominated solutions, here defined as follows.

\begin{definition} \label{def:locally_nondominated}
A point $\mathbf{f}(\mathbf{x}^\star) \in \mathbf{f}(X)$ is locally nondominated if there exists an $\varepsilon > 0$ such that there is no $\mathbf{x} \in X \bigcap B_n(\mathbf{x}^\star,\varepsilon)$, i.e., no $\mathbf{x} \in X$ in a ball of radius $\varepsilon$ around $\mathbf{x}^\star$, with $f_i(\mathbf{x}) \leq f_i(\mathbf{x}^\star)$ for all $i \in \left\{1,\ldots,m \right\}$ and $f_j(\mathbf{x}^\star) < f_j(\mathbf{x})$ for at least one $j \in \left\{ 1, \ldots, m \right\}$.
\end{definition}

Suresh~\cite{Suresh2010} proposes a tracing method based on topological derivatives to trace the Pareto curve between volume and compliance.\
Zhao~\cite{Zhao2013} uses a similar idea but employs a meshless Galerkin method.\
Takalloozadeh and Ho Yoon~\cite{Takalloozadeh2017} expand the tracing method of Suresh to compliance and stress minimization under thermal loads.\
Arifin et al.~\cite{Arifin2025} use weighted-sum in combination with tracing to design robust metalenses.\
Luchini et al.~\cite{Luchini2025} combine a gradient-informed binary topology optimization algorithm with a tracing approach.\
Although tracing methods are numerically efficient ways to quickly generate a series of at least locally nondominated points, the approaches in the literature are undermined by their inability to ensure a length scale.\
That is, they can result in designs with very small holes that violate the length scale.\
Furthermore, the nondominated set is built up of pieces with differing topology.\
As switching from one piece to another requires more than an infinitesimal change, tracing methods risk getting stuck in a bad locally nondominated set.\\

Multiobjective topology optimization has also been attempted with evolutionary and genetic algorithms~\cite{Kunakote2011,Simonetti2019,Lim2020}.\
This however prohibits the use of high-performing gradient-based methods and is therefore not further discussed or studied in this work.\\

Finally, some works also adapt and improve upon existing scalarizations techniques before applying them to topology optimization.\
Ryu et al.~\cite{Ryu2021} develop an adaptive weighted-sum scalarization and use clustering of the optimal designs to provide insight.\
Crescenti et al.~\cite{Crescenti2021} use the smart Normal Constraint method to optimize for compliance and natural frequency.\\

The above literature overview reveals a wide gap between what is commonly used in the topology optimization community (weighted-sum, $\varepsilon$-constraint) and the state-of-the-art in multiobjective optimization.\
Furthermore, cross-fertilization is relatively limited: attempts to apply more advanced scalarization methods are sparse and, as the numerical examples will show, sometimes based on imprecise assumptions about the nature of the Pareto frontiers that they aim to approximate.\

\section{Numerical examples}
\label{sec:numerical_examples}

This section provides several numerical examples.\
\Cref{sec:truss_opt} studies a simple truss optimization problem to illustrate how local optima cause local Pareto fronts.\
\Cref{sec:VC_Canti,sec:CC_Canti} consider topology optimization of compliance and volume.\
 \Cref{sec:CS_LBeam,sec:CD_Canti} study maximum stress and dynamic compliance minimization.\\

\begin{figure}
    \centering
    \includegraphics[width=0.7\linewidth]{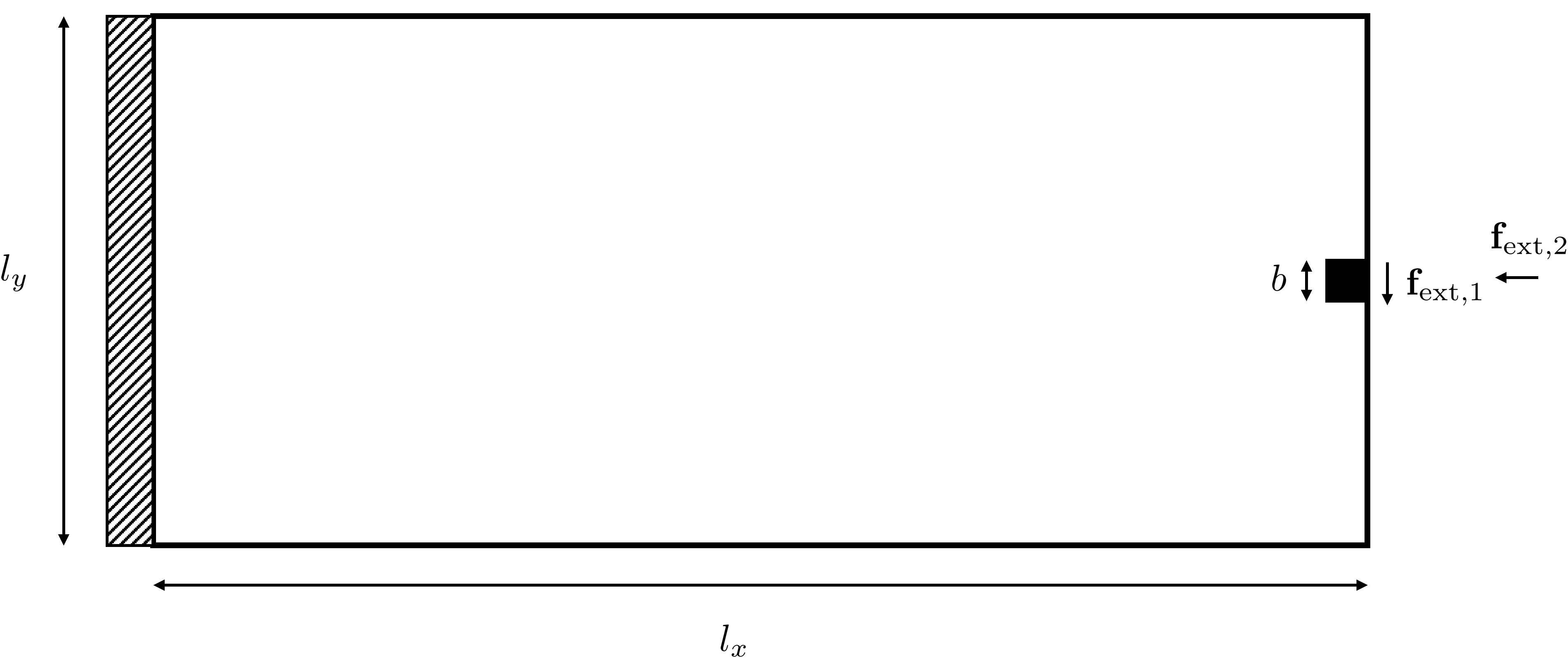}
    \caption{Cantilever problem setup.}
    \label{fig:Canti_setup}
\end{figure}

The compliance examples consider the $l_x$ by $l_y$ cantilever of \Cref{fig:Canti_setup}.\
It is fixed at the left hand side and loaded with two separate loads, a vertical and horizontal force with corresponding force vectors $\mathbf{f}_{\mathrm{ext},1}$ and $\mathbf{f}_{\mathrm{ext},2}$, displacement vectors $\mathbf{u}_1$ and $\mathbf{u}_2$ and compliance functions $c_1$ and $c_2$.\
The load is distributed over a length of $b \in \mathbb{R}$ and the elements in the $b$ by $b$ region next to the load are kept solid.\
Further details are in \ref{app:ComplianceMinimization}.\
Stress and dynamic compliance minimization are detailed in \ref{app:StressMinimization}and \ref{app:DynamicComplianceMinimization}.\
The treatment of the nondifferentiable $\max$ operator of the robust filtering approach is laid out in \ref{app:SlackVariableTreatment}.\\

All examples consider bi-objective problems.\
That is, $m=2$ and so the objective is to minimize $\mathbf{f}(\mathbf{x})=\left( f_1(\mathbf{x}), f_2(\mathbf{x})\right)^\top$.\
\mc{Since the effectiveness of the employed scalarization methods depends on the choice of scaling and scalarization parameters, this is now elaborated upon.\
\paragraph{Scaling} The scalarization methods are applied to scaled versions $f_1^S, f_2^S$ of the objectives $f_1,f_2$, to ensure their values do not exceed $1$.\
This is achieved via the extreme points, which are calculated as follows: for $\mathbf{x}_A \in \argmin \left\{  f_1(\mathbf{x}) \, | \, \mathbf{x} \in X \right\}$  and $\mathbf{x}_B \in \argmin \left\{  f_2(\mathbf{x}) \, | \, \mathbf{x} \in X \right\}$, we set $(f_1^{\min}, f_2^{\star})=\mathbf{f}(\mathbf{x}_A)$ and  $(f_1^{\star}, f_2^{\min})=\mathbf{f}(\mathbf{x}_B)$.\
Then, $f_{i}^S (\mathbf{x}) = f_i(\mathbf{x}) / f_i^\star$.\
Note that these optimization problems are solved with local solvers and thus only generate an approximation of the true value.\ }

\mc {\paragraph{Scalarization} The scalarization parameters are chosen as follows.
\begin{itemize}
\item The weighted-sum scalarization \ref{eq:WS_formulation} minimizes $w_1 f_{1}^S (\mathbf{x}) + w_2 f_{2}^S (\mathbf{x})$, where $w_1$ is evenly spaced in $[0,1]$ and $w_2=1-w_1$.\
\item The $\varepsilon$-constraint scalarization \ref{eq:EC_formulation} solves $\min_{\mathbf{x} \in X}\, f_1^S(\mathbf{x}) \, \mathrm{s.t.} \, f_2^S(\mathbf{x}) \leq \varepsilon_2$, with $\varepsilon_2$ uniformly in $\left[f_2^{\min}/f_2^{\star},1 \right]$, and vice versa for the $\min f_2(\mathbf{x})$ variant.\
\item The Pascoletti-Serafini scalarization \ref{eq:PS_formulation} constructs the line connecting $(f_1^{\min}/f_1^{\star},1)$ and $(1, f_2^{\min}/ f_2^{\star})$, also known as the Convex Hull of Individual Minima (CHIM).\
Then, it chooses $\mathbf{r} \in \mathbb{R}_2^{++}$ perpendicular to the CHIM, i.e.,
\begin{equation}
\mathbf{r} = \left( 1-\frac{f_2^{\min}}{f_2^{\star}},  1-\frac{f_1^{\min}}{f_1^{\star}} \right)^\top,
\end{equation}
and $\mathbf{a}$ equidistantly on the CHIM.\
\end{itemize}
}

\mc{All scalarization parameter choices can be translated to more than two objectives.\
For \ref{eq:WS_formulation}, the weights can be sampled from the unit simplex $\left\{ \mathbf{w}: w_i \in [0,1],\, \sum_{i=1}^m w_i = 1 \right\}$.\
For \ref{eq:EC_formulation}, the $\varepsilon$ values should be chosen in a higher-dimensional plane.\
For \ref{eq:PS_formulation}, $\mathbf{a}$ and $\mathbf{r}$ can be chosen on and perpendicular to the CHIM, respectively, provided the individual minima are affinely independent.\
However, any extension to $m > 2$ is nontrivial since covering the true nondominated set requires considering many scalarization parameters.\
We refer the reader to relevant literature on this topic~\cite{Messac2003,Eichfelder2008,MuellerGritschneder2009,Eichfelder2021}.\ 
The $m>2$ is not further treated to limit the scope of this paper.\ }

\begin{figure}
    \centering
    \includegraphics[width=0.95\linewidth]{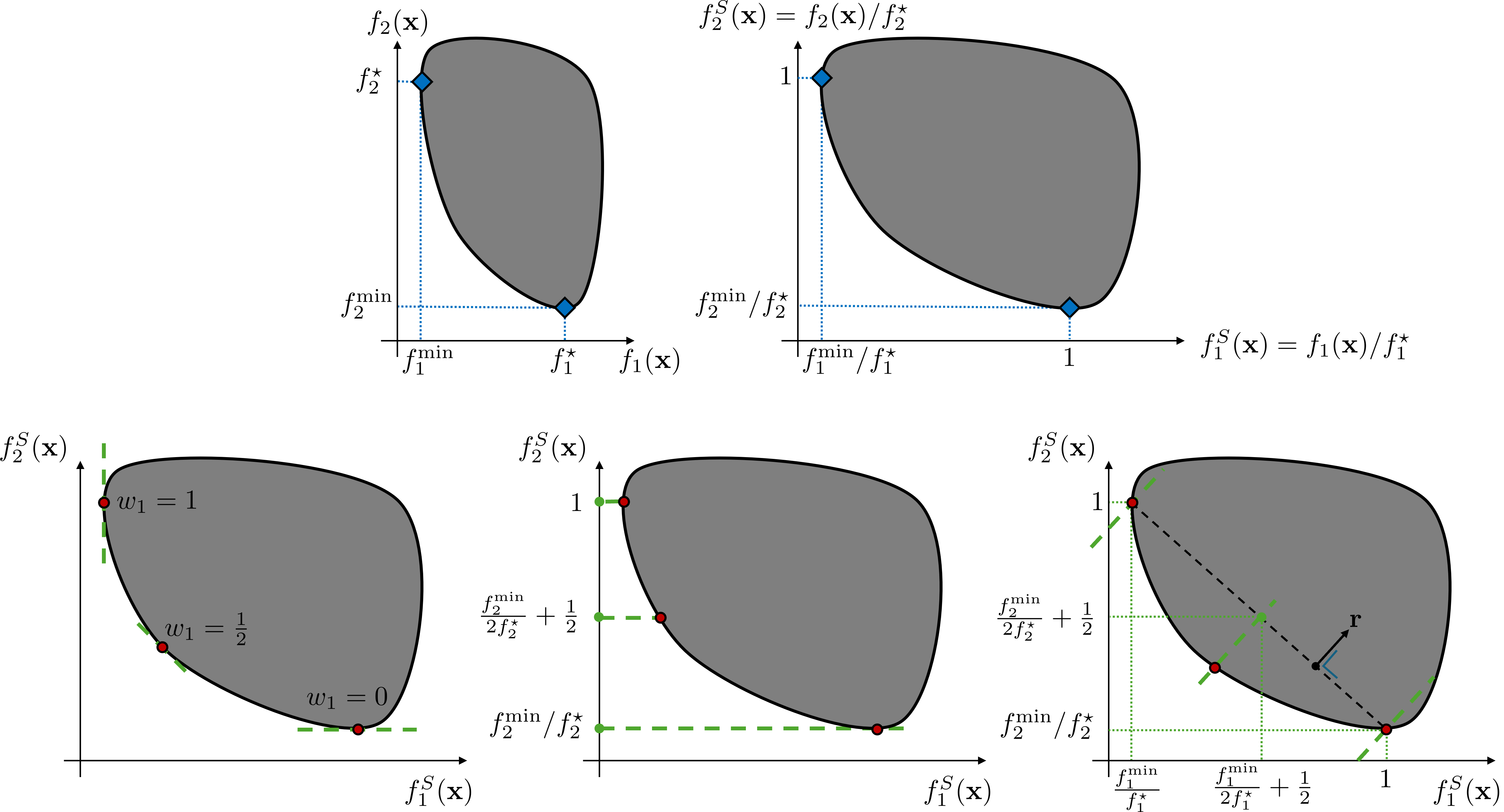}
    \caption{\mc{Illustration on the choice of scalarization parameters. Top: scaling of the objective space, by finding the points indicated by the blue diamonds and dividing through $f_1^{\star}$ and $f_2^\star$. Bottom: illustration of scalarization with three points, for \ref{eq:WS_formulation} (left), \ref{eq:EC_formulation} ($\min f_1(\mathbf{x})$ variant, middle) and \ref{eq:PS_formulation} (right). Note that the axes show the scaled objectives. The red points are obtained after solving with the scalarization parameters indicated geometrically in green. For \ref{eq:WS_formulation}, the green dashed lines indicate the slope of $y = \frac{w_1}{w_2} x + C$. For \ref{eq:EC_formulation}, the green points indicate $\varepsilon_2$ values in the $f_2^S$ axis, with horizontal dashed lines to indicate where they intersect with $\mathbf{f}^S (X)$. For \ref{eq:PS_formulation}, $\mathbf{a}$ is chosen uniformly on the CHIM and $\mathbf{r}$ is perpendicular to the CHIM.}}
    \label{fig:ScalarizationParameterIllustration}
\end{figure}
 
\subsection{Truss optimization} \label{sec:truss_opt}

\mc{The definitions of (weak) efficiency and nondominance in \Cref{sec:MOO} rely on a notion of global optimality.\
Unfortunately, because topology optimization is inherently nonconvex, it is difficult to obtain global optima.\
To simplify the setting and illustrate the concept of local Pareto fronts, this section considers a truss optimization example.} \\

\begin{figure}[h]
\centering
\begin{subfigure}[t]{0.33\linewidth}
  \centering
  \includegraphics[width=\textwidth]{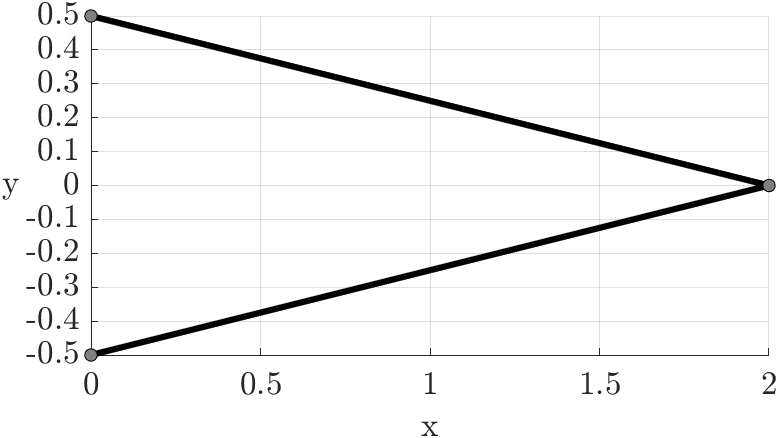}
 \caption{Two-bar truss}
  \label{fig:2bar}
\end{subfigure}%
\begin{subfigure}[t]{0.33\linewidth}
  \centering
  \includegraphics[width=\textwidth]{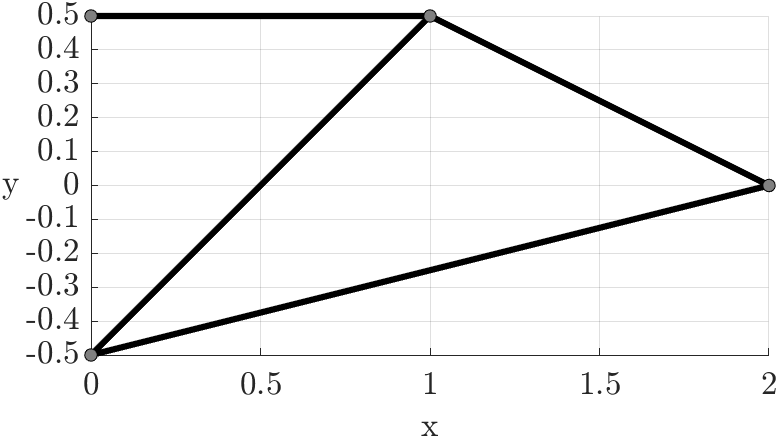}
 \caption{Four-bar truss}
  \label{fig:4bar}
\end{subfigure}%
\begin{subfigure}[t]{0.33\linewidth}
  \centering
  \includegraphics[width=\textwidth]{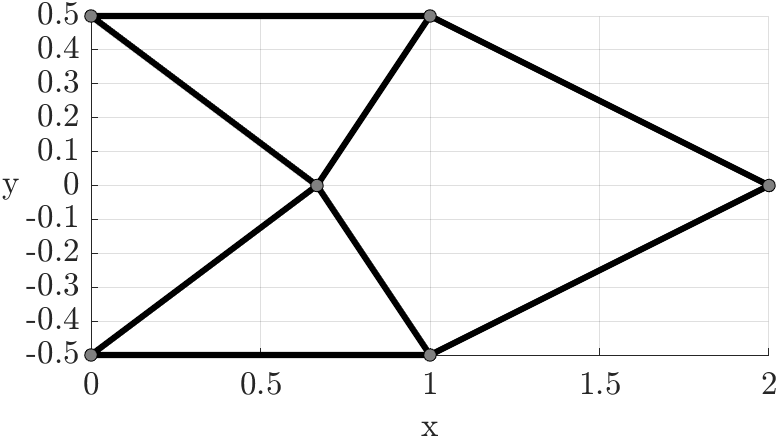}
 \caption{Eight-bar truss}
  \label{fig:8bar}
\end{subfigure}%
\caption{Three truss layouts.}
\label{fig:groundstructures}
\end{figure}

\Cref{fig:groundstructures} shows three truss ground structures: a symmetric two- and eight-bar truss and an asymmetric four-bar truss.\
They represent three stable and non-crossing truss topologies for the considered boundary conditions.\
Note that although there exist other truss topologies, e.g. with five or six members, these were omitted for brevity.\\

Consider compliance and volume minimization of each of these truss layouts under a vertical load at the node on the far right (at $(x,y)=(2,0)$).\
The bar areas, collected in the bar area vector $\mathbf{a}_{\mathrm{bar}}$, must be optimized and a volume and length scale constraint is added.\
This leads to the following formulation.

\begin{equation} \label{eq:truss_opt_A}
\begin{aligned}
\min_{\mathbf{a}_{\mathrm{bar}}} \quad &\left(c(\mathbf{a}_{\mathrm{bar}}), v(\mathbf{a}_{\mathrm{bar}}) \right)^\top \\
 \text{s.t.} \quad &v(\mathbf{a}_{\mathrm{bar}}) \leq v_{\max} \\
 &  \mathbf{a}_{\min} \leq \mathbf{a}_{\mathrm{bar}} \\
\end{aligned}
\end{equation}

For this section, $c$ and $v$ are redefined as the compliance and volume function dependent on the bar areas $\mathbf{a}_{\mathrm{bar}}$.\
The truss compliance is defined as $c(\mathbf{a}_{\mathrm{bar}})=\mathbf{u}^\top \mathbf{f}_{\mathrm{ext}}$ with $\mathbf{f}_{\mathrm{ext}}$ the external force vector and $\mathbf{u}$ the displacement vector.\
The latter is implicitly dependent on $\mathbf{a}_{\mathrm{bar}}$ via the state equation $\mathbf{K}(\mathbf{a}_{\mathrm{bar}})\mathbf{u}=\mathbf{f}$, where $\mathbf{K}=\sum_i a_i\mathbf{K}_i$ is a linearly weighted sum of elemental stiffness matrices.\
The truss stiffness $k(\mathbf{a}_{\mathrm{bar}})$ is defined as the inverse of the compliance: $k(\mathbf{a}_{\mathrm{bar}})=1/c(\mathbf{a}_{\mathrm{bar}})$.\

The volume is $v(\mathbf{a}_{\mathrm{bar}})=\mathbf{a}_{\mathrm{bar}}^\top \mathbf{l}$, with $ \mathbf{l}$ containing the bar lengths.\
The lower bound on the areas is uniform, i.e., $\mathbf{a}_{\min} =a_{\min}\mathbf{1}$ with $\mathbf{1}$ a vector of ones.\
It enforces a minimum length scale of $a_{\min} > 0$.\\

The trusses are statically determinate and no member can vanish due to the length scale ($a_{\min} > 0$), so $\mathbf{K}(\mathbf{a}_{\mathrm{bar}})$ is invertible and the compliance is convex in the bar areas.\
The volume constraint and compliance objective are linear in the bar areas and the length scale is a simple box constraint.\
Hence, the truss optimization problem of \Cref{eq:truss_opt_A} is convex~\cite{Svanberg1981}: local optima are also global optima.\\

\begin{figure}[h]
\centering
\begin{subfigure}[t]{0.275\linewidth}
  \centering
  \includegraphics[width=\textwidth]{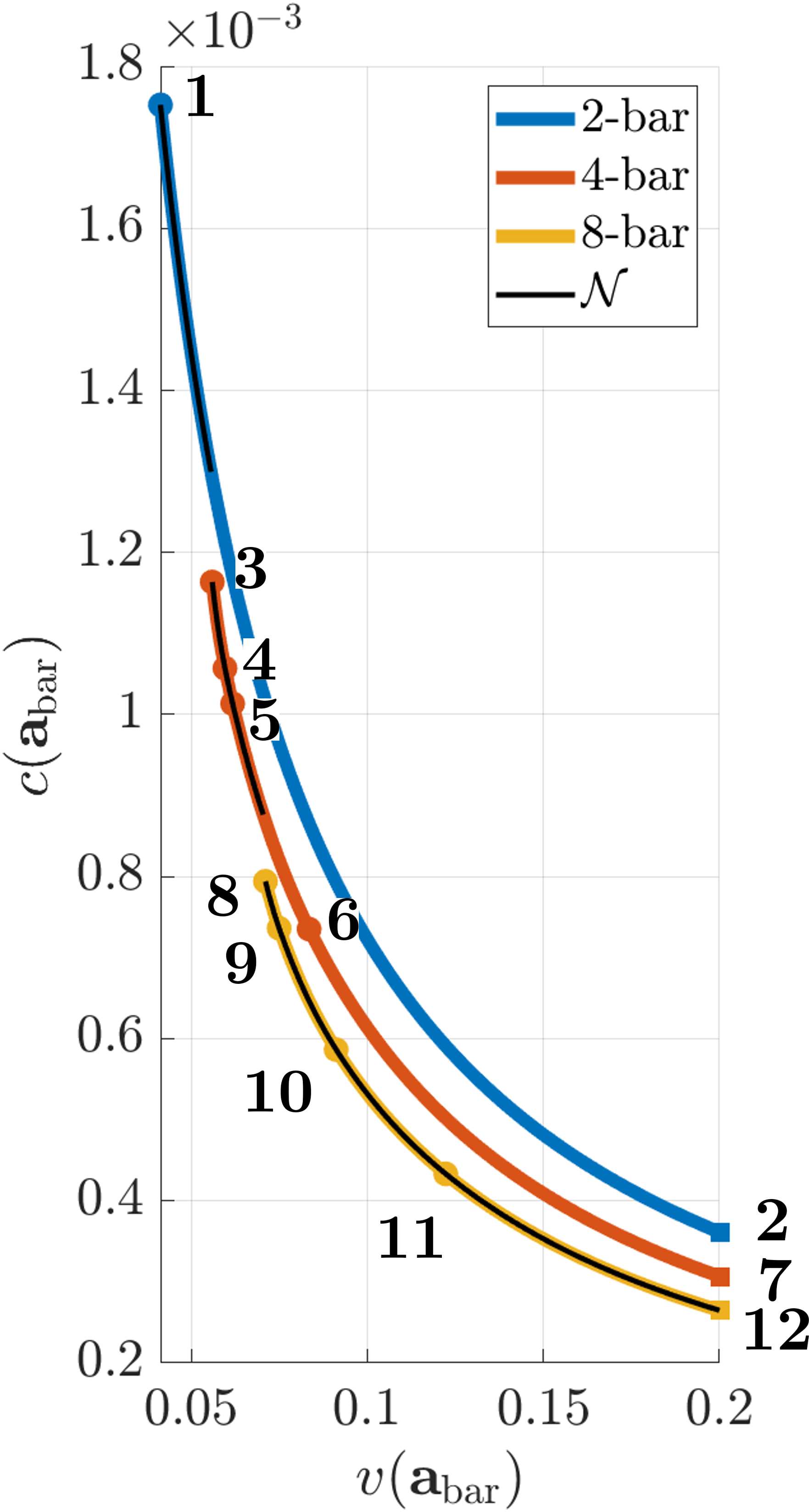}
  \caption{Volume versus compliance.}
  \label{fig:truss_opt_A2}
\end{subfigure}%
\begin{subfigure}[t]{0.275\linewidth}
  \centering
  \includegraphics[width=\textwidth]{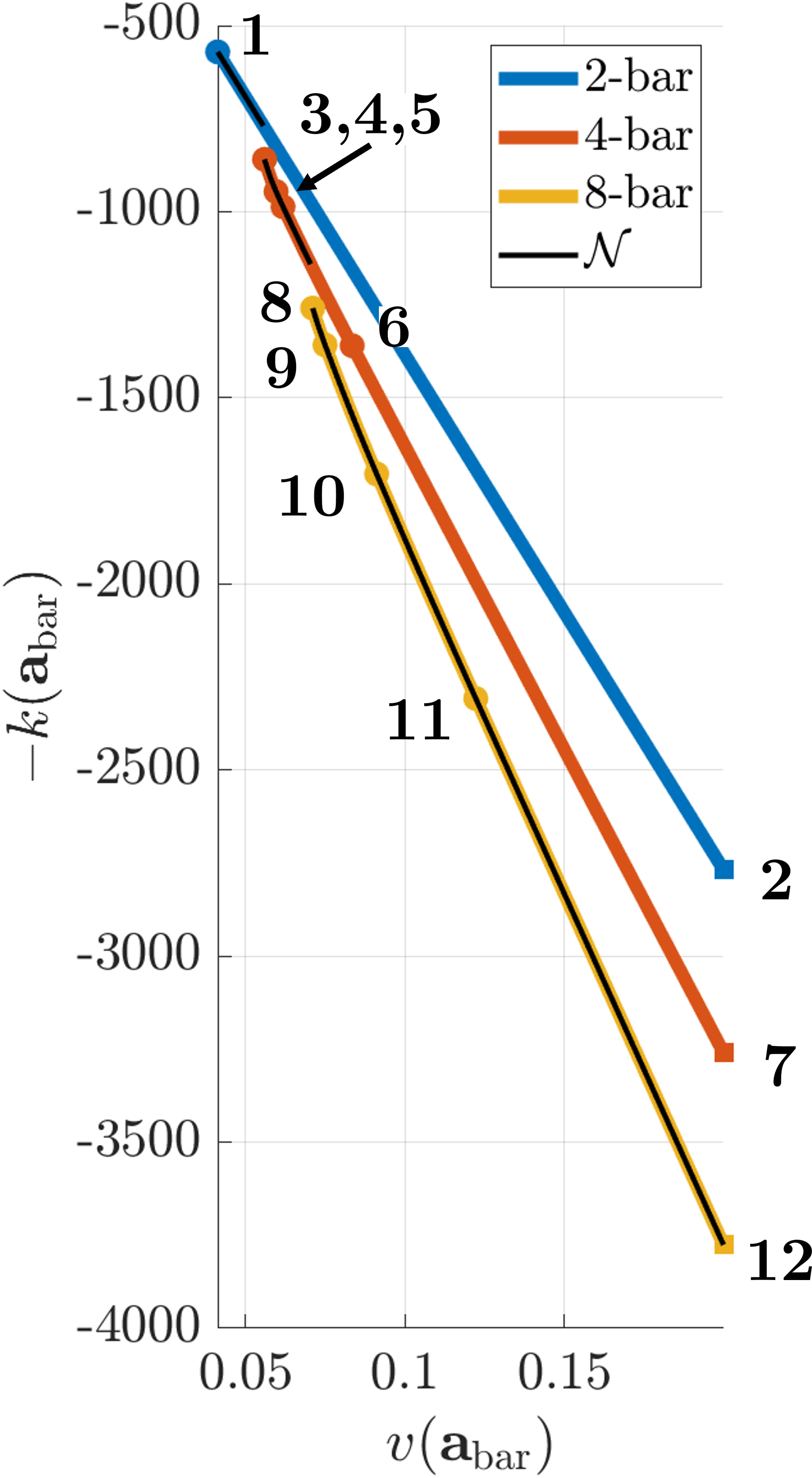}
  \caption{Volume versus stiffness.}
  \label{fig:truss_opt_A3}
\end{subfigure}%
\begin{subfigure}[t]{0.45\linewidth}
  \centering
  \includegraphics[width=\textwidth]{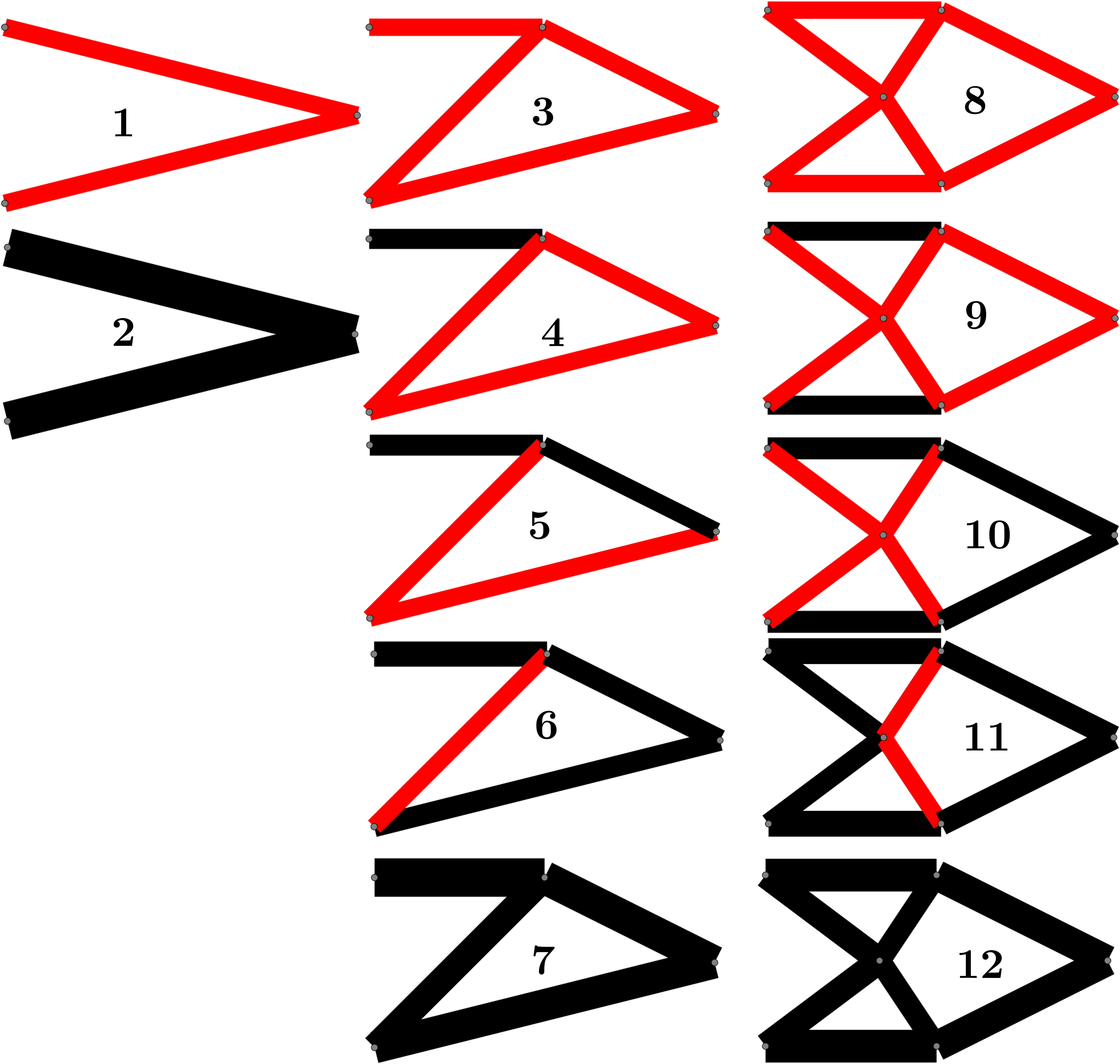}
\caption{Designs corresponding to points and squares on the left figures. Red bars have area $a_{\min}$.}
    \label{fig:truss_opt_A1}
\end{subfigure}%
\caption{Nondominated sets and corresponding designs when optimizing for the bar areas $\mathbf{a}_{\mathrm{bar}}$. The solid black line denotes the (exact) solution of the superproblem.}
\label{fig:truss_opt_A}
\end{figure}

\Cref{fig:truss_opt_A} shows the nondominated points when minimizing both volume and compliance, for each of the truss layouts.\
They are obtained by finding the individual minimum for the volume ($v_{\min}$ occurs at $\mathbf{a}_{\mathrm{bar}}=\mathbf{a}_{\min}$) and doing successive $\varepsilon$-constraint scalarizations between $v_{\min}$ and $v_{\max}$.\
For a fixed truss layout (two-bar, four-bar or eight-bar), the optimization problem of \Cref{eq:truss_opt_A} is from now on referred to as a truss subproblem.\
If the choice of the groundstructure is also considered, the problem is referred to as the truss superproblem.\
As stated, each subproblem is convex and hence \mc{its solutions are global minima}.\\

\Cref{fig:truss_opt_A2,fig:truss_opt_A3} plot the nondominated points of each subproblem in the volume-compliance and volume-stiffness space, respectively.\
The corresponding designs are in \Cref{fig:truss_opt_A1}.\
\mc{The curves represent the local Pareto fronts of the superproblem, which is defined more rigorously below.\
In what follows, the individual properties of the curves, which are local Pareto fronts of the superproblem, are first discussed.\
Then, their composition into the nondominated set of the superproblem is considered.} \\

\mc{The three volume-stiffness curves are continuous, piecewise linear.\
Therefore, the three volume-compliance curves are continuous, piecewise hyperbola.\
This structure arises from the linear behavior of the truss elements.\
In the absence of $a_{\min}$, the optimal bar diameters can be linearly scaled without losing global optimality: halving $v(\mathbf{a}_{\mathrm{bar}})$ also halves each optimal bar diameter, which in turn halves the stiffness and doubles the compliance.\
The \emph{piecewise} linearity is then caused by the length scale: when an optimal bar diameter would go below $a_{\min}$, it instead remains at $a_{\min}$.\
The other bars still scale proportionally, leading to a different but still linear $(v(\mathbf{a}_{\mathrm{bar}}), \,k(\mathbf{a}_{\mathrm{bar}}))$ behaviour.\
This successive activation of the length scale constraint is reflected in \Cref{fig:truss_opt_A1}: as the volume decreases, an increasing number of bars reach $a_{\min}$, as indicated in red.}\\ 

\mc{The solution of the superproblem is found by comparing the solutions of the subproblems.\
As the subproblems are solved exactly, owing to their convexity, the obtained nondominated set $\mathcal{N}$ is not an approximation but an exact solution and drawn with a solid black line in \Cref{fig:truss_opt_A2,fig:truss_opt_A3}.\
Crucially, $\mathcal{N}$ is not continuous: along the front, the optimal ground structure changes and this causes discontinuities in compliance and stiffness.\
In particular, there are two discontinuities.\
First, around $v(\mathbf{a}_{\mathrm{bar}})=0.07$ the optimal groundstructure switches from the eight-bar symmetric truss to the asymmetric four-bar truss.\
Second, around $v(\mathbf{a}_{\mathrm{bar}})=0.055$ the groundstructure again switches from the asymmetric four-bar truss to the symmetric two-bar truss.\
Both changes in the optimal groundstructure occur when a local Pareto front terminates, which occurs when all bar diameters reach the minimum length scale $a_{\min}$.}\\

\begin{figure}[h]
\centering
\begin{subfigure}[t]{0.275\linewidth}
  \centering
  \includegraphics[width=\textwidth]{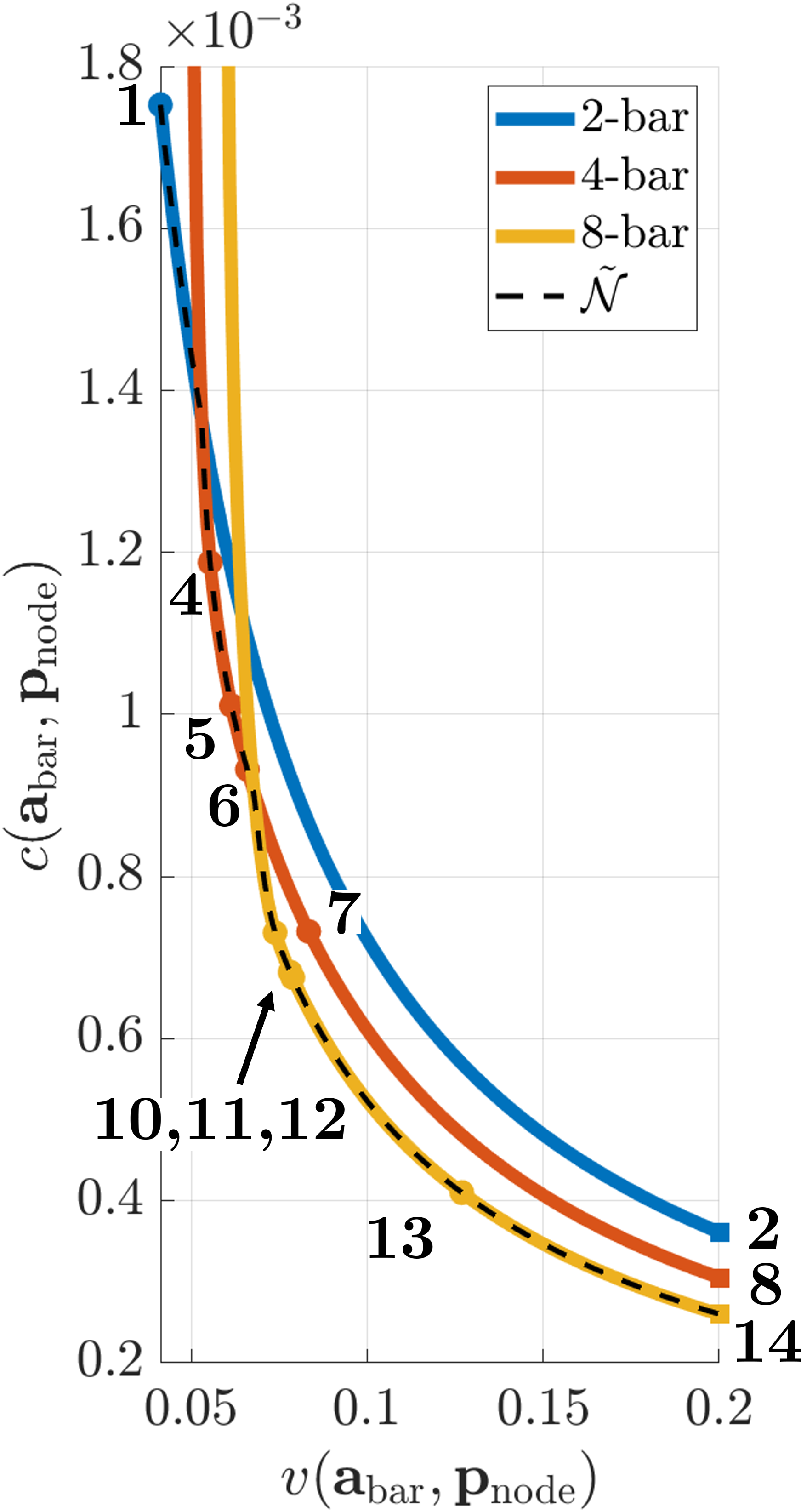}
  \caption{Volume versus compliance.}
  \label{fig:truss_opt_PA2}
\end{subfigure}%
\begin{subfigure}[t]{0.275\linewidth}
  \centering
  \includegraphics[width=\textwidth]{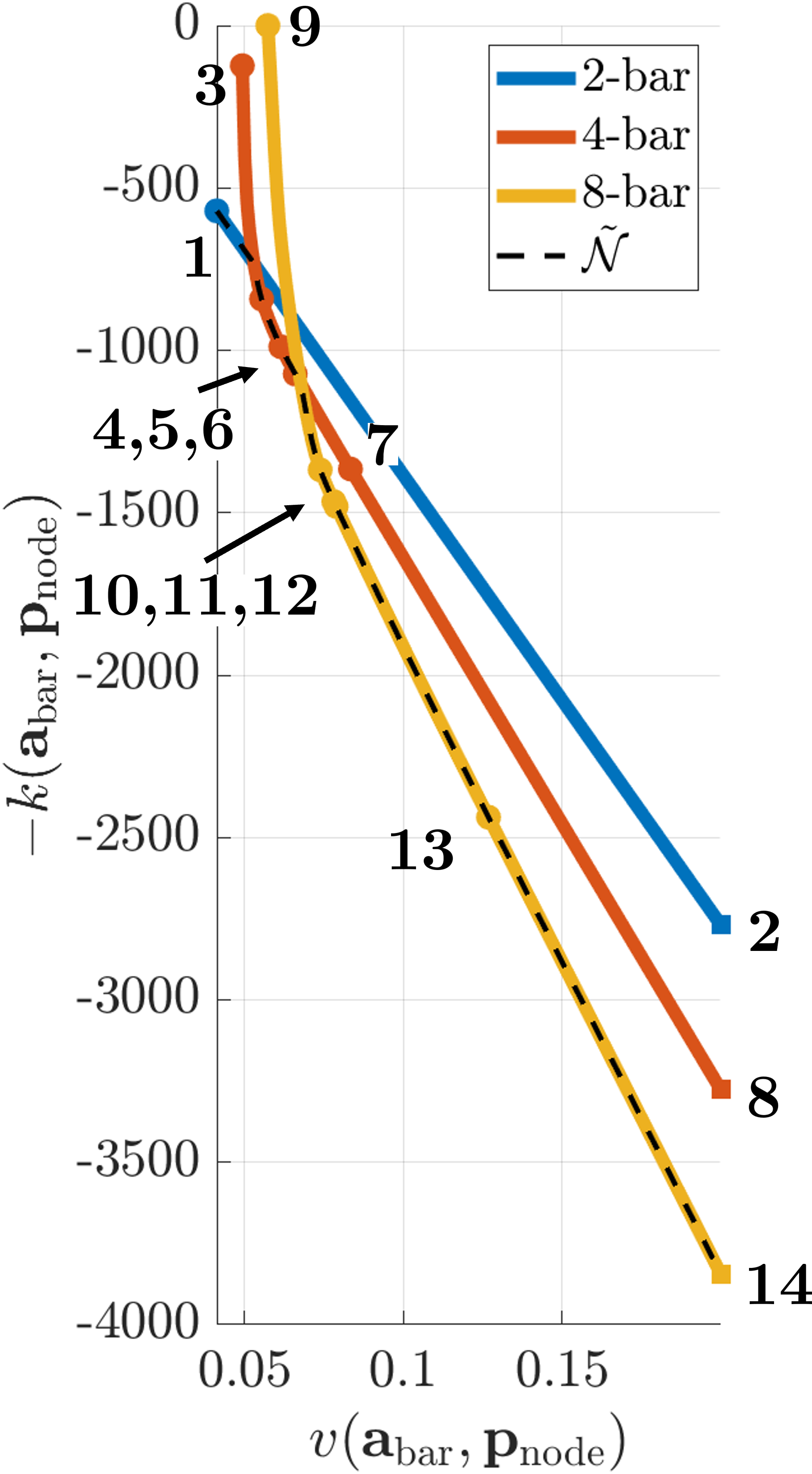}
  \caption{Volume versus stiffness.}
  \label{fig:truss_opt_PA3}
\end{subfigure}%
\begin{subfigure}[t]{0.45\linewidth}
  \centering
  \includegraphics[width=\textwidth]{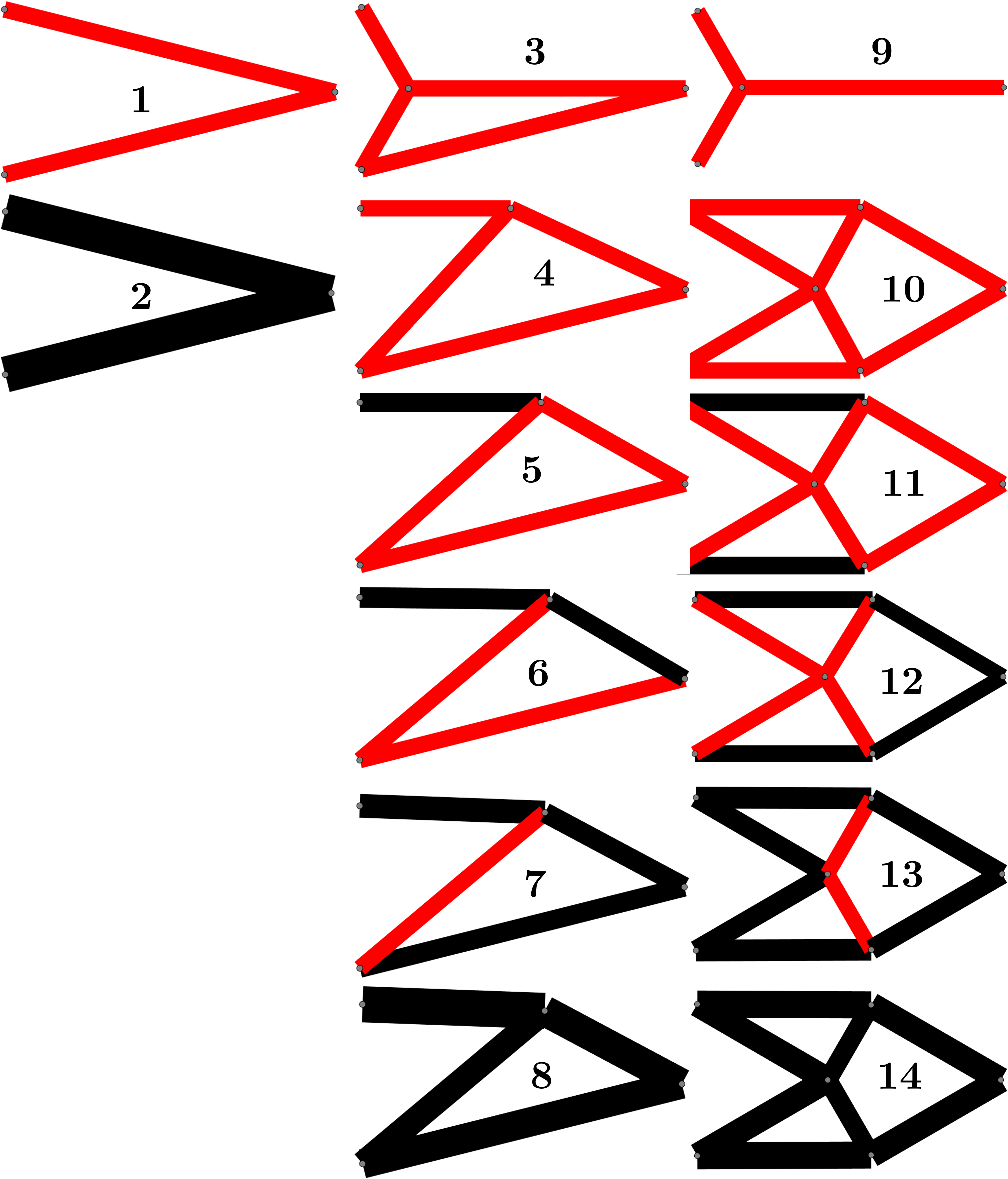}
\caption{Designs corresponding to points and squares on the left figures. Red bars have area $a_{\min}$.}
\end{subfigure}%
\caption{Nondominated sets and corresponding designs when optimizing for bar areas \emph{and} node positions $\mathbf{p}_{\mathrm{node}}$. \mc{The dashed black line denotes the approximate solution of the superproblem.}}
\label{fig:truss_opt_PA}
\end{figure}

\mc{The nondominated set in \Cref{fig:truss_opt_A2,fig:truss_opt_A3} is nonconvex and discontinuous.\
With increased design freedom, however, the latter can be avoided.\
Consider to this end an adaptation of \Cref{eq:truss_opt_A} where the truss node positions $\mathbf{p}_{\mathrm{node}}$ are also optimization parameters.\
The numerical results for this case are in \Cref{fig:truss_opt_PA}.\
Note that, due to the nonlinear dependence of the objectives on $\mathbf{p}_{\mathrm{node}}$, the subproblems are no longer convex and the nondominated set of the superproblem $\tilde{\mathcal{N}}$ is only an approximation.\
It is therefore drawn with a dashed line.}\\

\mc{The additional design freedom extends the local Pareto frontiers: instead of terminating when all bars reach $a_{\min}$, volume is further reduced by moving the nodes such that the bars shorten.\
This leads to highly compliant truss layouts and causes the local frontiers to intersect.\
As a result, $\tilde{\mathcal{N}}$ is now continuous.\
However, the corresponding design still undergoes the same two discontinuous topology changes: from symmetric eight-bar to asymmetric four-bar to symmetric two-bar.\
At these groundstructure transitions, the nondominated set remains continuous but is not differentiable due to a discontinuous derivative.}\\

\mc{The remaining numerical examples in this section consider topology optimization.\
Given the truss examples in this section, it should come as no surprise that \emph{multiobjective topology optimization problems with length-scale constraints have Pareto frontiers composed of segments of local  frontiers, each corresponding to a distinct topological layout}.\
Consequently, Pareto frontiers in topology optimization may be nonconvex and exhibit discontinuities in both value and derivative.}\\

For clarity, we here provide a definition of a local Pareto frontier.\
\begin{definition} \label{def:local_frontier}
Consider a connected subset $\mathcal{E}_{\mathrm{loc}}(\mathbf{f}, X) \subseteq X \subseteq \mathbb{R}^n$ such that (i) its image $\mathbf{f}(\mathcal{E}_{\mathrm{loc}})$ is locally nondominated and (ii) for all $\mathbf{x}_1,\mathbf{x}_2 \in \mathcal{E}_{\mathrm{loc}}$ there exists a continuous map $\gamma: [0,1] \rightarrow X$ such that $\gamma(0)=\mathbf{x}_1,\gamma(1)=\mathbf{x}_2$ and $\gamma(t) \in \mathcal{E}_{\mathrm{loc}}$ for all $t \in [0,1]$.\
Then $\mathbf{f}(\mathcal{E}_{\mathrm{loc}})$ is a local Pareto frontier.
\end{definition}

\mc{Before moving to the strongly nonconvex topology optimization problems, it is important to emphasize that the discontinuous behaviour of the frontier indeed stems from the imposed length scale.\
In the absence of a length scale, the optimization problem becomes convex and the efficient and nondominated sets are connected sets.\
Along the frontier, design and performance behave smoothly.\
For example, in the above truss optimization case, letting $a_{\min} \rightarrow 0$ would result in the eight-bar truss remaining dominant.\
Likewise, for a variable-thickness sheet formulation ($P=1$), the optimization problem is also convex.\
When the upper thickness bound is inactive, the optimal thickness distribution can be rescaled while retaining global optimality.\
If regions attain the maximum thickness ($x=1$), the corresponding densification changes the curvature of the frontier but not its smoothness nor its continuity.\
The same observations hold for a single-load-case volume-compliance minimization, where a solution with a rank-2 laminate constitutes the global optimum of the convex problem.\
The observed discontinuities in the efficient and nondominated set in the examples that follow are thus a direct result of the length scale.\ }

\subsection{Topology optimization of volume versus compliance} \label{sec:VC_Canti}

This section considers the cantilever example of \Cref{fig:Canti_setup}, but only focuses on the volume $v(\mathbf{x}_d)$ and the compliance $c_1(\mathbf{x}_e)$ under the first load case, for the remainder of this section simply referred to as $c(\mathbf{x}_e)$.\
Since compliance can become unbounded, it is bounded by $5 c_{\mathrm{solid}}$, with $c_{\mathrm{solid}}=c(\mathbf{1})$ the compliance of the fully solid structure.\
\mc{The $5c_{\mathrm{solid}}$ value is used to rescale the compliance objective.\ }
The volume fraction $v(\mathbf{x}_d)$ is already in $[0,1]$ and hence not normalized further.\ 
The optimization problem is written as 

\begin{equation}\label{eq:VC_Canti}
\begin{aligned} 
\min_{\mathbf{x}} \quad &\mathbf{f}(\mathbf{x})=\left(c(\mathbf{x}_e(\mathbf{x})),v(\mathbf{x}_d(\mathbf{x})) \right)^\top \\
\text{s.t.} \quad &c(\mathbf{x}_e(\mathbf{x})) \leq 5c_{\mathrm{solid}} \\
&\mathbf{0} \leq \mathbf{x} \leq \mathbf{1} \\
\end{aligned}
\end{equation}
where $x_i=1$ for the elements in the $b$ by $b$ square (see \Cref{fig:Canti_setup}).\
\mc{The $\varepsilon$-constraint and Pascoletti-Serafini scalarizations must be preceded by a volume and a compliance minimization run to find the extreme points of the Pareto frontier, which are used to scale and bound their scalarization parameters.\
The former yields the point $(v(\mathbf{x_d}),c(\mathbf{x_e}))=(v_{\min}, 5c_{\mathrm{solid}})$.\
The latter is omitted, since it always yields the solid structure with $(v(\mathbf{x_d)},c(\mathbf{x_e}))=(1,c_{\mathrm{solid}})$.\ }

\subsubsection{Weighted-sum scalarization} \label{sec:WS100}

\begin{figure}[!htb]
\centering
\begin{subfigure}[t]{0.45\textwidth}
  \centering
\captionsetup{width=0.95\textwidth}
  \includegraphics[width=0.95\textwidth]{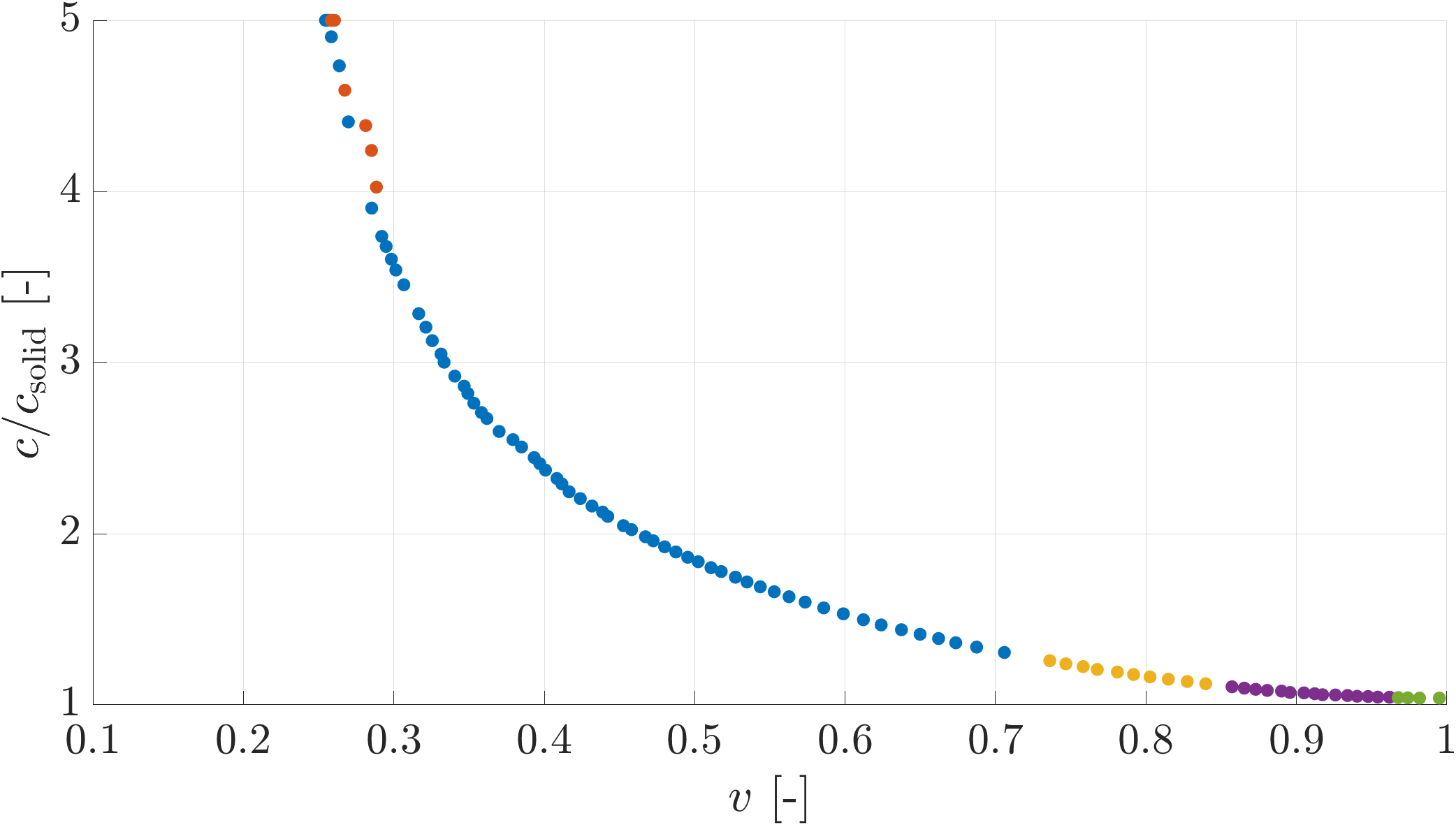}
  \caption{Pareto frontier for $R=0.1$.} \label{subfig:frontier_WS_R20}
\end{subfigure}%
\begin{subfigure}[t]{0.5\textwidth}
  \centering
\captionsetup{width=0.95\textwidth}
  \includegraphics[width=0.95\textwidth]{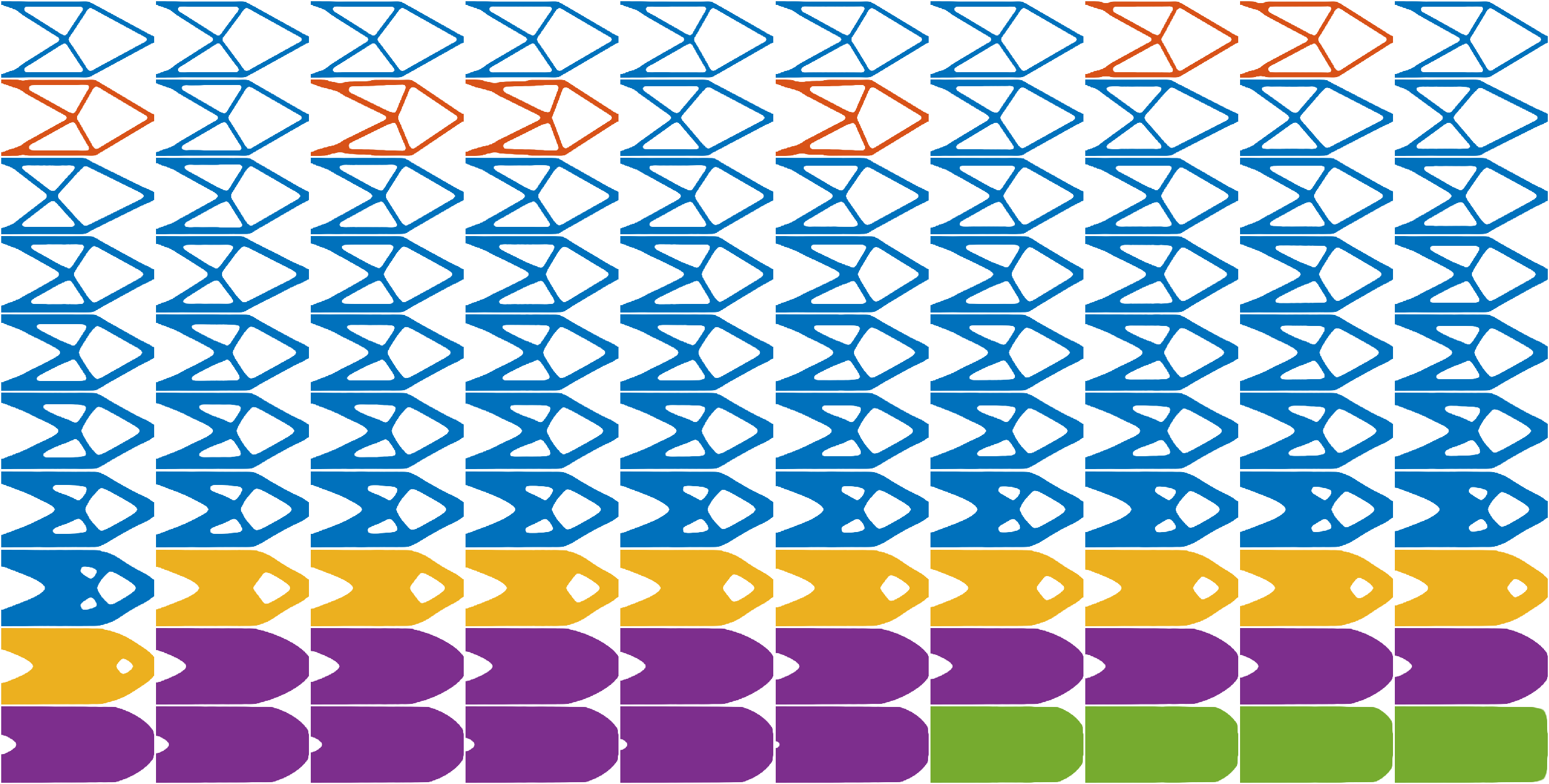}
  \caption{Designs for $R=0.1$. Blue and orange designs have the same layout and topology, but the orange designs are partially disconnected from the top and bottom edge.} \label{subfig:designs_WS_R20}
\end{subfigure}%
\\
\begin{subfigure}[t]{0.45\textwidth}
  \centering
\captionsetup{width=0.95\textwidth}
  \includegraphics[width=0.95\textwidth]{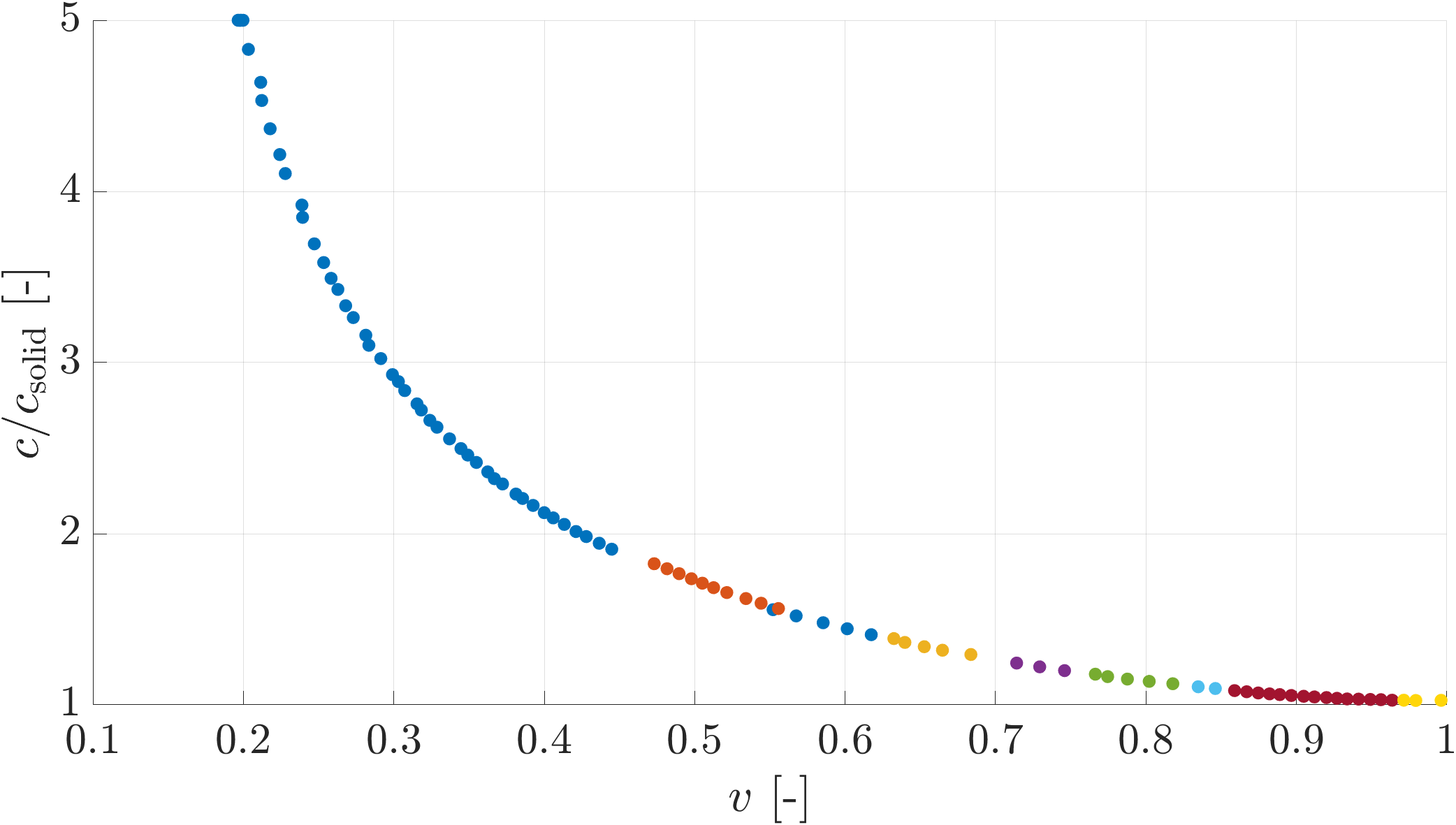}
  \caption{Pareto frontier for $R=0.05$.} \label{subfig:frontier_WS_R10}
\end{subfigure}%
\begin{subfigure}[t]{0.5\textwidth}
  \centering
\captionsetup{width=0.95\textwidth}
  \includegraphics[width=0.95\textwidth]{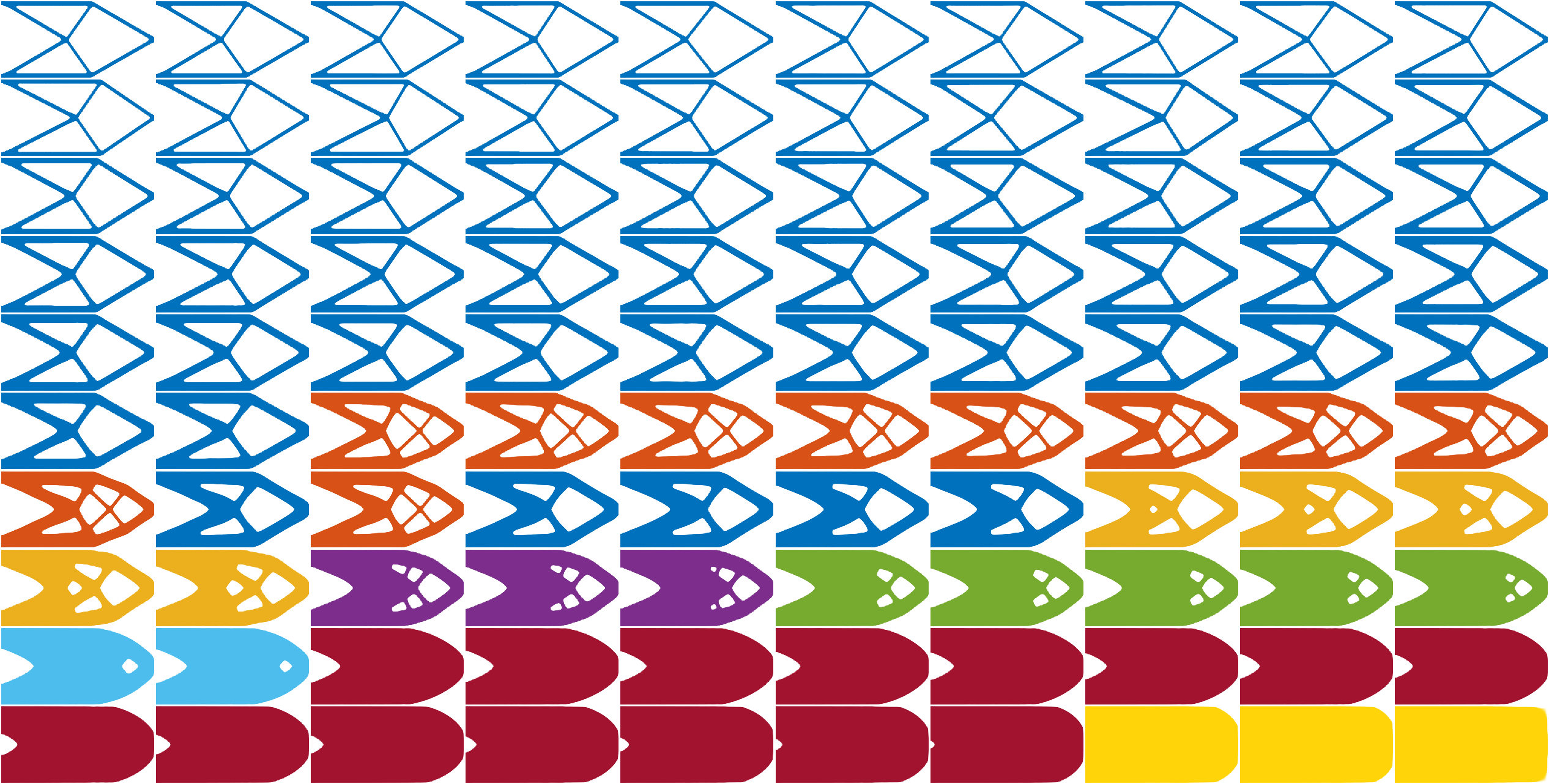}
  \caption{Designs for $R=0.05$.} \label{subfig:designs_WS_R10}
\end{subfigure}%
\\
\begin{subfigure}[t]{0.45\textwidth}
  \centering
\captionsetup{width=0.95\textwidth}
  \includegraphics[width=0.95\textwidth]{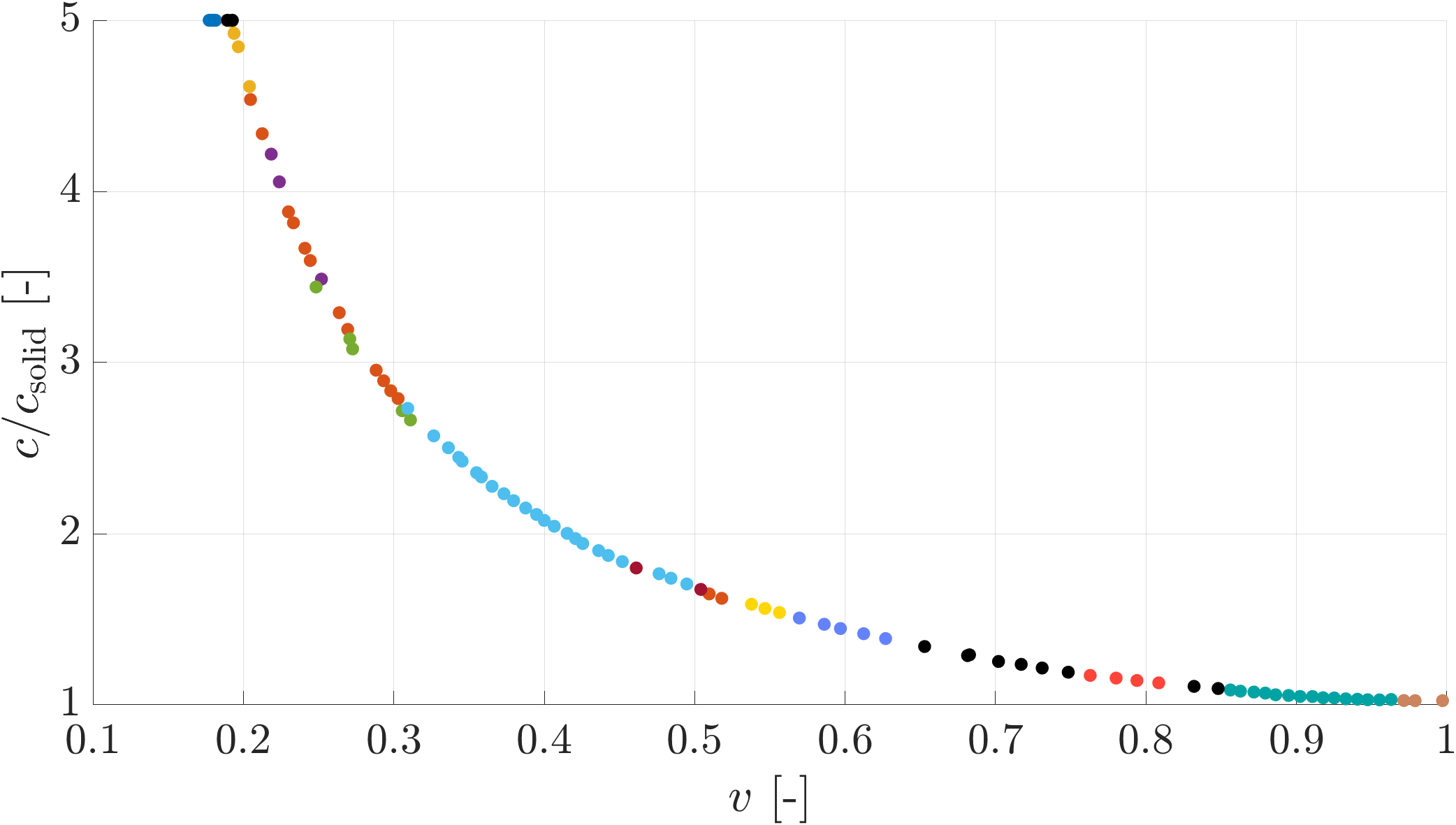}
  \caption{Pareto frontier for $R=0.025$.} \label{subfig:frontier_WS_R5}
\end{subfigure}%
\begin{subfigure}[t]{0.5\textwidth}
  \centering
\captionsetup{width=0.95\textwidth}
  \includegraphics[width=0.95\textwidth]{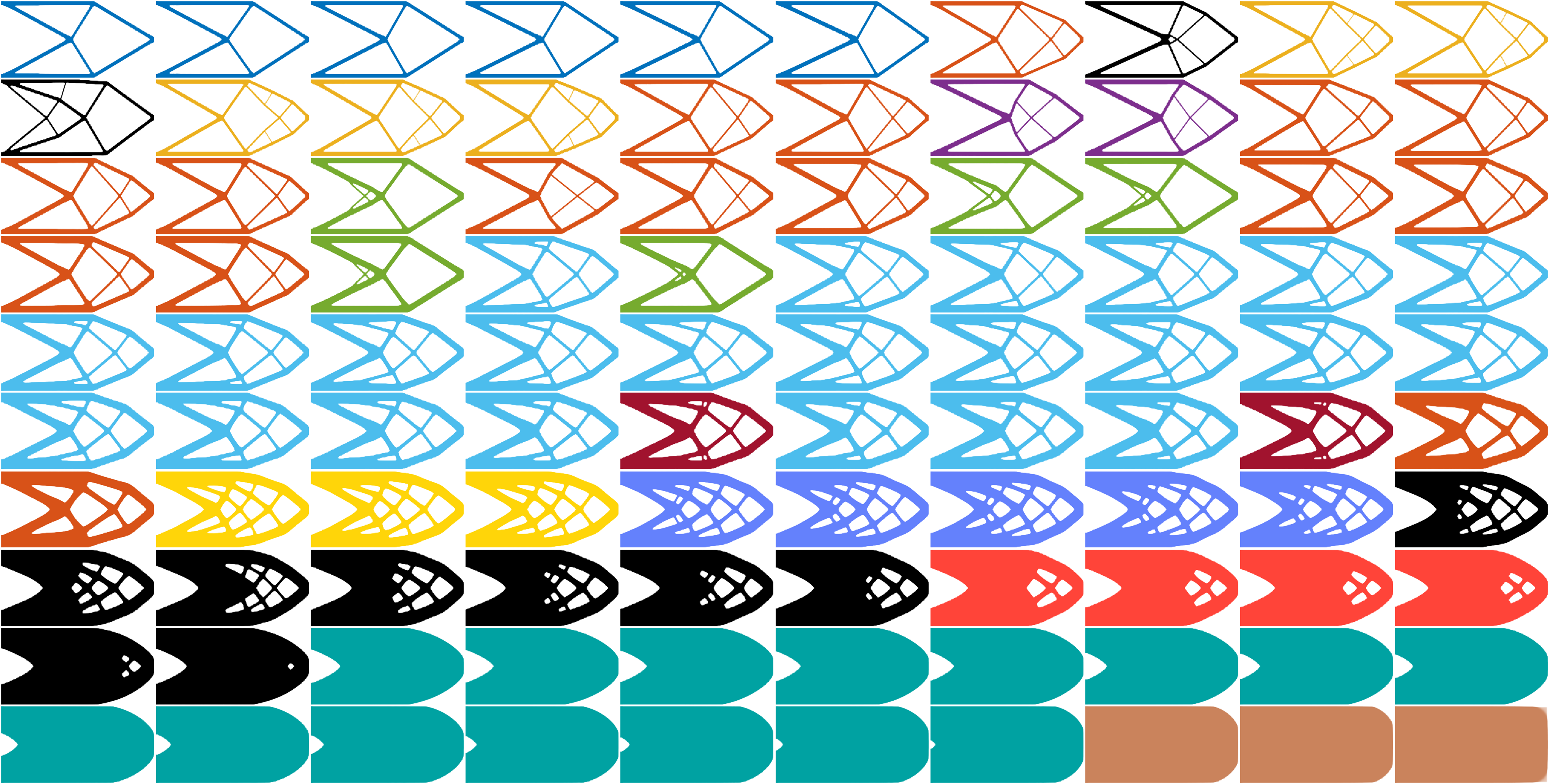}
  \caption{Designs for $R=0.025$.} \label{subfig:designs_WS_R5}
\end{subfigure}%
\caption{Pareto curves (left) and designs (right) for a $2$-by$1$ cantilever optimized for volume and compliance with weighted-sum scalarization. Three filter radii ($R=0.1$, $0.05$ and $0.025$ element widths) (top, middle, bottom). Designs are blueprint ($\eta=0.5$), ordered by volume, from left to right and top to bottom, and colored by beam layout and topology. Black is reserved for designs alone in their class. }
\label{fig:VC_Canti_WS100}
\end{figure}

The simplest way to determine efficient solutions of \Cref{eq:VC_Canti} is with weighted-sum scalarization.\
\Cref{fig:VC_Canti_WS100} shows the nondominated set approximation and corresponding designs for $100$ weights distributed uniformly in $[0,1]$, for the cantilever with $l_x=2$, $l_y=1$  meshed with $n_x=400$ by $n_y=200$ elements.\
Three filter radii are considered, $R=0.025$, $R=0.05$ and $R=0.1$, i.e., $5$, $10$ and $20$ element widths respectively, to show the effect of the minimum length scale.\
\mc{The coloring is used to indicate designs with a different topology and beam layout, as this signifies a discontinuity in the design and possibly objective space.\
Additionally, \Cref{subfig:designs_WS_R20} also uses a different color (blue and orange) to indicate whether the 8-bar truss designs stick to the top and bottom edges.\ }

\mc{Several observations can be made.\
These are divided into observations regarding (i) approximation quality obtained with the weighted-sum scalarization and (ii) the properties of the underlying Pareto frontier.} 

\mc{\paragraph{Approximation quality} The weighted-sum scalarization provides a rather poor approximation.\
In particular, we attribute two deficiencies to the scalarization method itself.\ }
\mc{
\begin{itemize}
\item First, the sample density is highly irregular: in some regions the points are close together, in other regions there are large spacings.\
For example, the blue points in \Cref{subfig:frontier_WS_R10} cluster around $v = 0.45$ but are more spaced out around $v=0.6$.\
This phenomenon can be at least partially attributed to the fact that weighted-sum scalarization can be understood as a parametrization of the Pareto frontier based on its slope, which naturally oversamples regions with high curvature if weights are uniformly sampled.\
\item Second, the approximations contain several large holes.\
This is most striking in \Cref{subfig:frontier_WS_R20}, which has an empty region around $v(\mathbf{x}_d)=73$\% where the design goes from one to three internal holes.\
This is attributed to the nonconvexity of the optimization problem and the lack of any control on the obtained local optimum: small changes in the scalarization parameter (i.e., in the weights) lead to other local frontiers, which have other slopes.\
For example, in \Cref{subfig:frontier_WS_R10} the design switches from the blue 8-bar truss to the orange 12-bar truss around $v \approx 0.45$, leaving a significant hole in the frontier approximation even though the scalarized weights only changed slightly.\
Then, around $v=0.55$, the blue 8-bar truss is recovered again.\
\end{itemize}
}
\mc{The dominated orange points in \Cref{subfig:frontier_WS_R20} correspond to eight-bar trusses whose bars are partially detached from the edges of the design domain.\
Numerical experiments shown later in this work use different methods but also observe this behaviour, and hence we do not attribute this deficiency to the weighted-sum scalarization.\
This deficiency is simply attributed to convergence of the local solver to a bad local optimum, and not to the weighted-sum scalarization.} \\ 

\mc{Summarizing, the weighted-sum scalarization provides irregular approximations due to (i) a lack of control regarding the local frontier a given weight $(w_1, w_2)$ converges to, and (ii) the slope-based spacing of the obtained points.\
The lack of control is true for any local solver, but it compounds the effect of the latter.}

\mc{
\paragraph{Pareto front properties} Despite the poor approximation quality, two general conclusions about the Pareto frontier in volume-compliance minimization can be made.\
First, the local Pareto frontiers seem to be convex.\
Any nonconvexities, such as for the blue points in \Cref{subfig:frontier_WS_R20} around $v=0.4$, are slight and could be attributed to numerical issues, such as premature convergence or simulation accuracy.\
Second, the length scale has a direct effect on the number of local optima, and hence local frontiers.\
The experiments at the smallest length scale, $R=0.025$, admit twelve local frontiers with more than one sample and eleven designs (in black) with a unique topology.\ }

\mc{Further conclusions regarding the frontier's convexity and continuity are complicated by the approximation quality.\
The next section turns to a continuation procedure to shed further light on these properties.\ }

\subsubsection{Tracing via epsilon-constraints} \label{sec:trace_VC}

The previous section observed that the weighted-sum scalarization has little control over the optimization path and eventual topology of the obtained design.\
This complicates conclusions regarding the convexity and continuity of the Pareto frontier.\
To obtain more control on design topology and illustrate local frontiers, this section combines $\varepsilon$-constraint scalarization with a crude tracing methodology.\
Three designs are selected from \Cref{subfig:designs_WS_R20} and used as initial guess of optimization runs with successively tightening bounds on either volume or compliance.\
The $\beta$ projection parameter is kept at its maximum value to ensure the procedure acts as a level set method that cannot create new holes.\\

\begin{figure}
    \centering
    \includegraphics[width=\linewidth]{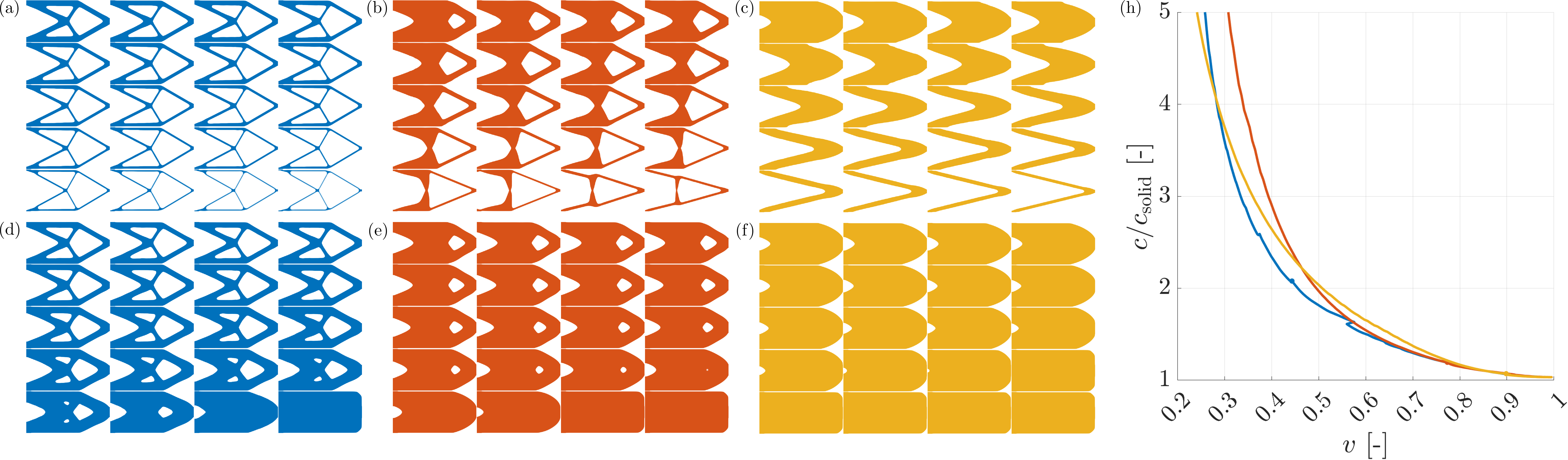}
    \caption{Local Pareto frontiers corresponding to three distinct topologies, found with a crude tracing method. Subfigures (a-f) show the designs corresponding to each subfrontier, starting with the original design top left and the final design at bottom right. The top (a-c) shows the tracing towards minimal volume, the bottom (d-f) shows tracing towards maximal volume. Subfigure (g) shows the subfrontiers in the image space.}
    \label{fig:trace_VC}
\end{figure}

\Cref{fig:trace_VC} shows the obtained local Pareto frontiers and corresponding designs.\
\mc{The local frontiers are drawn with solid lines to indicate the (relatively) smooth underlying design changes.\
The discontinuities of each local frontier are attributed to the crude continuation methodology, which always starts from a slightly infeasible guess to push designs to lower volume and compliance.\
More importantly, the experiment reveals a strong similarity to the truss example of \Cref{fig:truss_opt_PA}.\
That is, the Pareto frontier is smooth but experiences a discontinuous jump in the design space: the design switches topology, going from yellow at high compliance to blue, then orange and finally, again to yellow.\
Although each design change is discontinuous, the frontier remains continuous.\
Furthermore, the discontinuous change in the frontier's slope gives rise to a nonconvex region.\
Overall, the local frontiers in \Cref{fig:trace_VC} show striking resemblance to the truss example of \Cref{fig:truss_opt_PA}.\ }

\mc{The continuation procedure is imperfect in its control of the topology.\ 
For example, the blue design goes from three internal holes to one and then zero.\
The blue curve is thus again misleading, as this topological transition cannot be smooth.\
Nevertheless, the employed continuation procedure remains a valuable exploratory tool and should be further investigated in future research.\ }

A comparison can be made with the work of Ryu et al.~\cite{Ryu2021}, which proposes an automatic clustering scheme to group designs into categories, just like the manual coloring in \Cref{fig:VC_Canti_WS100}.\
However, our results show that the premise underlying such clustering somewhat misrepresents the structure of the frontier.\
Although designs can indeed be grouped together based on topology, it is shown in this work that these clusters overlap significantly in the image space.\
Furthermore, designs in the same cluster are not just point clouds in the image space but instead portions of local Pareto frontiers, which have a much richer mathematical structure.\\

\mc{The employed continuation procedure can also be compared to the tracing approaches in the topology optimization literature,} most notably the work of Suresh~\cite{Suresh2010}.\
There, the topological derivative approach creates locally nondominated points lying on smooth curves in both the image and design space.\
The results of this section show that this \mc{is} misleading.\
The smoothness observed in the design space is caused by the absence of an explicit length scale.\
That is, the produced designs can have infinitesimally small holes, in practice only limited by the discretization.\
This is undesirable in itself, but it also smooths the discontinuous derivatives of the Pareto frontier and hence hides this underlying property.\
Additionally, tracing methods based on topological derivatives therefore carry a risk: they may become stuck in a local Pareto frontier, from which they cannot escape.\
Just like there exists low-performing local optima in topology optimization (see e.g.~\cite{DeWeer2025}), there can also exist low-performing local frontiers.\
\mc{Future research on continuation or tracing methodologies should aim to avoid or at least acknowledge these pitfalls.\ }

\subsubsection{Comparing weighted-sum, $\varepsilon$-constraint and Pascoletti-Serafini scalarization} \label{sec:WS_EC_PS_VC}

This section compares the weighted-sum, $\varepsilon$-constraint and Pascoletti-Serafini scalarization for the 2-by-1 cantilever, meshed with $400$ by $200$ elements and with a filter radius of $R=0.1$ (i.e., $20$ element widths).\
\mc{The upper bound of $5c_{\mathrm{solid}}$ is again placed on the compliance and used to rescale the compliance.}\\

\begin{figure}
    \centering
    \includegraphics[width=0.9\linewidth]{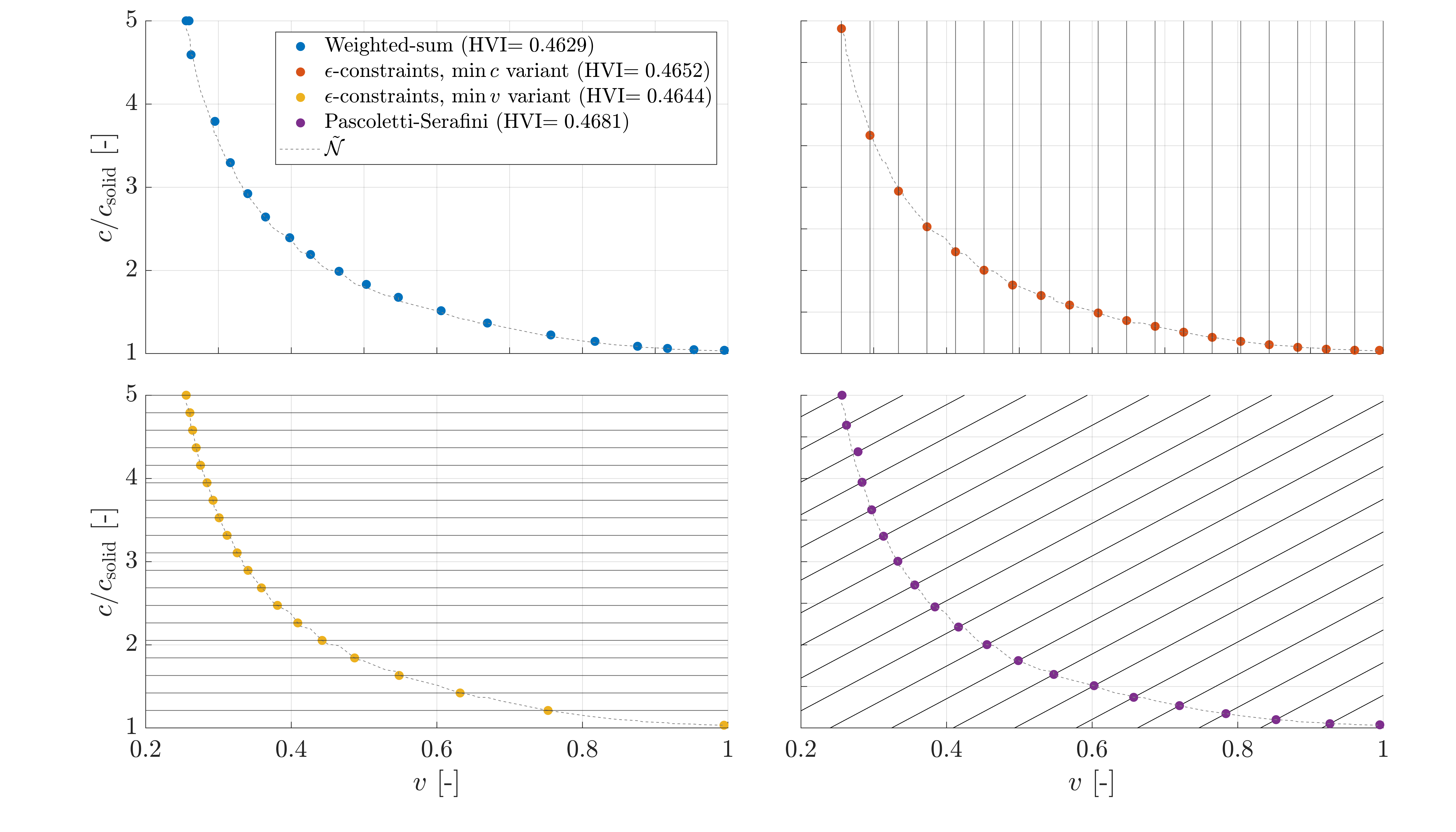}
    \caption{Comparison between weighted-sum, two $\varepsilon$-constraint variants and Pascoletti-Serafini scalarization for the volume-vs-compliance 2-by-1 cantilever. \mc{The dashed grey line connects the nondominated points of all four scalarizations combined. It does not imply a smooth transformation from one solution to another but merely illustrates the outer hull of the approximation $\tilde{\mathcal{N}}$. The legend quantifies the hypervolume indicator HVI.}}
    \label{fig:VC_Canti}
\end{figure}

\mc{\Cref{fig:VC_Canti} compares the Pareto frontier approximation, with each scalarization using 20 points.\
The weighted-sum scalarization produces an irregular approximation, with the same deficiences noted before.\
This is again attributed to (i) the lack of control offered by this scalarization and (ii) the local Pareto frontiers which strongly influence the optimization path.\
The two $\varepsilon$-constraint variants and the Pascoletti-Serafini scalarization produce more regularly spaced points.\
Due to the shape of the frontier, the volume minimization variant has tightly spaced points at low volumes and loose spacing at high volumes.\
Conversely, the compliance minimization variant has tight spacing at high volumes and loosely spaced points at low volumes.\
The Pascoletti-Serafini scalarization provides an almost even spacing.\
This is attributed to the choice of $\mathbf{r}$ as perpendicular to the CHIM, which ensures \ref{eq:PS_formulation} acts as a ``weighted sum'' of the two $\varepsilon$-constraint variants.\
In other words, due to the rescaling and choice of $\mathbf{a}$ and $\mathbf{r}$, \ref{eq:PS_formulation} produces point whose spacing is less sensitive to the shape of the frontier.\
This is also supported by the hypervolume indicator HVI quantified in the legend of \Cref{fig:VC_Canti}, which computes the area of the region weakly dominated by a pointwise approximation (see~\cite{HVI} and \ref{app:HVI}).\
With \ref{eq:PS_formulation}, the largest area is captured (HVI$=0.4681$), followed by the two $\varepsilon$-constraint variants (HVI$=0.4652$ and $0.4644$) and then the weighted-sum scalarization (HVI=$0.4629$).}\

\subsubsection{Asymmetry under symmetric boundary conditions}\label{sec:asym_VC_Canti}

\mc{Until now, each topology optimization run started with a symmetric initial guess and, with very few exceptions, also returned a symmetric solution.\
However, the truss example of \Cref{sec:truss_opt} shows that, even with symmetric boundary conditions, an asymmetric design can outperform a symmetric one for low volume fractions.\
As discussed, the superiority of the asymmetric four-bar truss is caused by the length scale: local frontiers corresponding to symmetric structures are terminated or pushed to become less efficient when bars reach the minimum area $a_{\min}$.\
This leaves regions where asymmetric solutions can shine, because they can have an intermediate numbers of beams.}\\

\mc{Previous studies have also focused on the optimality of asymmetric truss structures: see} Stolpe~\cite{Stolpe2010}, Rozvany~\cite{Rozvany2011}, Cheng and Liu~\cite{ChengLiu2011}, Guo et al.~\cite{Guo2012} and Richardson et al.~\cite{Richardson2013}, with differing conclusions depending on the exact setting and assumptions.\
\mc{Symmetry in topology optimization was also studied by White and Voronin~\cite{WhiteVoronin2018}, who indicate that (at $v_{\max}=0.5$) the symmetry of a cantilever's optimal layout depends on its length.}\\

\mc{The goal of this section is to (i) highlight the optimality of asymmetric solutions, (ii) demonstrate how a random sampling methodology produces hard-to-interpret results and (iii) identify how this methodological shortcoming could be improved in future research.\
In particular, we call for an approach that uses continuation of the scalarization parameters to explore the local frontiers.}\\

\begin{figure}[!htb]
\centering
\begin{subfigure}[t]{0.5\textwidth}
  \centering
\captionsetup{width=0.95\textwidth}
  \includegraphics[width=0.95\textwidth]{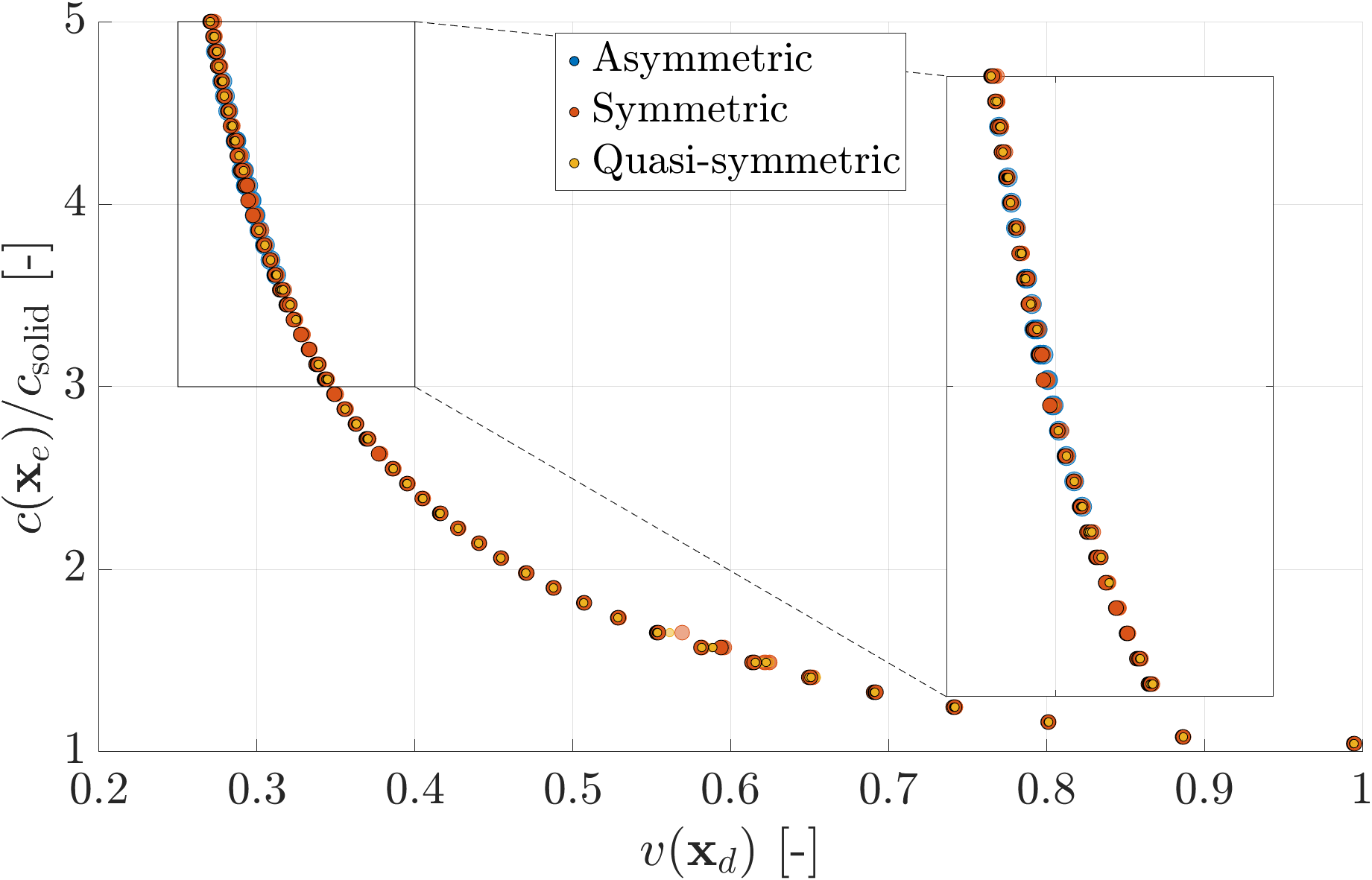}
  \caption{$l_x=2$} \label{subfig:asym_200}
\end{subfigure}%
\begin{subfigure}[t]{0.5\textwidth}
  \centering
\captionsetup{width=0.95\textwidth}
  \includegraphics[width=0.95\textwidth]{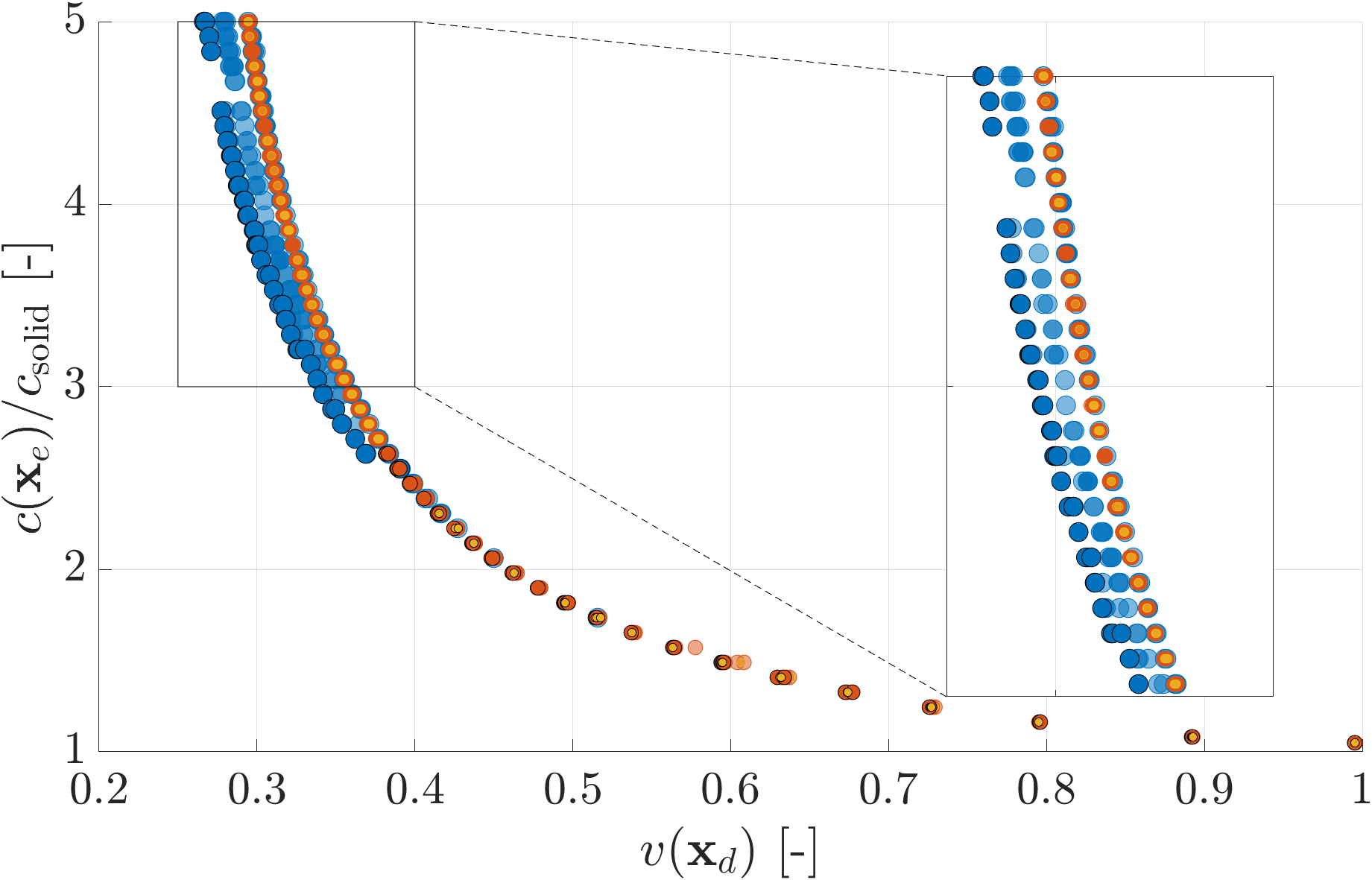}
  \caption{$l_x=2.5$} \label{subfig:asym_250}
\end{subfigure}%
\\
\begin{subfigure}[t]{0.5\textwidth}
  \centering
\captionsetup{width=0.95\textwidth}
  \includegraphics[width=0.95\textwidth]{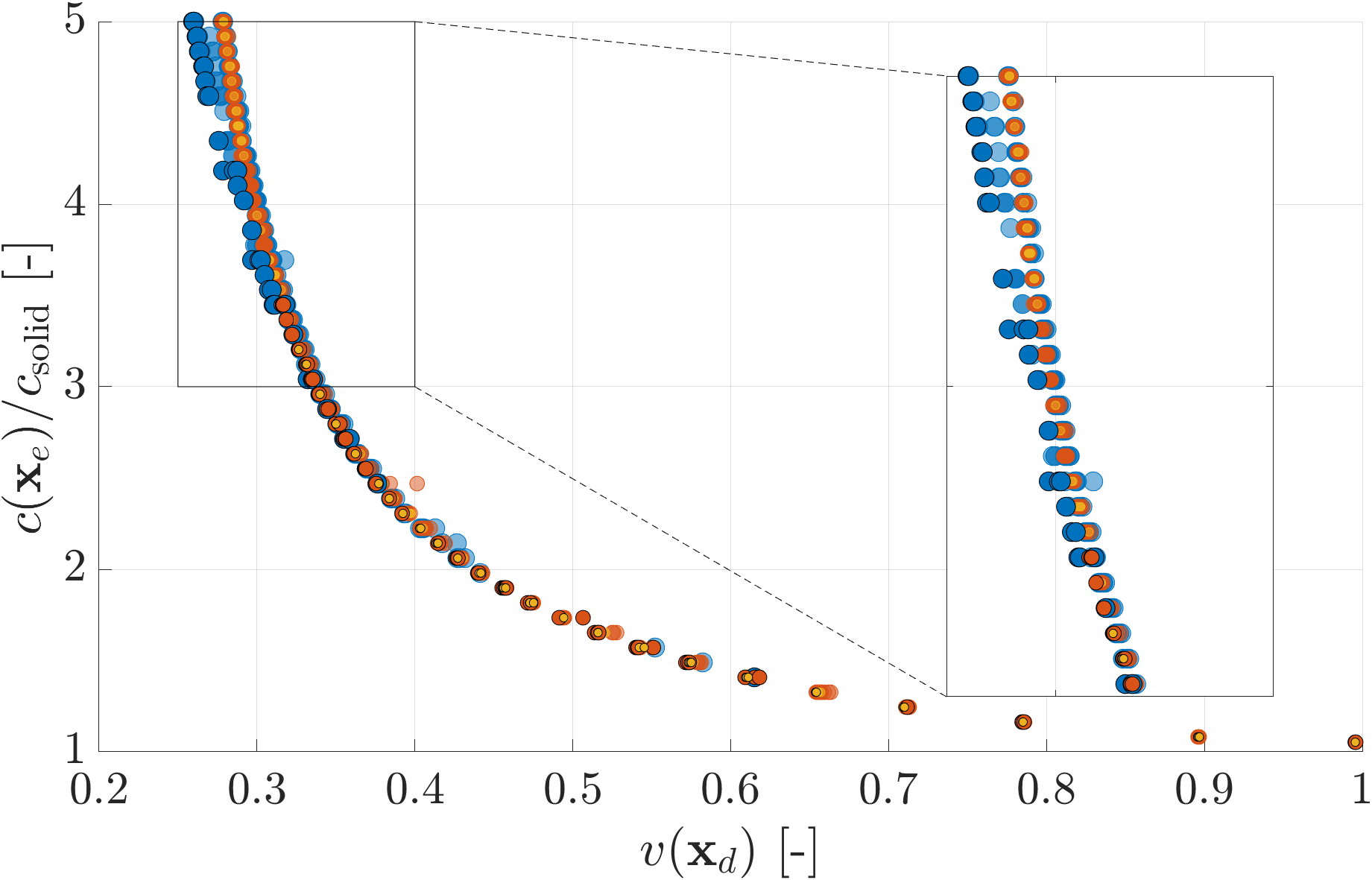}
  \caption{$l_x=3$} \label{subfig:asym_300}
\end{subfigure}%
\begin{subfigure}[t]{0.5\textwidth}
  \centering
\captionsetup{width=0.95\textwidth}
  \includegraphics[width=0.95\textwidth]{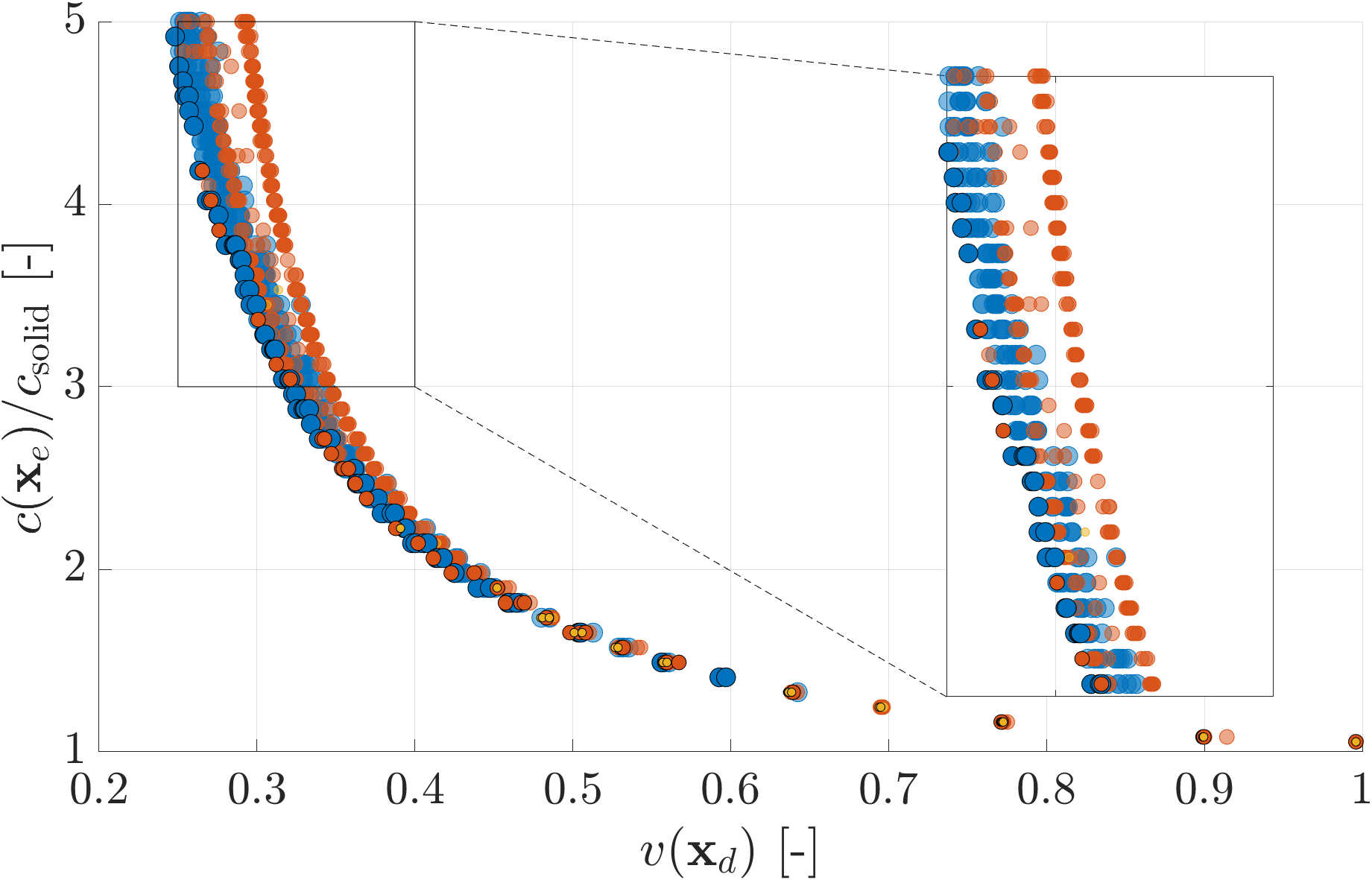}
  \caption{$l_x=3.5$} \label{subfig:asym_350}
\end{subfigure}%
\caption{Pareto frontiers for volume and compliance minimization of an $l_x$ by $1$ cantilever, obtained with $\varepsilon$-constraint minimization ($\min v$ variant, $50$ points). Points with a black edge are nondominated, transparent points without black edge are dominated. Symmetric initial guesses result in symmetric designs and are colored orange. Asymmetric guesses that lead to asymmetric designs are colored blue. Asymmetric guesses that still lead to almost symmetric designs are labeled ``Quasi-symmetric'' and colored yellow.}
\label{fig:Asym_VC1}
\end{figure}

\mc{\Cref{fig:Asym_VC1} shows the points obtained when applying $\varepsilon$-constraint minimization to four cantilevers of varying length $l_x$.\
For each $l_x$, the study is performed $20$ times, each with a different initial guess.\
Half of the initial guesses are asymmetric while the other half are symmetric, including the uniform guess.\
Symmetric and quasi-symmetric points correspond to symmetric designs obtained from symmetric and asymmetric initial guesses, respectively.\ }
\ref{app:ComplianceMinimization} lists numerical details.\\

\mc{Despite the crude sampling methodology and lack of global optimality, two general phenomena can be observed from the frontier approximations in \Cref{fig:Asym_VC1}.\
First, the optimal design's symmetry is influenced by both the cantilever length and the volume.\
For $l_x=2$, all performances are practically nondominated, whether they correspond to symmetric designs or not.\
In contrast, for $l_x=2.5, \, 3$ and $3.5$, the symmetry changes depending on the volume.\
Second, the approximations show clear signs of local frontiers.\
For example, the $l_x=3.5$ cantilever contains a distinct local frontier corresponding to (dominated) symmetric designs.}\\

\mc{A notable observation in \Cref{fig:Asym_VC1} is the discontinuity of the Pareto frontier approximations.\
For example, \Cref{subfig:asym_250} shows a discontinuity near the top of the frontier, where a nondominated local frontier ends before another begins.\
However, the discontinuity in the approximation does not necessarily imply a discontinuity in the underlying Pareto frontier.\
For example, the absence of nondominated asymmetric designs in this region may result from dominated symmetric designs with a large basin of attraction, causing many optimization runs to converge to them~\cite{DeWeer2025}.\
Consequently, the question regarding continuity of the true Pareto frontiers of \Cref{fig:Asym_VC1} remains open.\
Answering this question requires greater control of the path taken by the optimizer and a more exhaustive sampling of the design space.\
This could, for instance, be attempted in a way similar to the continuation approach employed in \Cref{sec:trace_VC}, but is for now deferred to future work.} \\

\subsection{Compliance versus compliance} \label{sec:CC_Canti}

This section studies compliance minimization of the cantilever of \Cref{fig:Canti_setup} under a vertical and compressive load case.\
To limit compliance, an upper bound of $10$ times the solid compliance is placed for both load cases.\
Volume is constrained by $v_{\max}$ to prevent convergence to the trivial solution at $v(\mathbf{x})=100\%$.\
This yields the following formulation.\

\begin{equation}
\begin{aligned}
\min_{\mathbf{x} \in [0,1]^n} \quad &(c_1(\mathbf{x}_e), c_2(\mathbf{x}_e))^\top\\
\text{s.t.} \quad &v(\mathbf{x}_d) \leq v_{\max} \\
&c_1(\mathbf{x}_e) \leq 10 c_1(\mathbf{1}), \, c_2(\mathbf{x}_e) \leq 10 c_2(\mathbf{1})
\end{aligned}
\end{equation}

The Pareto frontier is approximated with 30 points of weighted-sum, $\varepsilon$-constraint and Pascoletti-Serafini scalarization.\
\mc{The extreme points are found by minimizing $c_1(\mathbf{x}_e)$ and $c_2(\mathbf{x}_e)$ under the volume constraint.\
They are used to rescale the objectives and choose the scalarization parameters.\
Details are in \ref{app:ComplianceMinimization}.}\\

\mc{The subsequent discussion first treats the $v_{\max}=50\%$ case and then turns to the more challenging $v_{\max}=30\%$ case.\
For each case, the ``true'' frontier is first examined by combining the solutions of all four scalarizations, after which the scalarization methods are compared.}\\

\begin{figure}[!htb]
\centering
\begin{subfigure}[t]{0.55\textwidth}
  \centering
\captionsetup{width=0.95\textwidth}
  \includegraphics[width=0.95\textwidth]{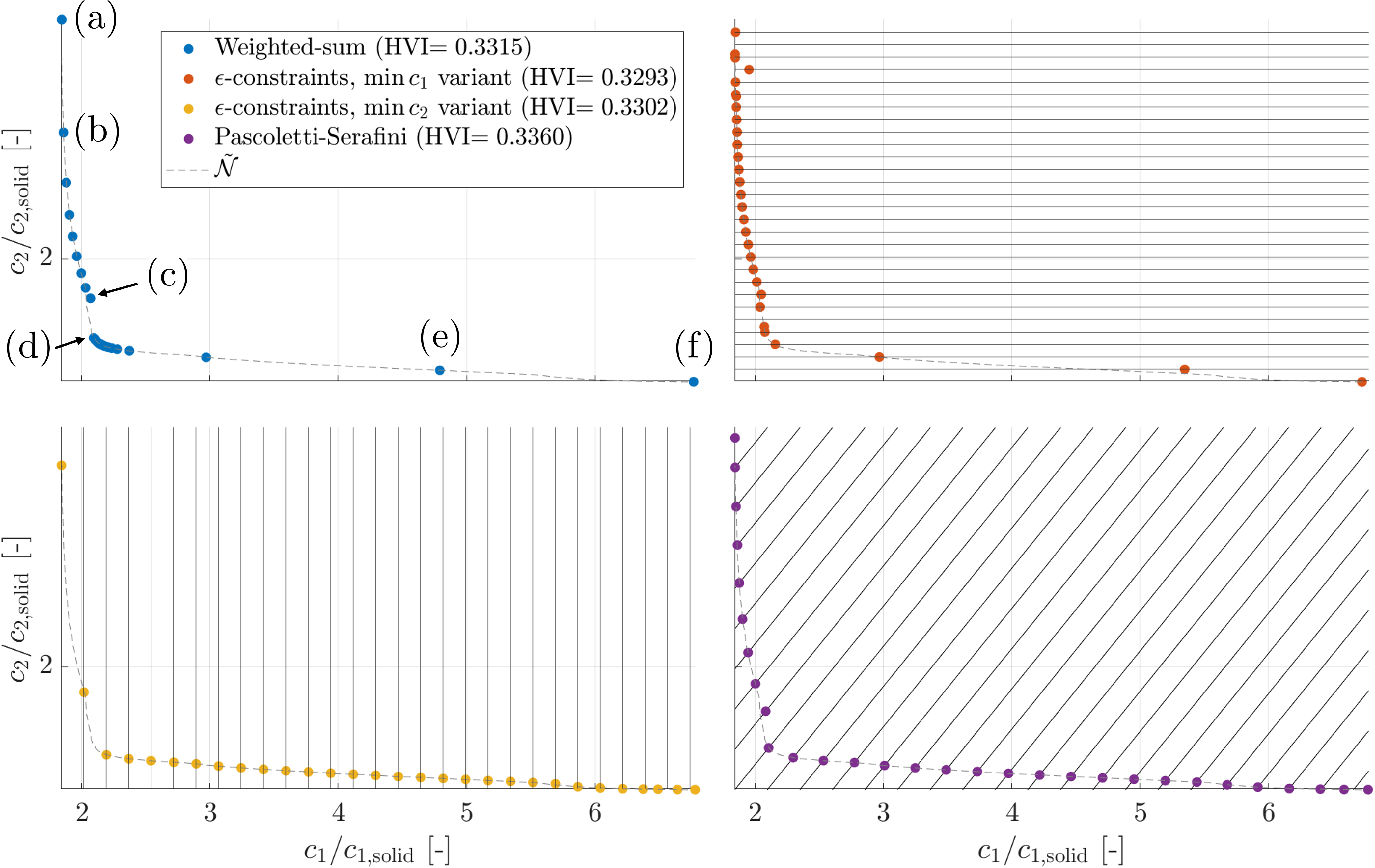}
\end{subfigure}%
\begin{subfigure}[t]{0.45\textwidth}
  \centering
\captionsetup{width=0.95\textwidth}
  \includegraphics[width=0.95\textwidth]{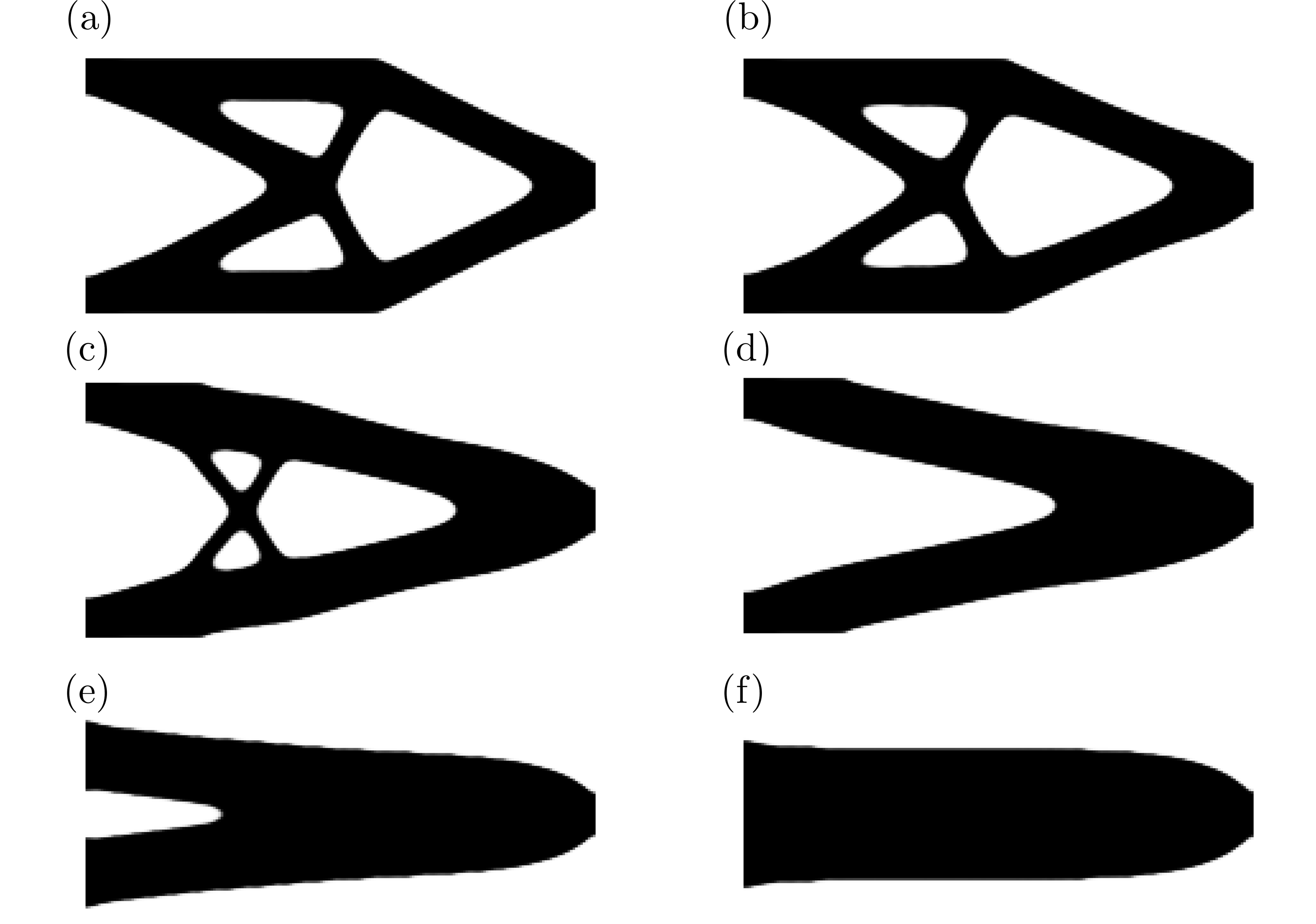}
\end{subfigure}%
\caption{Pareto frontier approximations and selection of some designs for cantilever compliance minimization under two load cases, at $v_{\max}=50\%$. \mc{Due to the changes in topology, the dashed grey line does not always represent a smooth transformation. It merely illustrates the outer hull of the four combined scalarized approximations. The legend includes the hypervolume indicator for each approximation.}}
\label{fig:CC_Canti_V50}
\end{figure}

\mc{\Cref{fig:CC_Canti_V50} shows the $v_{\max}=50\%$ case.\
Between (a) and (c), the design is an eight-bar truss that gradually adopts a more pointed shape to carry the compressive load.\
No change in topology occurs, so solutions (a), (b) and (c) seem to belong to the same local frontier.\
Between (c) and (d), the topology changes: the crosslike infill is removed, yielding a two-bar truss.\
The change in topology corresponds to a change of local frontier, which in turn causes a nonconvexity.\
The remainder of the frontier, from design (d) to (f), sees the two bars gradually merge into a single bar.\
Around (d), the slope of the frontier drops from nearly vertical to nearly horizontal, remaining so during the transition to (e) and (f).\
Note that the continuity of the frontier is unclear in the nonconvex region between (c) and (d).\
Only a limited number of nondominated designs are found between (c) and (d), with varying topologies.\
A further investigation should rely on a continuation of the scalarization parameters, but is here left to future work.}\\ 

\mc{The shape of the frontier affects the scalarization methods in different ways.\
The weighted-sum scalarization produces a cluster of solutions in the high-curvature region around point (d), as this region spans a large range of slopes.\
Additionally, \ref{eq:WS_formulation} leaves a gap in the nonconvex region.\
The $\varepsilon$-constraint scalarizations each provide a good approximation of one ``leg'' of the frontier, but provide poor approximations of the other.\
This is attributed to the high-curvature of the frontier and the uniform distribution of $\varepsilon$ parameters.\
Finally, the Pascoletti-Serafini scalarization provides a uniform approximation due to the choice of $\mathbf{a}$ and $\mathbf{r}$.\
Overall, as in \Cref{sec:WS_EC_PS_VC},  the approximation quality of \ref{eq:PS_formulation} is more robust to the shape of the frontier.\
This is also supported by the hypervolume indicator, which is largest for \ref{eq:PS_formulation}.}\\

\mc{The $v_{\max}=30\%$ case is more challenging than the $50\%$ case because (i) the maximum bound on the first compliance objective is active and (ii) some of the symmetric designs become highly suboptimal.\
Accordingly, the subsequent study includes a KKT-based convergence criterion, detailed in \ref{app:ComplianceMinimization}, and uses a symmetric and asymmetric initial guess.}\\

\begin{figure}[!htb]
\centering
\begin{subfigure}[t]{0.6\textwidth}
  \centering
\captionsetup{width=0.95\textwidth}
  \includegraphics[width=0.95\textwidth]{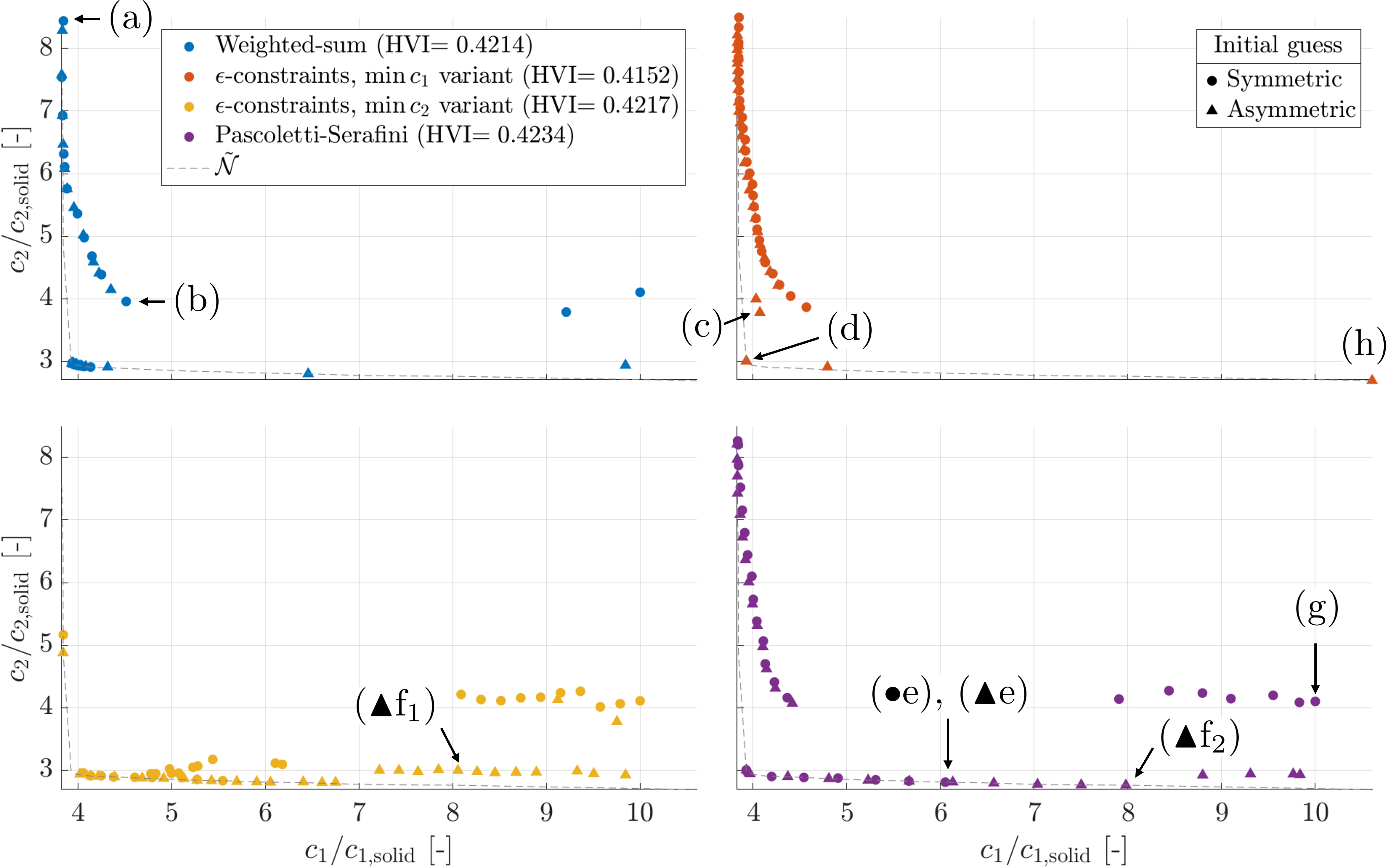}
\end{subfigure}%
\begin{subfigure}[t]{0.4\textwidth}
  \centering
\captionsetup{width=0.95\textwidth}
  \includegraphics[width=0.95\textwidth]{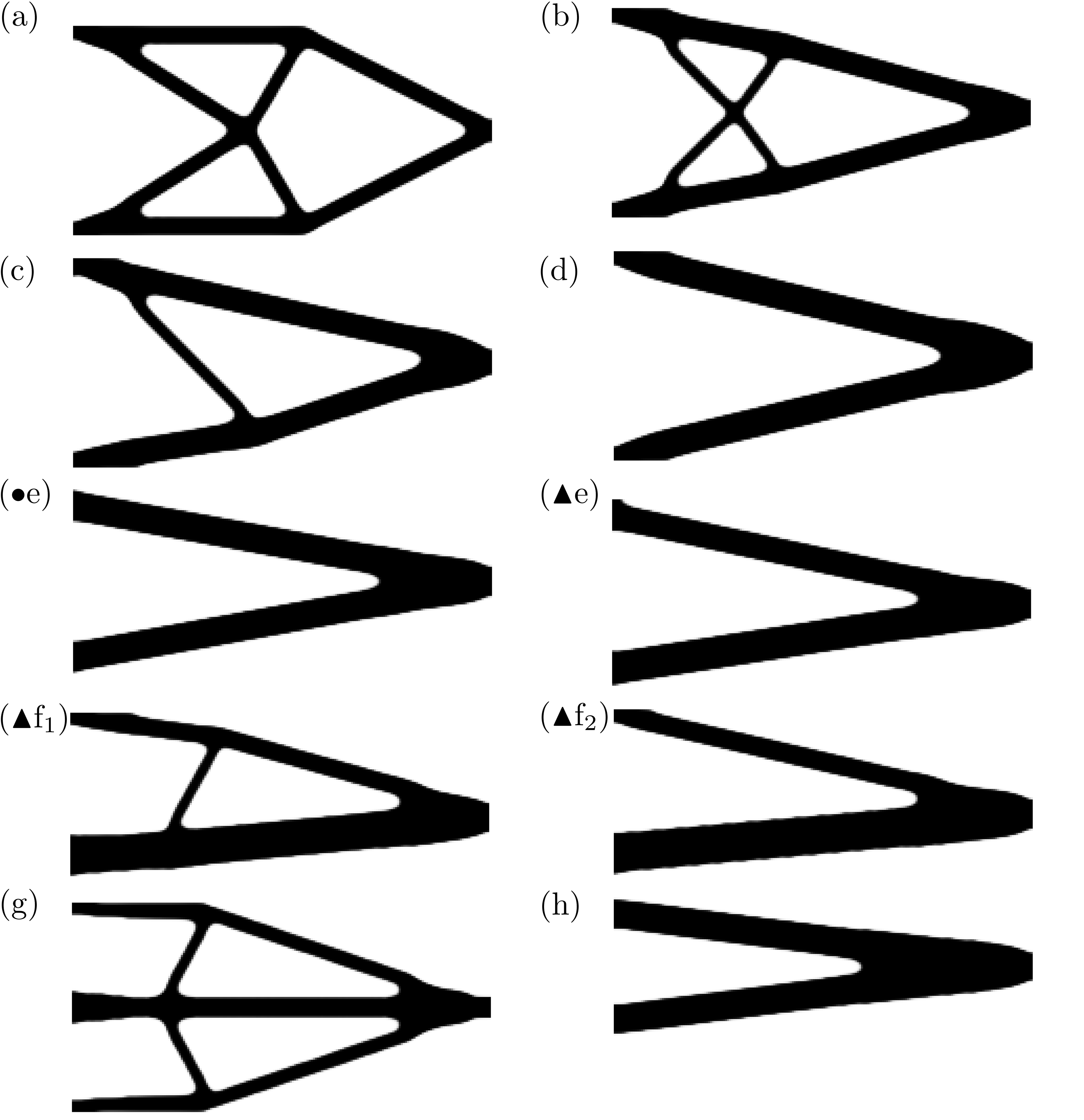}
\end{subfigure}%
\caption{Pareto frontier approximations and selection of some designs for cantilever compliance minimization under two load cases, at $v_{\max}=30\%$. The dashed grey line is again merely illustrative. \mc{The circles and triangles indicate the symmetry of the initial guess. The legend shows the hypervolume indicators HVI for each approximation.}}
\label{fig:CC_Canti_V30}
\end{figure}

\mc{\Cref{fig:CC_Canti_V30} shows the obtained results.\
The transition from design (a) to (b) corresponds to the symmetric eight-bar truss gradually adopting a more pointed shape.\
However, in contrast to the $v_{\max}=50\%$ case, the corresponding local frontier is largely dominated, first by the asymmetric design (c) and then by the two-bar truss (d).\
Then, the transition from (d) to either ($\Huge \bullet$e) or ($\blacktriangle$e) represents a contraction of the two-bar truss in a symmetric and asymmetric way, respectively.\
The asymmetric case can be further extended to designs ($\blacktriangle \mathrm{f}_2$) and eventually $(h)$, although note that the latter violates the compliance bound due to premature convergence.\
In contrast, no nondominated designs are found above $c_1/c_{\mathrm{solid}}=8$: the symmetric designs around (g) represent a local frontier dominated by the local frontier from ($\blacktriangle \mathrm{f}_2$) to (h).}\\

\mc{The approximation quality obtained by each scalarization method is strongly affected by both the shape of the frontier and the initial guess.\
Weighted-sum scalarization again clusters in the high-curvature region and provides only four samples in the ``high'' $c_1$ region, of which two are strongly dominated.\
The $\varepsilon$-constraint variants again provide good approximations of one leg of the frontier, with the $\min c_2$ variant benefiting the most from the asymmetric initial guess to approximate the high $c_1$ region.\
The approximation with \ref{eq:PS_formulation} is qualitatively superior to the others.\
This is also verified quantitatively by computing the hypervolume indicators: \ref{eq:PS_formulation} covers the largest area (HVI$=0.4234$), followed by the $\min c_2$ variant of \ref{eq:EC_formulation} (HVI$=0.4217$), then \ref{eq:WS_formulation} and finally the $\min c_1$ variant of \ref{eq:EC_formulation}.\
However, it is arguably more important to conclude that no scalarization provides a truly satisfactory approximation.\
There remain open questions as to (i) whether the frontier is continuous and (ii) whether the local frontiers really terminate in the locations indicated by the scalarized approximations.\
An approach that follows up the scalarization with a continuation of the local frontiers could shed further light on these questions, but is left to future work.\ }

\subsection{Compliance versus maximum stress} \label{sec:CS_LBeam}

\begin{figure}[!h]
    \centering
    \includegraphics[width=0.4\linewidth]{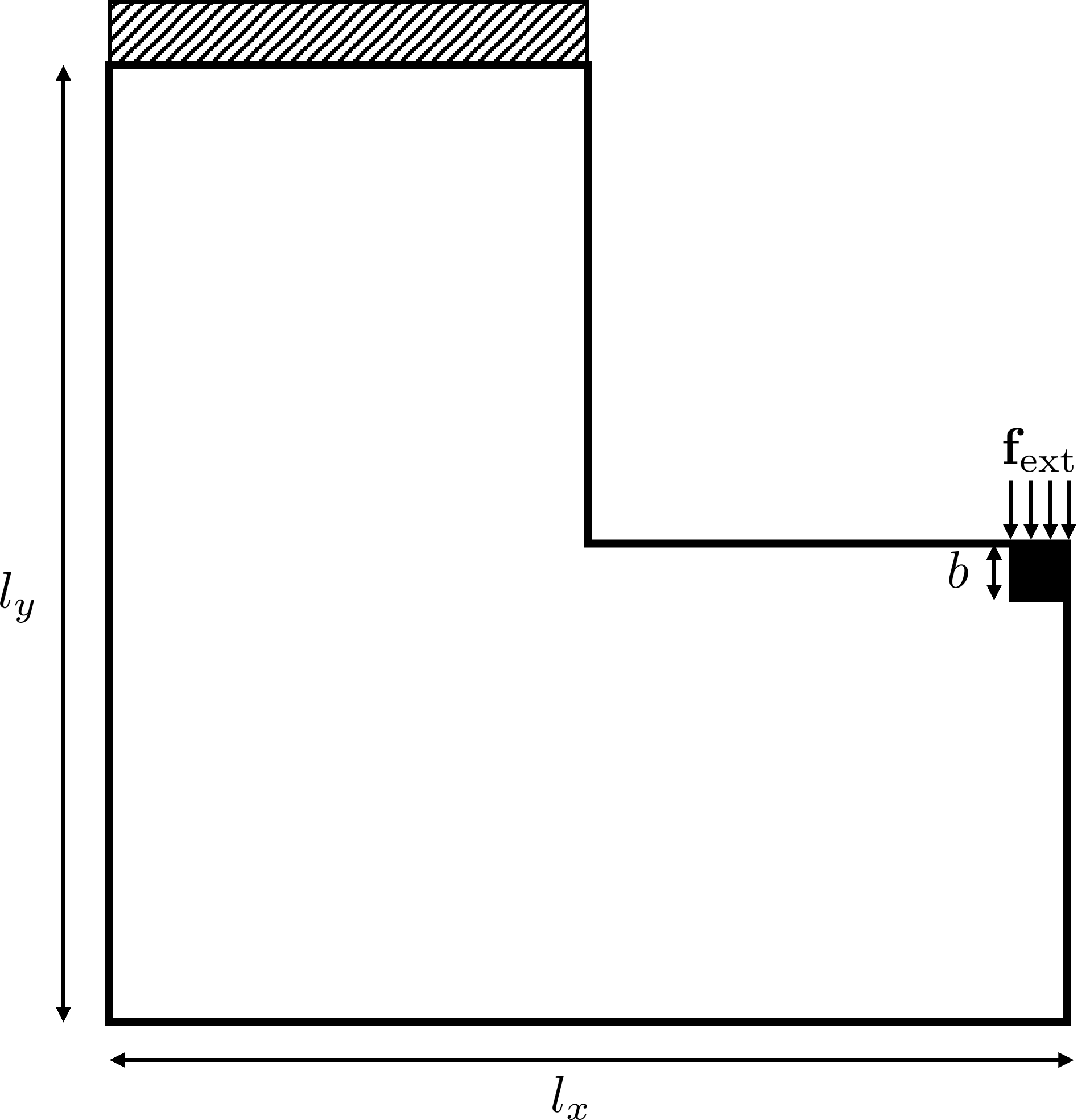}
    \caption{L-shaped problem setup.}
    \label{fig:LBeam}
\end{figure}

To study the effect of more difficult objective functions, this section considers the trade-off between compliance and the maximum stress.\
Consider the L-shaped design domain of \Cref{fig:LBeam}, with $l_x=l_y=1$ and meshed with square elements of size $1/200$ by $1/200$ for a total of $n=\frac{3}{4}200^2=30000$.\
The $b$ by $b$ region consists of $10$ by $10$ solid elements on which a downward force of size $1$ is applied.\ 
The volume fraction is bounded from above by $v_{\max}=50\%$.\
The objectives are the compliance under the applied load and the maximum stress.\\

After application of the load, each element $k$ in the L-shaped design domain has an associated Von Mises stress $\sigma_{k,\mathrm{VM}}$, defined at the element center.\
The elemental stresses are aggregated into the aggregated stress function $s(\mathbf{x}) \approx \max(\{\sigma_{1,\mathrm{VM}}, \ldots, \sigma_{n,\mathrm{VM}} \})$.\
This work employs the function of Kreisselmeier and Steinhauser~\cite{KreisselmeierSteinhauser1979} to provide a differentiable approximation of the maximum stress among the elements.\
Due to the filtering approach, the maximum must also be taken over the eroded, dilated and blueprint design (see \Cref{eq:robustness}).\
To circumvent the nondifferentiability of the $\max$ operator, additional optimization variables are used, with details in \ref{app:SlackVariableTreatment}.\\

The compliance and stress minimization problem is formulated as

\begin{equation}
\begin{aligned}
\min_{\mathbf{x} \in [0,1]^n} \quad &(c(\mathbf{x}_e), \max \left\{ s(\mathbf{x}_e),s(\mathbf{x}_b),s(\mathbf{x}_d) \right\} )^\top \\
\mathrm{s.t.} \, &v(\mathbf{x}_d) < v_{\max}=0.5,
\end{aligned}
\end{equation}
where $\mathbf{x}$ is fixed at $1$ for those entries corresponding to the $100$ elements in the $b$ by $b$ solid region.\
\mc{As before, separate compliance and stress minimization runs are performed to obtain $(c_{\min},s^\star)$ and $(c^\star, s_{\min})$, respectively, to both rescale the objectives and select the values of $\mathbf{a}$, $\mathbf{r}$ and $\boldsymbol{\varepsilon}$.}\
Continuation and material parameters are in \ref{app:StressMinimization}.\\

\begin{figure}[!htb]
\centering
\begin{subfigure}[t]{0.4\textwidth}
  \centering
\captionsetup{width=0.95\textwidth}
  \includegraphics[width=0.95\textwidth]{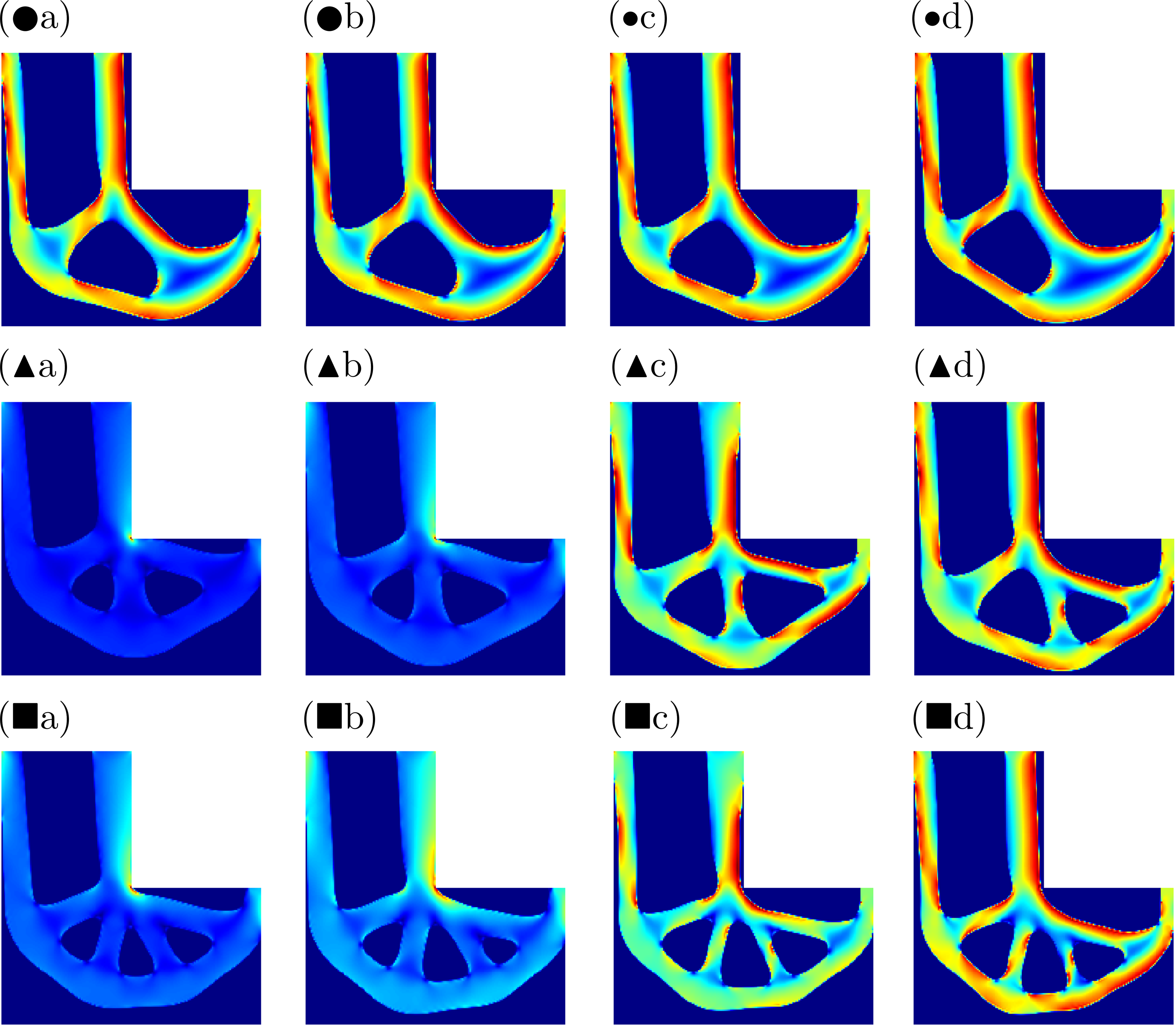}
  \caption{Selected designs, with separate colormaps indicating the relaxed Von Mises stress (see $\ref{app:StressMinimization}$).}
\end{subfigure}%
\begin{subfigure}[t]{0.6\textwidth}
  \centering
\captionsetup{width=0.95\textwidth}
  \includegraphics[width=0.95\textwidth]{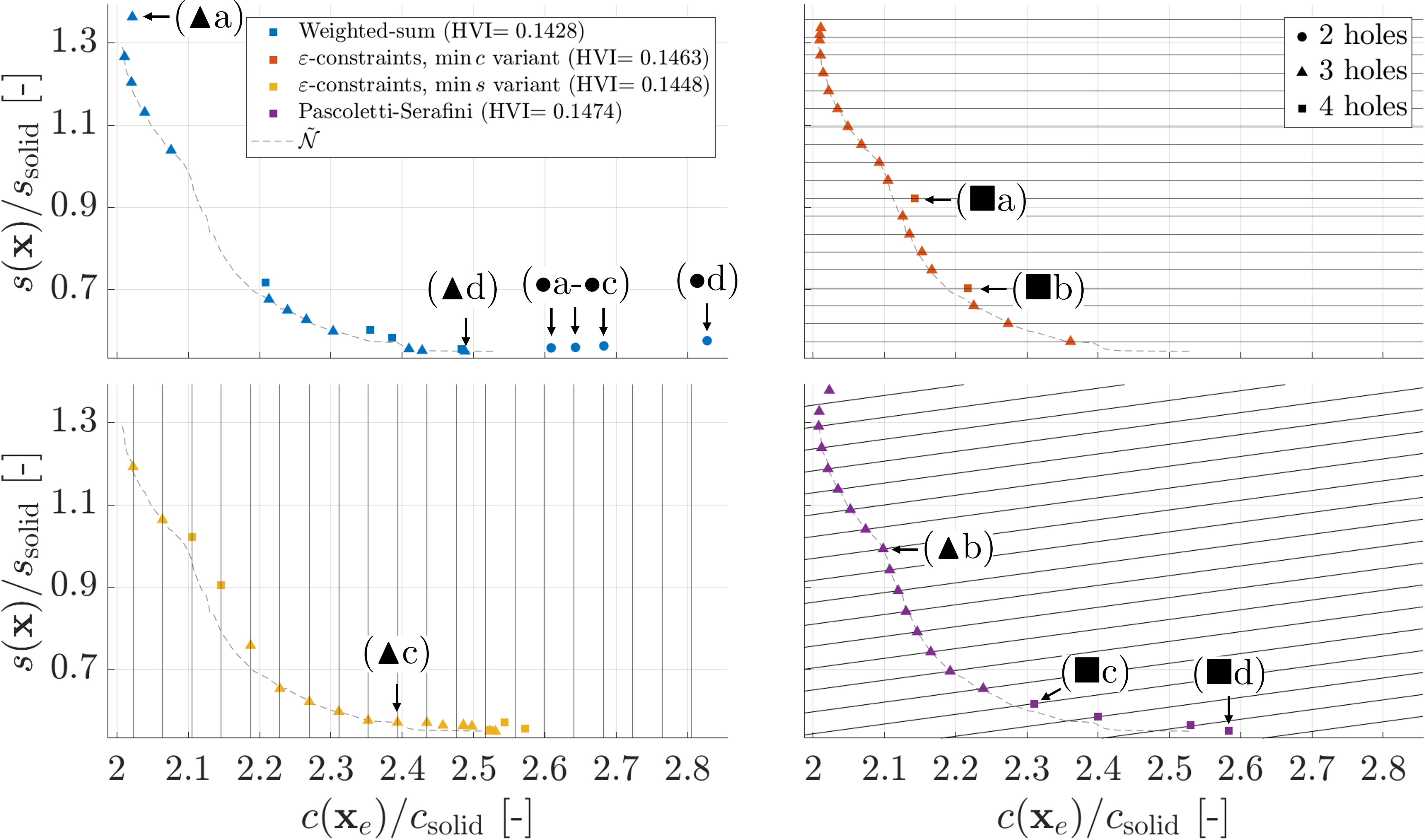}
  \caption{Nondominated sets for each scalarization. Black lines denote values of constant $\varepsilon$ or $\mathbf{p} = \left\{ \mathbf{a}+t\mathbf{r}  \,|\, t \in \mathbb{R}\right\}$ to illustrate the scalarization parameters.}
\end{subfigure}%
\caption{Maximum stress and compliance minimization results with weighted-sum, $\varepsilon$-constraints ($\min s(\mathbf{x})$ and $\min c(\mathbf{x})$ variants) and Pascoletti-Serafini scalarization. \mc{The legend mentions the hypervolume indicator (HVI). Three different topologies are obtained: (i) designs with 2 holes are marked with a circle, (ii) designs with 3 holes are marked with a triangle and (iii) designs with 4 holes are marked with a square. Due to these different topologies, the dashed grey line does not correspond to a smooth design space transformation. Nonetheless, the portions connecting two designs with the same topology likely represent a local frontier.}}
\label{fig:CS_LBeam}
\end{figure}

\mc{\Cref{fig:CS_LBeam} shows the Pareto frontier approximation obtained with twenty points of the weighted-sum, $\varepsilon$-constraint and Pascoletti-Serafini scalarization.\
The following discussion first treats the general characteristics of the  designs and frontiers before comparing the scalarization methods.}\\

\mc{The obtained designs can be grouped into three topologies: those with two, three and four holes.\
The designs with two holes are only found by the weighted-sum scalarization and represent the highest compliance but are all dominated.\
Furthermore, they also dominate each other: going from design $(\Huge \bullet \mathrm{a})$ to $(\Huge \bullet \mathrm{d})$, the inner bar moves upward which increases both the stress \emph{and} the compliance.\ }
 
\mc{The designs with three holes are generally superior.\
Moving from design $(\blacktriangle \mathrm{a})$ to $(\blacktriangle \mathrm{d})$, the bars lengthen and detach from the high-stress corner of the L-shaped beam.\
Interestingly, design $(\blacktriangle \mathrm{b})$ appears to lie at the crossing of two convex regions of the two-hole designs' local frontier.\
This means the three-hole local frontier is nonconvex and nonsmooth at $(\blacktriangle \mathrm{b})$.\ }

\mc{The designs with four holes are mostly dominated by the three-hole designs.\
Moving from $(\blacksquare \mathrm{a})$ to $(\blacksquare \mathrm{d})$, the same progression is observed as the three-hole designs: the bars lengthen and detach from the high-stress corner of the L-shaped beam.\
Near the lower stress region of the frontier, the three- and four-hole designs show about equal performance.} \\

\mc{Overall, the three-hole designs are dominant or close to dominant in every region of the Pareto frontier.\ 
It therefore seems like the frontier is dominated by a single topology, which implies a single local frontier.\
Nonetheless, the local frontier corresponding to the three-hole (and four-hole) designs is nonconvex: when the designs detach from the high-stress corner (at point $\blacktriangle \mathrm{b}$), the Pareto frontier has an inflection.\
This affects the approximation quality of the scalarization methods, which are now compared.\
}

\mc{The weighted-sum scalarization produces an approximation with a large gap, which is attributed to the nonconvexity of the frontier in this region.\
Additionally, \ref{eq:WS_formulation} produces several dominated designs near the extreme points of the frontier: four low-stress two-hole designs and one low-compliance three-hole design.\
Of these, points $\blacktriangle \mathrm{a}$ and $\Huge \bullet \mathrm{d}$ correspond to $(c_{\min},s^\star)$ and $(c^\star, s_{\min})$, respectively, which are used to rescale the objectives and select values for $\mathbf{a}$, $\mathbf{r}$ and $\boldsymbol{\varepsilon}$.\
This affects the approximation quality of both \ref{eq:EC_formulation} and \ref{eq:PS_formulation}.\
In particular, the $\min c(\mathbf{x})$ variant of the $\varepsilon$-constraint scalarization does not cover the low-stress region and clusters near the low-compliance region due to inactive constraints.\
Conversely, the $\min s(\mathbf{x})$ variant clusters in the low-stress region and does not cover the low-compliance region.\
Finally, the Pascoletti-Serafini scalarization produces the best Pareto frontier approximation: it covers the complete frontier and produces only two dominated points.\
This is also reflected by the hypervolume metric, which indicates that \ref{eq:PS_formulation} covers the largest area (HVI$=0.1474$), followed by the \ref{eq:EC_formulation} variants (HVI$=0.1463$ and $0.1448$) and then \ref{eq:WS_formulation} (HVI$=0.1428$).\
The exceptional robustness of \ref{eq:PS_formulation} is attributed to $\mathbf{a}$ and $\mathbf{r}$ lying uniformly on and perpendicular to the CHIM, respectively, and the fact that the CHIM mirrors the slope of the Pareto frontier.\ }

\subsection{Compliance versus dynamic compliance} \label{sec:CD_Canti}

This section applies weighted-sum, $\varepsilon$-constraints and Pascoletti-Serafini scalarization to the minimization of compliance $c(\mathbf{x})$ and dynamic compliance $d(\mathbf{x})$ of a $2$-by-$1$ cantilever with volume fraction $v(\mathbf{x}) \leq 0.5$.\
\mc{The most significant methodological difference with the previous examples is the different rescaling of the dynamic compliance objective.\
After obtaining $(c_{\min},d^{\star})$ and $(c^\star,d_{\min})$, it was noticed that $d_{\min} \approx d^{\star}$.\
This strongly affects the scalarization methods.\
For example, for \ref{eq:PS_formulation} the line connecting $(c_{\min},d^{\star})$ and $(c^{\star},d_{\min})$ is almost horizontal and so $\mathbf{r}$ is almost vertical.\
As a result, \ref{eq:PS_formulation} behaves almost as the $\min d(\mathbf{x})$ variant of the $\varepsilon$-constraint scalarization.\
To provide a better comparison, this section therefore normalizes the dynamic compliance with $d^S (\mathbf{x}) = \frac{d(\mathbf{x})-d_{\min}}{d^\star - d_{\min}}$ and adapts the scalarization methods accordingly.\ }
For brevity, additional numerical details are moved to \ref{app:DynamicComplianceMinimization}.\\

\begin{figure}[!htb]
\centering
\begin{subfigure}[t]{0.65\textwidth}
  \centering
\captionsetup{width=0.95\textwidth}
  \includegraphics[width=0.95\textwidth]{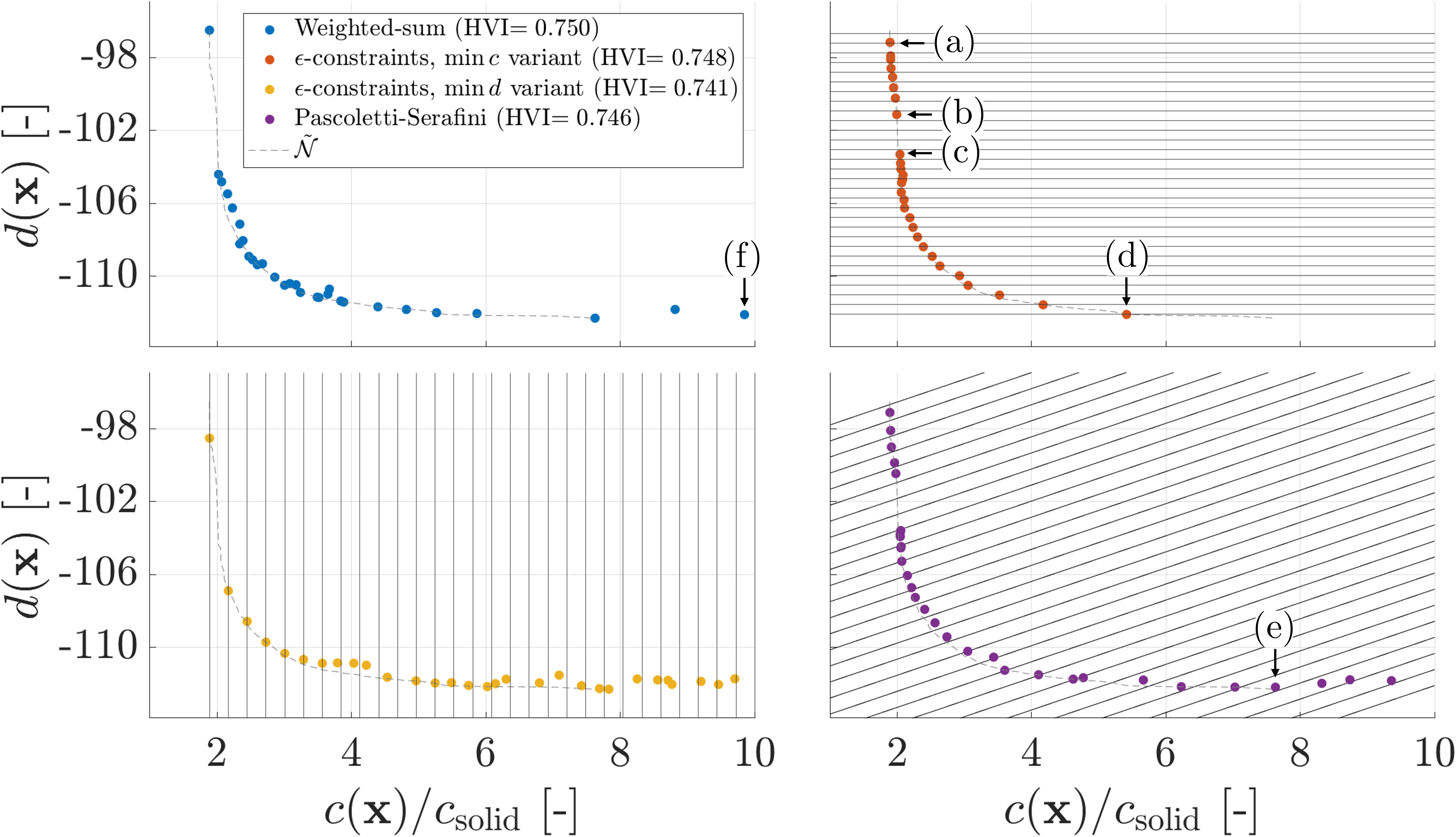}
  \caption{Nondominated set approximations. The grey dashed line $\tilde{\mathcal{N}}$ denotes the combined approximation, but does not imply intermediate points exist. Black lines denote values of constant $\varepsilon$ or the set $\left\{ \mathbf{a}+t\mathbf{r}  \,|\, t \in \mathbb{R}\right\}$. \mc{The legend shows the hypervolume indicators (HVI).}}
  \label{subfig:CD_Canti_A}
\end{subfigure}%
\begin{subfigure}[t]{0.35\textwidth}
  \centering
\captionsetup{width=0.95\textwidth}
  \includegraphics[width=0.95\textwidth]{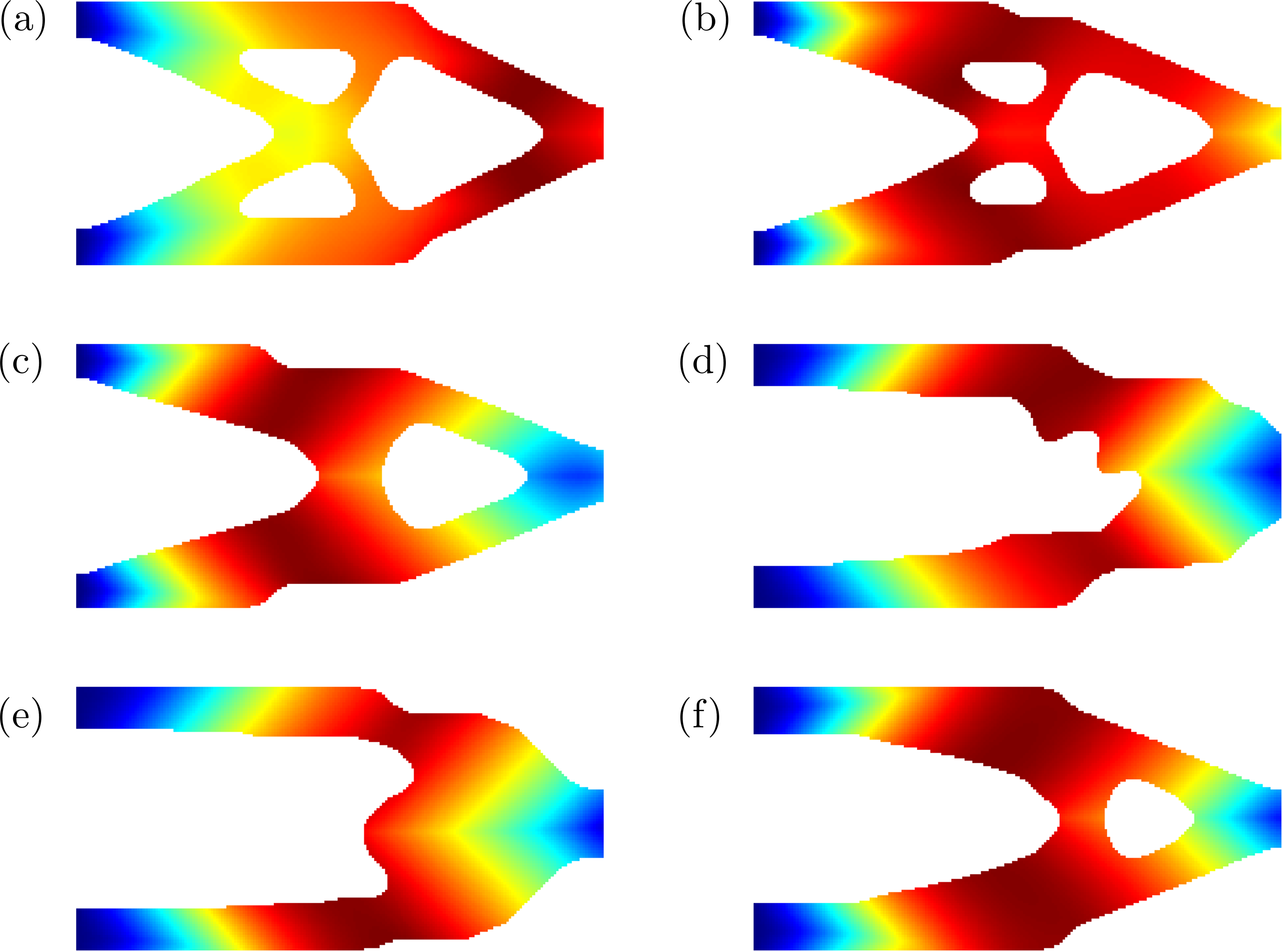}
    \caption{\mc{Corresponding designs. Colorscale indicates the relative magnitude of the dynamic displacement at the loading frequency.}}
\end{subfigure}%
\caption{Static and dynamic compliance minimization results with weighted-sum, $\varepsilon$-constraint ($\min c(\mathbf{x})$ and $\min d(\mathbf{x})$ variants) and Pascoletti-Serafini scalarization.}
\label{fig:CD_Canti}
\end{figure}

\mc{\Cref{fig:CD_Canti} shows the obtained results.\
The following discussion first treats the characteristics of the frontier and designs and then compares the scalarization methods.}\\

\mc{The frontier consists of two separate parts.\
The first, from point (a) to (b), corresponds to statically stiff designs with high dynamic compliance.\
The second part, from (c) to (d) and then (e), corresponds to statically compliant but dynamically stiff designs.\
Between points (d) and (e), the frontier is especially flat.\
The low dynamic compliance in this region stems from an antiresonance placed at the loading frequency.\
The region from (e) to (f) is not part of the nondominated set, as these are all dominated by point (e) except (e) itself.\
Finally, between points (b) and (c) a gap occurs for all scalarization methods.\
This could be attributed to either (a) the lack of feasible points in this region, or (b) the existence of dynamically stiff local optima with a large basin of attraction~\cite{DeWeer2025}.}\\

\mc{Based on the above, the scalarization methods are now compared.\
The weighted-sum scalarization clusters most points around  the bend of the frontier, between (c) and (d) in \Cref{subfig:CD_Canti_A}.\
This is in line with previous observations and again linked to the slope-based interpretation of the weighted-sum scalarization.\
A non-uniform choice of the weights could alleviate the clustered spacing.\ }

\mc{The $\varepsilon$-constraint scalarization variants each provide a good approximation of one leg of the frontier but undersample the other.\
This is also in line with previous observations.\
The static compliance minimization variant ($\min c(\mathbf{x})$) clusters around (c), where many of the optimization runs do not have an active dynamic compliance constraint.\
Furthermore, since the dominated point (f) is used to choose the $d(\mathbf{x}) \leq \varepsilon_d$ constraint bounds, the approximation does not cover point (d) to (e).\ 
The dynamic compliance minimization variant (``$\min d(\mathbf{x})$'') yields an irregular spacing for the low dynamic-compliance region, possibly due to the many local optima and the flatness of this part of the frontier.\
Finally, the Pascoletti-Serafini strikes a balance between the two $\varepsilon$-constraint variants.\
It also clusters around (c), although the formulation of \ref{eq:PS_formulation} ensures at least one of the $\mathbf{a}+t\mathbf{r}$ constraints is active.\
Overall, the Pascoletti-Serafini scalarization is most robust against the shape of the frontier.} \\

\mc{The hypervolume indicators show that \ref{eq:WS_formulation} covers the largest area, closely followed by $\min c$ variant and then the \ref{eq:PS_formulation} formulation.\
It is however unclear which factor contributes the most to this result.\
The $\min d$ variant covers the lowest area, which is likely due to the fact that a large amount of its points are dominated.\ }

\section{Conclusions} \label{sec:Conclusions}

This work studies topology optimization problems with multiple objective functions.\
Despite a long coexistence, a literature overview (\Cref{sec:MOTO}) shows little cross-fertilization and a widening gap between the topology and multiobjective optimization fields.\
This work aims to address this gap by studying topology optimization problems involving compliance, volume, maximum stress, dynamic compliance.\
For several bi-objective numerical examples, the nondominated set is approximated with the weighted-sum (\ref{eq:WS_formulation}), $\varepsilon$-constraint (\ref{eq:EC_formulation}) and Pascoletti-Serafini (\ref{eq:PS_formulation}) scalarization.\\

The first contribution of this paper (\Cref{sec:MOO}) is to refer to theoretical results of the field of multiobjective optimization regarding the chosen scalarizations.\
Although weighted-sum and $\varepsilon$-constraint scalarization are widely used in the topology optimization community, theory shows that (i) weighted-sum cannot guarantee a complete recovery of the front and (ii) $\varepsilon$-constraint scalarization may lead to infeasible optimization problems.\
In contrast, Pascoletti-Serafini scalarization has all the theoretical properties a scalarization method should have: a suitable variation of $(\mathbf{a},\mathbf{r}) \in \mathbb{R}^m \times \mathbb{R}^m_{++}$ yields all nondominated points and each choice yields a weakly efficient solution.\\

The second contribution is a novel insight into the structure of the Pareto front in topology optimization.\
\mc{As illustrated for a truss example in \Cref{sec:truss_opt}, local optima in single-objective topology optimization extend to local Pareto fronts in multiobjective topology optimization.\
Each local front is associated to locally optimal designs with a fixed topology, since topological changes are inherently discontinuous: under an imposed length scale, holes must appear and disappear suddenly.\
Then, just as the global solutions form a subset of the local solutions, the global Pareto front is a subset of the local fronts.\
Consequently, the Pareto front consists of local front segments.\
In contrast to the smooth, often convex fronts reported in the literature, this insight naturally admits nonconvexities and discontinuities.}\\

\mc{The numerical examples demonstrate how local fronts shape the global front.\ 
For volume and compliance minimization, the front is continuous but not smooth: when two local fronts cross, the slope of the front changes discontinuously (see \Cref{sec:trace_VC}).\
For minimization of compliance under two load cases, a discontinuity is observed at low volume fractions: the sudden termination of a dominant local front causes a gap in the Pareto front.\
Therefore, even though the local fronts in volume and compliance minimization are observed to be convex, their composition does not necessarily yield a convex global front.\
Next, the stress minimization example of \Cref{sec:CS_LBeam} shows a single dominant local front.\
Nonetheless, the front is nonconvex due to an inflection when the designs detach from the sharp corner of the L-shaped design domain.\
Finally, for the minimization of static and dynamic compliance in \Cref{sec:CD_Canti}, the Pareto front contains a gap that may arise either from the sudden termination of a local front or from local optima with a large basin of attraction.\
Note however, that the analysis is complicated by the inherent nonconvexity of topology optimization, which prohibits certainty regarding global optimality.\
Further research, e.g, with continuation and design space sampling, can be done to draw firmer conclusions and further explore the local fronts that make up the Pareto front.}\\

\mc{The third contribution is a numerical comparison of the three chosen scalarization methods, and how they are affected by the a priori unknown shape of the Pareto front.\
All scalarization methods are sensitive to the range of the objectives and hence rescaled by obtaining approximate values for the extent of the front.\
The choice of scalarization method had little influence on convergence speed but a pronounced effect on approximation quality was observed.}\\

\mc{The numerical examples show that the approximation quality of the weighted-sum scalarization is strongly dependent on the shape of the Pareto front.\
When it consists of more than one local front segment, each segment exhibits a different sample density, resulting in highly irregular spacing.\
When the front is nonconvex, \ref{eq:WS_formulation} leaves gaps.\
When the front has regions with high curvature, it oversamples these regions and neglects the remainder.\
A uniformly-sampled weighted-sum scalarization is therefore best suited for problems with a convex front, or as an initial exploratory tool.\
Increasing the number of samples cannot resolve the fundamental limitations of \ref{eq:WS_formulation}.\
Instead, approximation quality is more effectively improved by changing to other scalarization methods, such as \ref{eq:EC_formulation} and \ref{eq:PS_formulation}.}\\

\mc{The $\varepsilon$-constraint scalarization uses a uniform $\varepsilon$ distribution, which generally leads to more regularly spaced solutions.\ 
The scalarization is sensitive to inaccuracies in the estimated individual minima.\
Since these values determine distribution of the $\varepsilon$-parameters, inaccurate estimates can lead to either a loss of coverage of the front or a clustering due to inactive constraints.\
Additionally, for high-curvature front, one leg is oversampled while the other is undersampled.\
However, in contrast to \ref{eq:WS_formulation}, \ref{eq:EC_formulation} can approximate the nonconvex region.\
An adaptive $\varepsilon$-selection could further improve the approximation quality~\cite{Eichfelder2008}.}\\

\mc{Of the considered scalarization methods, the Pascoletti-Serafini scalarization is the least sensitive to the shape of the Pareto front.\
It can provide an almost uniform distribution of solutions for Pareto fronts with nonconvexities and high curvature.\
This robustness is attributed to its flexible formulation, which behaves as a weighted average of two $\varepsilon$-constraints.\
Nevertheless, even \ref{eq:PS_formulation} cannot fully overcome the challenges posed by the local fronts.}\\

\mc{In addition to the scalarization methods, the numerical experiments employ design space sampling (\Cref{sec:asym_VC_Canti,sec:CC_Canti}) and continuation of the scalarization parameters (\Cref{sec:trace_VC}) as exploratory tools.\
To obtain a more complete picture of the Pareto front, future research should shift from approximating isolated Pareto-optimal points to systematically extending local optima into local fronts, revealing how they collectively form the global Pareto front.}

\section*{Statements and Declarations}

\subsection*{Acknowledgements}

\mc{The authors would like to express their gratitude to the referees for their valuable comments and suggestions, which greatly improved the manuscript.\ }

\subsection*{Funding}
Resources and services were provided by the VSC (Flemish Supercomputer Center), funded by the Research Foundation - Flanders (FWO) and the Flemish Government.\
The research stay of T. De Weer was funded by FWO grant V456025N.\
Financial support from the Villum Foundation through the Villum Investigator Project AMSTRAD (VIL54487) is gratefully acknowledged.\

\subsection*{Conflict of Interest}
On behalf of all authors, the corresponding author states that there is no conflict of interest.

\subsection*{Author Contributions}
All authors contributed to the study conceptualization and design, commented on previous versions and read and approved the final version.\
Data collection, analysis and writing was performed by the corresponding author.\

\subsection*{Ethics approval and Consent to participate}

Not applicable.

\subsection*{Data Availability and Replication of Results}

Methodological details are described to facilitate the replication of the results.\
Data will be made available on request.

\appendix

\section{$\max$ function treatment and interface with the Method of Moving Asymptotes} \label{app:SlackVariableTreatment}

The numerical examples consider two types of objective functions.\
Type A objectives are differentiable and denoted as $f_i^A$, $i=1,\ldots,m_A$.\
They are dealt with in the usual way.\
Type B objectives contain a $\max$ operator and are hence not differentiable.\
They are denoted as $f_j^B=\max \left\{ f_{j,e},f_{j,b},f_{j,d} \right\}$, $j=1,\ldots,m_B$, where the $e$, $b$ and $d$ subscripts denote eroded, blueprint and dilated.\
\mc{Note that both objectives correspond to the rescaled objectives $f_i^S$ discussed at the beginning of \Cref{sec:numerical_examples}.\ }

The weighted-sum and $\varepsilon$-constraint scalarizations are combined in the following optimization problem.

\begin{equation} \label{eq:WS_EC}
\begin{aligned}
\min_{\mathbf{x} \in X, z_j} \quad &\sum_{i=1}^{m_A} w_i^A f_i^A(\mathbf{x}) +  \sum_{j=1}^{m_B} w_j^B z_j  \\
\text{s.t.} \quad &f_i^A(\mathbf{x}) \leq \varepsilon_i^A, \quad i=1,\ldots,m_A \\
&z_j^B \leq \varepsilon_j^B, \quad j=1,\ldots,m_B \\
&f_{j,e}(\mathbf{x}) \leq z_j, f_{j,b}(\mathbf{x}) \leq z_j, f_{j,d}(\mathbf{x}) \leq z_j  \quad j=1,\ldots,m_B \\
\end{aligned}
\end{equation}

The type B objectives are thus replaced by additional variables $z_j$ and the $\max$ operator is enforced via the last constraint.\
Weighted-sum scalarizations omits the first two constraints and varies $w_i^A$ and $w_j^B$ whereas $\varepsilon$-constraint scalarization sets one of the weights to $1$ and varies the values of $\varepsilon_j^A$ and $\varepsilon_j^B$.\

Pascoletti-Serafini scalarization is implemented via the following formulation.

\begin{equation} \label{eq:PS_EC}
\begin{aligned}
\min_{\mathbf{x} \in X, z_j, t} \quad &t  \\
\text{s.t.} \quad & f_i^A(\mathbf{x}) - a_i^A - r_i^A t \leq 0, \quad i=1,\ldots,m_A \\
&z_j - a_j^B - r_j^B t \leq 0, \quad j=1,\ldots,m_B \\
&f_{j,e}(\mathbf{x}) \leq z_j, f_{j,b}(\mathbf{x}) \leq z_j, f_{j,d}(\mathbf{x}) \leq z_j  \quad j=1,\ldots,m_B \\
\end{aligned}
\end{equation}

The first two constraints enforce $\mathbf{a}+t\mathbf{r} - \mathbf{f}(\mathbf{x}) \geq 0$, where $ a_i^A, r_i^A$ and $ a_j^B, r_j^B$ are entries of $\mathbf{a}$ and $\mathbf{r}$ corresponding to the type A and B objectives, respectively.\
The last constraint again enforces the $\max$ operator.\

The interface with the Method of Moving Asymptotes~\cite{Svanberg1987} is straightforward for \Cref{eq:WS_EC}: the additional variables are added to the optimization variable vector and $c_i=1000$, $d_i=1$.\
For numerical stability, $a_0$ is set to 1 and all $a_i$ are set to 0, as suggested by Svanberg~\cite{Svanberg1987}.\
The asymptote control parameters are $s_{\mathrm{init}} = 0.2$,  $s_{\mathrm{decr}} = 0.65$,  $s_{\mathrm{incr}} = 1.07$.\
All examples employ an outer move limit with $\mu = 0.2$, except the dynamic compliance minimization routines which use $\mu = 0.1$.\

For \Cref{eq:PS_EC}, MMA's internal ``artificial'' variable $z$ is used to represent $t$ whereas the additional variables $z_j$ are treated as for \Cref{eq:WS_EC}.\
Thus, $a_0$ is again $1$ but all $a_i$ are set to the values of the corresponding entries in $\mathbf{r}$.\
Note that MMA's internal variable $z$ is bound to be positive whereas $t$ can become negative.\
To circumvent this, Pascoletti-Serafini's $\mathbf{a}$ vector is updated to $\mathbf{a}-\mathbf{r}$ whenever $t < 0.1$.\ 
All other MMA parameters ($c_i$, $d_i$, $s_{\mathrm{init}}$, $\ldots$) are the same as for \Cref{eq:WS_EC}.\

\mc{MMA requires properly scaled objectives and constraints.\
For the objectives, this is already satisfied by the use of the rescaled objectives $f_i^S$.\
The $\varepsilon$ constraints are scaled by dividing with the $\varepsilon$ value.\ 
For example, the $f_i^A(\mathbf{x}) \leq \varepsilon_i^A$ constraints in \Cref{eq:WS_EC} are written as $f_i^A(\mathbf{x}) / \varepsilon_i^A -1 \leq 0$.\ }

\section{Compliance minimization parameters}\label{app:ComplianceMinimization}

All topology optimization examples in \Cref{sec:VC_Canti,sec:CC_Canti} use isotropic material with a Young modulus of $1\,\mathrm{Pa}$ and a Poisson ratio of $0.3$.\
SIMP penalization is used with $P$ going from $1$ to $3$ in $N_P$ steps.\
The minimum Young modulus $\delta$ is $1\mathrm{e}{-6}$.\
Whenever $P$ is increased, the projection parameter $\beta$ also increases, going from $1$ to $\beta_{\mathrm{lim}}=2R$ in $N_{\beta}=N_P$ steps~\cite{daSilva2019}.\
Here, $R$ is the filter parameter of the PDE filter with consistent boundary conditions~\cite{Lazarov2011,Wallin2020}.\
The robust approach is used to enforce a minimum length scale~\cite{wang2011projection,qian2013topological}, where $\eta_b=0.5$ and $\eta_{e,d}=\eta_b \pm 0.1$.\
The stopping criterion, used for both termination and continuation, is $5$ consecutive iterations with a relative objective change below $\mathrm{rtol}$ and nonviolated constraints.\
The $b$-by-$b$ square of solid elements in the cantilever is kept at size $b=10$ elements.\
Unless specified otherwise, the initial guess is a uniform density $x=0.5$ at every element.\
Remaining parameters are listed in \Cref{tab:compliance_minimization_parameters} for numerical example.\
\mc{The $\mathrm{rtol}$ parameter for the final row, corresponding to the compliance minimization results of \Cref{sec:CC_Canti}, is $1\mathrm{e}{-5}$ for the $v_{\max}=0.5$ case.\
Due to the challenging nature of the $v_{\max}=0.3$ case, the $\mathrm{rtol}$-based relative objective-change stopping criterion is replaced by a KKT-based criterion.\
Specifically, convergence is declared when the KKT residual remains below $1.e-4$ for five consecutive iterations, or below $1.e-3$ for twenty-five consecutive iterations.\ }

\begin{table}[h!]
\caption{Parameters for examples involving compliance and volume minimization. Filter radius $R$ is expressed in element widths.}
\label{tab:compliance_minimization_parameters}
\begin{tabular}{l c c c c c}
\toprule
\textbf{Section} & $N_P$ & $\mathrm{rtol}$ & $R$ &  $n_x$ & $n_y$\\
\midrule
\Cref{sec:WS100} & $20$ & $1\mathrm{e}{-4}$ & $5$, $10$, $20$ &  $400$ & $200$ \\[0.5ex]
\Cref{sec:trace_VC} & n.a. & $1\mathrm{e}{-4}$ & $20$ &  $400$ & $200$ \\[0.5ex]
\Cref{sec:WS_EC_PS_VC} & $20$ & $1\mathrm{e}{-4}$ & $20$ &  $400$ & $200$ \\[0.5ex]
\Cref{sec:asym_VC_Canti} & $20$ & $1\mathrm{e}{-4}$ & $10$ &  $200,\, 250,\, 300,\, 350$ & $100$ \\[0.5ex]
\Cref{sec:CC_Canti} & $50$ & $1\mathrm{e}{-5}$ ($v_{\max}=0.5$) & $10$ &  $200$ & $100$ \\[0.5ex]
\bottomrule
\end{tabular}
\end{table}

The asymmetry study of \Cref{sec:asym_VC_Canti} classifies designs based on the asymmetry metric $\left\lVert \mathbf{x}_b - \mathbf{x}_{b}^\star \right\rVert_1 / \left\lVert \mathbf{x}\right\rVert_1$, where $\mathbf{x}_b^\star$ contains the densities of the blueprint design flipped around the $\mathrm{x}$ axis.\ 
In this way, a fully symmetric design has an asymmetry metric of zero and a fully asymmetric design (solid above symmetry plane, void below) has a value of one.\
Optimization runs with a symmetric initial guess always return designs with a symmetry index equal to zero, to within optimization tolerance.\
For the asymmetric initial guesses, a study shows that the visual change from symmetry to asymmetry occurs roughly at an asymmetry index of around $0.3$.\
This cutoff is used to classify designs between asymmetric and quasi-symmetric in \Cref{sec:asym_VC_Canti}.\ 

\section{Stress minimization}
\label{app:StressMinimization}

After computation of the displacement vector $\mathbf{u}$, every element $k$ has a corresponding stress tensor (in Voigt notation)  $\bm{\sigma}_k \in \mathbb{R}^{3}$ defined at the element's Gauss point, i.e., the element center.\
It is computed as
\begin{equation}
\bm{\sigma}_k = \mathbf{C}_{\mathrm{solid}} \mathbf{B}_k \mathbf{u},
\end{equation}
where $ \mathbf{C}_{\mathrm{solid}} \in \mathbb{R}^{3 \times 3}$ is the constitutive matrix of the solid material and $\mathbf{B}_k \in \mathbb{R}^{3 \times n}$ is the global strain-displacement matrix of element $k$.\
The Von Mises stress in element $k$ is then
\begin{equation}
\sigma_{k, \mathrm{VM}} = \sqrt{\bm{\sigma}_k^\top \mathbf{M} \bm{\sigma}_k},
\end{equation}
with $M$ the following constant matrix

\begin{equation}
\mathbf{M} = 
\begin{bmatrix}
1 & -\frac{1}{2} & 0\\
-\frac{1}{2} & 1 & 0\\
0 & 0 & 3
\end{bmatrix}.
\end{equation}

To avoid singularities at low densities, we employ the $\varepsilon$-relaxed approach~\cite{Cheng1997} and define the relaxed Von Mises stress as 

\begin{equation}
\sigma_{k, \mathrm{rVM}} = \sigma_{k, \mathrm{VM}} \frac{x_{k,P}}{\varepsilon(1-x_{k,P})+x_{k,P}},
\end{equation}
where $x_{k,P}$ is the physical density at element $k$ and $\varepsilon=0.2$ to ensure stress accuracy, following da Silva et al.~\cite{daSilva2019}.\
Note that the $\varepsilon$ symbol clashes with its usage in the main body of the text and is only temporarily redefined here to keep consistency with topology optimization literature.\

da Silva et al.'s value for $\varepsilon$ is only valid when employing SIMP interpolation with $P=3$.\
This is achieved by starting at $P=1$ and going to $P=3$ in $20$ steps, requiring convergence (relative objective change $< 1\mathrm{e}{-3}$ for five iterations) each time.\
The Heaviside projection parameter $\beta$ is increased to $2R/h$, with $h$ the element width, in the same number of steps as well.\
For the L-shaped problems, the filter radius is $R=0.1$ (i.e., $20$ element widths).\

The relaxed Von Mises stress is aggregated via the function of Kreisselmeier and Steinhauser~\cite{KreisselmeierSteinhauser1979} to produce the aggregated stress, which approximates the maximum stress in the design as follows.
\begin{equation}
s_{agg}=s(\mathbf{x})= \sigma^{\max}_{\mathrm{rVM}} + \frac{1}{\mu_{KS}} \ln \left(  \sum_{k=1}^n \exp( \mu_{KS} (\sigma_{k, \mathrm{rVM}}-\sigma^{\max}_{\mathrm{rVM}})) \right)
\end{equation}
Here, $\mu_{KS}$ is a hyperparameter set to $5e7$ and $\sigma^{\max}_{\mathrm{rVM}}$ is the maximum relaxed Von Mises stress.\

A relative convergence tolerance criterion of $1\mathrm{e}{-4}$ is used.\
All other parameters are the same as the compliance minimization examples.\

\section{Dynamic compliance minimization}
\label{app:DynamicComplianceMinimization}

The dynamic compliance is computed by solving the following frequency-dependent system of equations.
\begin{equation} \label{eq:dynamic_Ku=f}
\mathbf{A}(\mathbf{x}, \omega) \mathbf{u}_d=\mathbf{f}_{\mathrm{ext}}.
\end{equation}

The dynamic system matrix $\mathbf{A}(\mathbf{x}, \omega)$ depends on the angular frequency $\omega=2 \pi f$, with $f$ the driving frequency, and is defined as

\begin{equation}
\mathbf{A}(\mathbf{x}, \omega) = (1+i\zeta)\mathbf{K}(\mathbf{x}) - \omega^2 \mathbf{M}(\mathbf{x}),
\end{equation}

where $\zeta$ is the structural damping coefficient set at $0.1$ for the examples in this work, $i=\sqrt{-1}$ is the imaginary unit and $\mathbf{M}(\mathbf{x})$ is the mass matrix.\
The latter is constructed via the finite element method, using a material density of $1000 \frac{\mathrm{kg}}{\mathrm{m}^3}$ and SIMP interpolation with $\delta=1\mathrm{e}{-9}$ and penalization $P=3$.\
The stiffness matrix $\mathbf{K}(\mathbf{x})$ is constructed as before, now with Young modulus $5\mathrm{GPa}$ and SIMP interpolation with $\delta=5\mathrm{kPa}$ and a varying penalization.\
The minimum SIMP values for the density and stiffness are chosen in this way to prevent spurious modes~\cite{pedersen2000maximization}.\

For simplicity, $\mathbf{f}_{\mathrm{ext}}$ is the force vector used in the (static) compliance scaled by $100$, modeling a sinusoidally varying force of $100$ Newton applied at the tip of the cantilever.\
The dynamic compliance is defined as $|\mathbf{u}_d^\top \mathbf{f}_{\mathrm{ext}}|$.\
Here, $| \cdot |$ is the complex norm.\
As it can take values with strongly differing orders of magnitude, we define the function $d_{\omega}(\mathbf{x}, \omega)=10 \log_{10} \left( |\mathbf{u}_d(\omega)^\top \mathbf{f}_{\mathrm{ext}}| \right)$.\
Note that this value is negative for the cantilevers studied in this work, which simply means that the dynamic compliance lies below $1$.\
For numerical stability, a frequency range around the target frequency of $100$ Hz is considered: $d_{\omega}(\mathbf{x}, \omega)$ is computed at $99$Hz, $100$Hz and $101$Hz and then averaged.\
Thus, $d(\mathbf{x}) = \frac{1}{3} \left( d_{\omega}(\mathbf{x}, 2\pi99) + d_{\omega}(\mathbf{x}, 2\pi100) + d_{\omega}(\mathbf{x}, 2\pi101) \right)$.
During optimization, $d(\mathbf{x})$ is rescaled by dividing with $0.1 d_{\mathrm{solid}}$, where $d_{\mathrm{solid}}=d(\mathbf{1})$ is the performance of the solid design.\

Dynamic compliance minimization uses the same stopping criterion as before (\ref{app:ComplianceMinimization}), but with a relative objective function tolerance of $1\mathrm{e}{-5}$.\
The SIMP penalization parameter $P$ goes from $1$ to $3$ in $50$ steps of $\Delta P=0.04$.\
The $\beta$ projection parameter goes from $1$ to $2R$, with $R$ the filter radius of $10$ element widths, in the same number of steps.\
The double filtering technique of Christiansen et al.~\cite{christiansen2015doublefilt} is used to prevent the blueprint, dilated and eroded designs from having a different topology.\
It involves an additional presmoothing and projection step with the same $\beta$ value and a halved filter radius of $5$ element widths.\\

After the final $P$ and $\beta$ values are reached and convergence is achieved, the continuation proceeds by explicitly penalizing the measure of non-discreteness $M_{\mathrm{nd}}(\mathbf{x})=\frac{100}{n} \sum_{i=1}^{n} 4 x_i (1-x_i)$ of the eroded, dilated and blueprint designs.\
This is done by adding $\frac{\lambda}{3} \left(M_{\mathrm{nd}} ( \mathbf{x}_b)+M_{\mathrm{nd}}(\mathbf{x}_e)+M_{\mathrm{nd}}(\mathbf{x}_d) \right)$ to the objective and doubling $\lambda$, starting from $1\mathrm{e}{-10}$, until $\lambda>0.1$.\
The MMA parameters are identical except for the outer move limit, which is set to $0.1$ instead of $0.2$.

\section{Hypervolume indicator}
\label{app:HVI}

\mc{The hypervolume indicator is computed via the definition in~\cite{HVI}.\
It measures the area of the region that is dominated by the approximation and bounded from above by a point $\mathbf{b}$ (denoted $\mathbf{r}$ in~\cite{HVI}).\
In this work, the indicator values are always computed for the rescaled function values, so we can choose $\mathbf{b} = (1,1)$.\ }

\clearpage
\bibliography{mybibfile}

\end{document}